\documentclass{amsart}
\usepackage{amssymb}
\newtheorem{theorem}{Theorem}[section]
\newtheorem{corollary}[theorem]{Corollary}
\newtheorem{lemma}[theorem]{Lemma}
\newtheorem{proposition}[theorem]{Proposition}
\theoremstyle{definition}

\newtheorem{remark}[theorem]{Remark}

\theoremstyle{remark}
\newtheorem*{acknowledgements}{Acknowledgements}

\newlength{\displayboxwidth}

\newcommand{\bs}{\boldsymbol{\sigma}}
\newcommand{\wlim}{\operatorname*{w-lim}}
\newcommand{\wstarlim}{\operatorname*{w*-lim}}
\newcommand{\chapter}{\section}
\newcommand{\chaptername}{section}
\newcommand{\Chaptername}{Section}
\newcommand{\proofheadfont}{\textit}
\newcommand{\proofindent}{\noindent}
\newbox\ipbox
\newcommand{\ip}[2]{\left\langle #1\mathrel{\mathchoice
{\setbox\ipbox=\hbox{$\displaystyle \left\langle\mathstrut #1#2\right\rangle$}
\vrule height\ht\ipbox width0.25pt depth\dp\ipbox}
{\setbox\ipbox=\hbox{$\textstyle \left\langle\mathstrut #1#2\right\rangle$}
\vrule height\ht\ipbox width0.25pt depth\dp\ipbox}
{\setbox\ipbox=\hbox{$\scriptstyle \left\langle\mathstrut #1#2\right\rangle$}
\vrule height\ht\ipbox width0.25pt depth\dp\ipbox}
{\setbox\ipbox=\hbox{$\scriptscriptstyle \left\langle\mathstrut #1#2\right\rangle$}
\vrule height\ht\ipbox width0.25pt depth\dp\ipbox}
} #2\right\rangle}
\newcommand{\diracb}[1]{\left\langle #1\mathrel{\mathchoice
{\setbox\ipbox=\hbox{$\displaystyle \left\langle\mathstrut #1\right.$}
\vrule height\ht\ipbox width0.25pt depth\dp\ipbox}
{\setbox\ipbox=\hbox{$\textstyle \left\langle\mathstrut #1\right.$}
\vrule height\ht\ipbox width0.25pt depth\dp\ipbox}
{\setbox\ipbox=\hbox{$\scriptstyle \left\langle\mathstrut #1\right.$}
\vrule height\ht\ipbox width0.25pt depth\dp\ipbox}
{\setbox\ipbox=\hbox{$\scriptscriptstyle \left\langle\mathstrut #1\right.$}
\vrule height\ht\ipbox width0.25pt depth\dp\ipbox}
}\right. }
\newcommand{\dirack}[1]{\left. \mathrel{\mathchoice
{\setbox\ipbox=\hbox{$\displaystyle \left.\mathstrut #1\right\rangle$}
\vrule height\ht\ipbox width0.25pt depth\dp\ipbox}
{\setbox\ipbox=\hbox{$\textstyle \left.\mathstrut #1\right\rangle$}
\vrule height\ht\ipbox width0.25pt depth\dp\ipbox}
{\setbox\ipbox=\hbox{$\scriptstyle \left.\mathstrut #1\right\rangle$}
\vrule height\ht\ipbox width0.25pt depth\dp\ipbox}
{\setbox\ipbox=\hbox{$\scriptscriptstyle \left.\mathstrut #1\right\rangle$}
\vrule height\ht\ipbox width0.25pt depth\dp\ipbox}
} #1\right\rangle}
\def\openone
{\mathchoice
{\hbox{\upshape \small1\kern-3.3pt\normalsize1}}
{\hbox{\upshape \small1\kern-3.3pt\normalsize1}}
{\hbox{\upshape \tiny1\kern-2.3pt\SMALL1}}
{\hbox{\upshape \Tiny1\kern-2pt\tiny1}}}
\numberwithin{equation}{section}
\newlength{\qedskip}
\newlength{\qedadjust}

\begin{document}
\title[Pure states on $\mathcal{O}_{d}$]{Pure states on $\mathcal{O}_{d}$}
\author[O. Bratteli]{Ola~Bratteli}
\address{Department of Mathematics\\
University of Oslo\\
PB 1053 -- Blindern\\
N-0316 Oslo\\
Norway}
\email{bratteli@math.uio.no}
\author[P. E. T. Jorgensen]{Palle~E.~T.~ Jorgensen}
\address{Department of Mathematics\\
The University of Iowa\\
14 MacLean Hall\\
Iowa City, IA 52242-1419\\
U.S.A.}
\email{jorgen@math.uiowa.edu}
\author[A. Kishimoto]{Akitaka~Kishimoto}
\address{Department of Mathematics\\
University of Hokkaido\\
Sapporo 060\\
Japan}
\email{kishi@math.hokudai.ac.jp}
\author[R. F. Werner]{Reinhard~F.~Werner}
\address{Institut f\"{u}r Mathematische Physik\\
Universit\"{a}t Braunschweig\\
Mendelssohnstr.\ 3\\
D-38016 Braunschweig\\
Germany}
\email{R.Werner@tu-bs.de}
\thanks{Research supported by the Norwegian Research Council, the University of Oslo,
and the Scandinavia--Japan Sasakawa Foundation.}
\dedicatory{Dedicated to Professor Erling St\o{}rmer
on the occasion of his sixtieth birthday}
\begin{abstract}
We study representations of the Cuntz algebras $\mathcal{O}_{d}$ and their
associated decompositions.
In the case that these representations are irreducible, their restrictions to
the gauge-invariant subalgebra $\operatorname*{UHF}\nolimits_{d}$ have an
interesting cyclic structure.
If $S_{i}$, $1\leq i\leq d$, are representatives of the Cuntz relations on a
Hilbert space $\mathcal{H}$, special attention is given to the subspaces which
are invariant under $S_{i}^{\ast}$. The applications include wavelet
multiresolutions corresponding to wavelets of compact support (to appear in
the later paper \cite{BEJ97}), and finitely correlated states on
one-dimensional quantum spin chains.
\end{abstract}
\maketitle
\tableofcontents

\chapter{\label{Intro}Introduction}

The aim of the present paper was at the outset threefold:

\begin{enumerate}
\item  To develop further and simplify the theory of finitely (and infinitely)
correlated states of the Cuntz algebra $\mathcal{O}_{d}$ given in
\cite{BrJo97a}.

\item  To apply this theory to analyze in detail the representations of
$\mathcal{O}_{N}$ coming from compactly supported wavelets constructed by
multiresolution wavelet analysis of scale $N$ \cite{BrJo97b}. The main idea is
that repeated applications of the adjoints of the Cuntz operators on any
trigonometric polynomial in $L^{2}\left(  \mathbb{T}\right)  $ in that case
ultimately maps the polynomial into a fixed finite-dimensional subspace
$\mathcal{K}\subset L^{2}\left(  \mathbb{T}\right)  $ of low-order
polynomials, and thus the results of the present paper apply. This application
will be postponed to the paper \cite{BEJ97}.

\item  To understand better the connection between the theory of finitely
correlated states on one-dimensional quantum spin chains developed in
\cite{FNW92,FNW94} and the
corresponding states on the Cuntz algebras.
\end{enumerate}

The setting and results (especially Theorem \ref{Thm5.1}) also serve as a
generalization of the single-operator commutant lifting theorem \cite{DMP68}
from one variable to several. In this setting, $\mathcal{O}_{d}$, for $d\geq
2$, is viewed as the multivariable version of the familiar $C^{\ast}$-algebra
generated by a single isometry.

Recall that if $d\in\left\{  2,3,\dots\right\}  $, the \emph{Cuntz algebra}
$\mathcal{O}_{d}$ is the universal $C^{\ast}$-algebra generated by elements
$s_{1},\dots,s_{d}$ subject to the relations
\begin{align*}
s_{i}^{\ast}s_{j}^{{}}  &  =\delta_{ij}\openone,\\
\sum_{j=1}^{d}s_{j}^{{}}s_{j}^{\ast}  &  =\openone.
\end{align*}

There is a canonical action of the group $U\left(  d\right)  $ of unitary
$d\times d$ matrices on $\mathcal{O}_{d}$ given by%
\[
\tau_{g}\left(  s_{i}\right)  =\sum_{j=1}^{d}\overline{g_{ji}}s_{j}%
\]
for $g=\left[  g_{ij}\right]  _{i,j=1}^{d}\in U\left(  d\right)  $. In
particular the \emph{gauge action} is defined by%
\begin{equation}
\tau_{z}\left( s_{i}\right) 
=zs_{i},\quad z\in\mathbb{T}\subset\mathbb{C}%
\mkern2mu%
. \label{eqGaugeAction}
\end{equation}
If $\operatorname*{UHF}\nolimits_{d}$ is the fixed point subalgebra under the
gauge action, then $\operatorname*{UHF}\nolimits_{d}$ is the closure of the
linear span of all Wick ordered monomials of the form%
\[
s_{i_{1}}^{{}}\cdots s_{i_{k}}^{{}}s_{j_{k}}^{\ast}\cdots s_{j_{1}}^{\ast}.
\]
$\operatorname*{UHF}\nolimits_{d}$ is isomorphic to the $\operatorname*{UHF}$
algebra of Glimm type $d^{\infty}$:%
\[
\operatorname*{UHF}\nolimits_{d}\cong M_{d^{\infty}}=\bigotimes_{1}^{\infty
}M_{d}%
\]
in such a way that the isomorphism carries the Wick ordered monomial above
into the matrix element%
\[
e_{i_{1}j_{1}}^{\left(  1\right)  }\otimes e_{i_{2}j_{2}}^{\left(  2\right)
}\otimes\dots\otimes e_{i_{k}j_{k}}^{\left(  k\right)  }\otimes\openone
\otimes\openone\otimes\cdots.
\]
The restriction of $\tau_{g}$ to $\operatorname*{UHF}\nolimits_{d}$ is then
carried into the action%
\[
\operatorname{Ad}\left(  g\right)  \otimes\operatorname{Ad}\left(  g\right)
\otimes\cdots
\]
on $\bigotimes_{1}^{\infty}M_{d}$. We define the canonical endomorphism
$\lambda$ on $\operatorname*{UHF}\nolimits_{d}$ (or on $\mathcal{O}_{d}$) by%
\begin{equation}
\lambda\left(  x\right)  =\sum_{j=1}^{d}s_{j}^{{}}xs_{j}^{\ast}
\label{eqNew1.1}
\end{equation}
and the isomorphism carries $\lambda$ over into the one-sided shift%
\[
x_{1}\otimes x_{2}\otimes x_{3}\otimes\cdots\longrightarrow\openone
\otimes x_{1}\otimes x_{2}\otimes\cdots
\]
on $\bigotimes_{1}^{\infty}M_{d}$.
(See
\cite{Cun77,Eva,BEGJ}.)

If $s_{i}\mapsto S_{i}\in\mathcal{B}\left(  \mathcal{H}\right)  $ is a
representation of the Cuntz relations on a Hilbert space $\mathcal{H}$, we
will consider the situation that there is a closed subspace $\mathcal{K}%
\subset\mathcal{H}$ such that%
\[
S_{i}^{\ast}\mathcal{K}\subset\mathcal{K}%
\]
for $i\in\mathbb{Z}_{d}$, and $\mathcal{K}$ is cyclic for the representation.
Thus, if $P\colon\mathcal{H}\rightarrow\mathcal{K}$ is the orthogonal
projection onto $\mathcal{K}$, we have%
\[
PS_{i}^{\ast}P=S_{i}^{\ast}P.
\]
In this situation, define $V_{i}\in\mathcal{B}\left(  \mathcal{K}\right)  $ by%
\[
V_{i}=PS_{i}=PS_{i}P.
\]
Then%
\[
\sum_{i\in\mathbb{Z}_{d}}V_{i}^{{}}V_{i}^{\ast}=\openone
\]
so the map $\bs\colon\mathcal{B}\left(  \mathcal{K}\right)  \rightarrow
\mathcal{B}\left(  \mathcal{K}\right)  $ defined by%
\begin{equation}
\bs\left(  X\right)  =\sum_{i\in\mathbb{Z}_{d}}V_{i}^{{}}XV_{i}^{\ast}
\label{eqNew1.2}
\end{equation}
is completely positive and unital. We show that the representation can be
completely recovered from $\left(  \mathcal{K},V_{1},\dots,V_{d}\right)  $ in
Theorem \ref{Thm2.1} and Theorem \ref{Thm5.1}, and the commutant of the
representation is isometrically order isomorphic to the fixed point set
$\mathcal{B}\left(  \mathcal{K}\right)  ^{\bs}=\left\{  A\in\mathcal{B}\left(
\mathcal{K}\right)  \mid\bs\left(  A\right)  =A\right\}  $ by Proposition
\ref{Pro4.1} and Theorem \ref{Thm5.1}. This fixed point set is not an algebra
in general, as is discussed in some detail in \Chaptername{} \ref{Ergodic}. In
particular the representation of $\mathcal{O}_{d}$ is irreducible if and only
if $\bs$ is ergodic in the sense that $\mathcal{B}\left(  \mathcal{K}\right)
^{\bs}=\mathbb{C}%
\mkern2mu%
\openone$. In \Chaptername{} \ref{Restriction} we assume that the
representation is irreducible and study its restriction to
$\operatorname*{UHF}\nolimits_{d}$ in the case that there is a normal $\bs
$-invariant state $\varphi$ on $\mathcal{B}\left(  \mathcal{K}\right)  $. Such
a state is automatically unique if it exists, and if $\mathcal{K}$ is
finite-dimensional it always exists. In this case we replace $\mathcal{K}$
with the smaller $S_{i}^{\ast}$-invariant space $E\mathcal{H}$, where $E$ is
the support projection of $\varphi$, replace $\varphi$ with its restriction to
$E\mathcal{B}\left(  \mathcal{K}\right)  E=\mathcal{B}\left(  E\mathcal{K}%
\right)  $, and we define a 
state $\psi$ on
$\mathcal{O}_{d}$
by%
\[
\psi\left(  s_{i_{1}}^{{}}\cdots s_{i_{n}}^{{}}s_{j_{m}}^{\ast}\cdots
s_{j_{1}}^{\ast}\right)  =\varphi\left(  ES_{i_{1}}^{{}}\cdots S_{i_{n}}^{{}%
}S_{j_{m}}^{\ast}\cdots S_{j_{1}}^{\ast}E\right)  .
\]
Then 
$\psi\circ
\lambda=\psi$. 
(Occasionally, we will identify $\psi$ with its
normal extension to $\mathcal{B}\left(  \mathcal{H}\right)  $, defined by
$\psi \left( X\right) =\varphi \left( EXE\right) $
for $X\in \mathcal{B}\left(  \mathcal{H}\right) $.
This extension is a type $\mathrm{I}$ factor
state with multiplicity $\dim\left(  E\right)  $.)
We show in Theorem \ref{Thm6.3} that the set of $t\in
\mathbb{T}$ such that $\psi\circ\tau_{t}=\psi$ is equal to the peripheral
point spectrum $\operatorname{PSp}\left(  \bs\right)  \cap\mathbb{T}$ of $\bs
$, and this set is a finite subgroup of $\mathbb{T}$. If $k$ is the
order of this subgroup, and $U\in\mathcal{B}\left(  \mathcal{H}\right)  $ is
the unitary operator such that 
$\tau_{\frac{1}{k}}$, corresponding to
$z=e^{i\frac{2\pi}{k}}$ in (\ref{eqGaugeAction}),
satisfies
$\tau_{\frac{1}{k}}=\operatorname{Ad}\left(
U\right)  $ 
with
$U^{k}=\openone$ ($U$ is unique up to a phase factor in
$\mathbb{Z}_{k}\subset\mathbb{T}$), and%
\[
U=\sum_{l\in\mathbb{Z}_{k}}e^{\frac{2\pi il}{k}}E_{l}%
\]
is the spectral decomposition of $U$, then the subalgebra $\operatorname*{UHF}%
\nolimits_{d}\subset\mathcal{O}_{d}$ acts irreducibly on each of the subspaces
$E_{l}\mathcal{H}$, the corresponding representations of $\operatorname*{UHF}%
\nolimits_{d}$ are irreducible and mutually disjoint, and are mapped
cyclically into each other by the endomorphism $\lambda$.

In particular, this means that the restriction of the representation to
$\operatorname*{UHF}\nolimits_{d}$ is irreducible if and only if the
peripheral point spectrum $\operatorname{PSp}\left(  \bs\right)
\cap\mathbb{T}$ of $\bs$ consists of the point $1$ alone. It is remarkable
that, if $\mathcal{K}$ is finite-dimensional, this is exactly the condition
ensuring that the translation-invariant state defined by $\left\{
\varphi,V_{1},\dots,V_{d}\right\}  $ on the two-sided one-dimensional quantum
chain $\bigotimes_{-\infty}^{\infty}M_{d}=\bigotimes_{\mathbb{Z}}^{{}}M_{d}$
is pure \cite{FNW92,FNW94}.
To be precise, this condition on $\left\{  \varphi,V_{1},\dots,V_{d}\right\}
$ is sufficient to ensure purity of $\omega$. It is not necessary for the
given $\left\{  \varphi,V_{1},\dots,V_{d}\right\}  $, but if $\omega$ is pure
and finitely generated, there exists 
some
$\left\{  \varphi,V_{1},\dots
,V_{d}\right\}  $ on a finite-dimensional $\mathcal{K}$, defining $\omega$, such
that the corresponding $\bs$ is ergodic and has trivial peripheral spectrum.
One source of the nonuniqueness of $\left\{  \varphi,V_{1},\dots
,V_{d}\right\}  $, and the corresponding non-necessity of the conditions on
this set, is the following: If $\mathcal{K}$ is replaced by $\mathcal{K}%
\otimes\mathcal{K}^{\prime}$, where $\mathcal{K}^{\prime}$ is a Hilbert space
of finite dimension $\geq2$, $V_{k}$ by $V_{k}\otimes\operatorname{id}$ and
$\varphi$ by $\varphi\otimes\varphi^{\prime}$, where $\varphi^{\prime}$ is a
faithful state on $\mathcal{B}\left(  \mathcal{K}^{\prime}\right)  $, then the
new data define exactly the same state as the old, but the fixed point set of
the new $\bs$ contains at least $\openone\otimes\mathcal{B}\left(
\mathcal{K}^{\prime}\right)  $. To avoid this kind of degeneracy, we make in
\Chaptername{} \ref{Quantum} the overall assumption that the operators $V_{1}%
,\dots,V_{d}$ on $\mathcal{K}$ (which does not need to be finite-dimensional)
generate a factor $\mathcal{M}$ with a faithful normal $\bs$-invariant state
$\varphi$, and that $\mathcal{B}\left(  \mathcal{K}\right)  ^{\bs}%
=\mathcal{M}^{\prime}$.
If in addition $\mathcal{M}$ is type $\mathrm{I}$, we prove that the
corresponding
translationally invariant state $\omega$ on $\bigotimes_{\mathbb{Z}}M_{d}$ is
pure if and only if $\operatorname{PSp}\left(  \bs|_{\mathcal{M}}\right)
\cap\mathbb{T}=\left\{  1\right\}  $. If $\mathcal{M}$ is a finite type
$\mathrm{I}$ factor, this is exactly the result in \cite{FNW94}. If
$\mathcal{M}$ is not type $\mathrm{I}$, this equivalence is no longer true,
but in that case we can prove that $\omega$ is pure if and only if $\omega$ is
a factor state, i.e., if and only if $\omega$ has the clustering property
$\lim_{\left|  n\right|  \rightarrow\infty}\omega\left(  x\lambda^{n}\left(
y\right)  \right)  =\omega\left(  x\right)  \omega\left(  y\right)  $ for each
pair $x,y\in\bigotimes_{\mathbb{Z}}M_{d}$.

For more background material on the representations
of  $\mathcal{O}_{d}$, see, e.g., \cite{BJP96,BrJo96b,BrJo97a}.
A representation $\left( S_{i}\right) $ of
$\mathcal{O}_{d}$ on $\mathcal{H}$ defines 
an endomorphism 
$\sigma \left( \,\cdot \,\right) =\sum _{i}S_{i}^{{}}\cdot S_{i}^{\ast}$
of $\mathcal{B}\left(  \mathcal{H}\right)  $,
and conversely. Moreover,
the connection between an
endomorphism $\sigma $,
corresponding to $\lambda $ in (\ref{eqNew1.1}),
and the associated completely
positive map $\bs $ in (\ref{eqNew1.2}) above, is given
by
\[
P\sigma \left( X\right) P=\bs \left( PXP\right) ,\quad X\in 
\mathcal{B}\left(  \mathcal{H}\right) .
\]
The lifting problem, addressed in 
\Chaptername{} \ref{General} below,
then concerns the 
reconstruction of the endomorphism $\sigma $,
or 
the associated $\mathcal{O}_{d}$-representation,
from some given completely positive
normal unital map $\bs $ of 
$\mathcal{B}\left(  \mathcal{K}\right) $.
 
Other somewhat related aspects of the representation theory of $\mathcal{O}%
_{d}$, and its restriction to $\operatorname*{UHF}\nolimits_{d}$, have 
been considered in \cite{LTW88,Spi90,FoLa97,Fow97,DaPi97b}.

\chapter{\label{General}General states on $\mathcal{O}_{d}$}

\setlength{\displayboxwidth}{\textwidth} \addtolength{\displayboxwidth
}{-1.25\leftmargini}
First some notation: Let $d\in\left\{  2,3,\dots\right\}  $ and let
$\mathbb{Z}_{d}$ be a set of $d$ elements. (The group structure of
$\mathbb{Z}_{d}$ is spurious for the purposes of this paper.) Let
$\mathcal{I}$ be the set of finite sequences $\left(  i_{1},\dots
,i_{m}\right)  $ where $i_{k}\in\mathbb{Z}_{d}$ and $m\in\left\{
1,2,\dots\right\}  $. We also include the empty sequence $\varnothing$ in
$\mathcal{I}$, and denote elements in $\mathcal{I}$ by $I,J,\dots$. If
$I=\left(  i_{1},\dots,i_{m}\right)  \in\mathcal{I}$ and $i\in\mathbb{Z}_{d}$,
we let $Ii$ denote the element $\left(  i_{1},\dots,i_{m},i\right)  $ in
$\mathcal{I}$, and $s_{I}=s_{i_{1}}s_{i_{2}}\cdots s_{i_{m}}\in\mathcal{O}%
_{d}$ and $s_{I}^{\ast}=s_{i_{m}}^{\ast}s_{i_{m-1}}^{\ast}\cdots s_{i_{1}%
}^{\ast}\in\mathcal{O}_{d}$. In particular $s_{\varnothing}^{{}}%
=\openone=s_{\varnothing}^{\ast}$.

The following theorem is a version of a result of Popescu \cite{Pop89}. It
generalizes \cite{BEGJ}. We give a streamlined proof which applies in this case.

\begin{theorem}
\label{Thm2.1}\raggedright Let $d\in\left\{  2,3,\dots\right\}  $. There is a
canonical one-one correspondence between the following objects.
\begin{equation}
\begin{minipage}[t]{\displayboxwidth}\raggedright States $\hat{\omega
}$ on $\mathcal{O}_{d}$.\end{minipage} \label{Thm2.1(1)}%
\end{equation}
\begin{equation}
\begin{minipage}[t]{\displayboxwidth}\raggedright Functions $C\colon
\mathcal{I}\times\mathcal{I}\rightarrow\mathbb{C}%
$ with the following properties: \end{minipage} \label{Thm2.1(2)}%
\end{equation}
\begin{enumerate}
\item \label{Thm2.1(2)(1)}$C\left(  \varnothing,\varnothing\right)  =1$,

\item \label{Thm2.1(2)(2)}for any function $\lambda\colon\mathcal{I}%
\rightarrow\mathbb{C}$ with finite support we have $\sum_{I,J\in\mathcal{I}%
}\overline{\lambda\left(  I\right)  }C\left(  I,J\right)  \lambda\left(
J\right)  \geq0$,

\item \label{Thm2.1(2)(3)}$\sum_{i\in\mathbb{Z}_{d}}C\left(  Ii,Ji\right)
=C\left(  I,J\right)  $ for all $I,J\in\mathcal{I}$.
\end{enumerate}
\begin{equation}
\begin{minipage}[t]{\displayboxwidth}\raggedright
Unitary equivalence classes of objects $\left( \mathcal{K},\Omega,V_{1}%
,\dots,V_{d}\right) $ where \end{minipage} \label{Thm2.1(3)}%
\end{equation}
\begin{enumerate}
\item \label{Thm2.1(3)(1)}$\mathcal{K}$ is a Hilbert space,

\item \label{Thm2.1(3)(2)}$\Omega$ is a unit vector in $\mathcal{K}$,

\item \label{Thm2.1(3)(3)}$V_{1},\dots,V_{d}\in\mathcal{B}\left(
\mathcal{K}\right)  $,

\item \label{Thm2.1(3)(4)}the linear span of vectors of the form $V_{I}^{\ast
}\Omega$, where $I\in\mathcal{I}$, is dense in $\mathcal{K}$,

\item \label{Thm2.1(3)(5)}$\sum_{i\in\mathbb{Z}_{d}}V_{i}^{{}}V_{i}^{\ast
}=\openone_{\mathcal{K}}$.
\end{enumerate}
The correspondence is given by%
\begin{equation}
\hat{\omega}\left(  s_{I}^{{}}s_{J}^{\ast}\right)  =C\left(  I,J\right)
=\ip{V_{I}^{\ast}\Omega}{ V_{J}^{\ast}\Omega}. \label{eq2.4}%
\end{equation}
\end{theorem}

\begin{proof}
It is immediate that if either $\hat{\omega}$ or $\left(  \mathcal{K}%
,\Omega,V_{1},\dots,V_{d}\right)  $ is given, and $C\left(  \,\cdot
\,,\,\cdot\,\right)  $ is defined by the relation (\ref{eq2.4}), then $C$
satisfies (\ref{Thm2.1(2)}). ((\ref{Thm2.1(2)(1)}) corresponds to the
normalization $\left\|  \hat{\omega}\right\|  =1=\hat{\omega}\left(
\openone\right)  $, or $\left\|  \Omega\right\|  =1$, (\ref{Thm2.1(2)(2)})
corresponds to positivity, and (\ref{Thm2.1(2)(3)}) to the relations $\sum
_{i}s_{i}^{{}}s_{i}^{\ast}=\openone$, $\sum_{i}V_{i}^{{}}V_{i}^{\ast}%
=\openone$.)

To go from the positive definite function $C$ in (\ref{Thm2.1(2)}) to the
object $\left(  \mathcal{K},\Omega,V_{1},\dots,V_{d}\right)  $ one uses the
usual Kolmogorov construction: one puts $\mathcal{K}$ equal to the completion
of the free vector space $\mathcal{L}\left(  \mathcal{I}\right)  $ of all
formal finite linear combinations $\sum_{I\in\mathcal{I}}\lambda\left(
I\right)  I$ (alias all functions $\lambda\colon\mathcal{I}\rightarrow
\mathbb{C}$ with finite support) with respect to the pre-inner product defined
by sesquilinearity from%
\[
\ip{I}{J}=C\left(  I,J\right)  ,
\]
after dividing out the vectors of zero norm. This gives a map $\Phi
\colon\mathcal{L}\left(  \mathcal{I}\right)  \rightarrow\mathcal{K}$, and one
defines $V_{i}$ by%
\[
V_{i}^{\ast}\Phi\left(  I\right)  =\Phi\left(  Ii\right)  .
\]
It is now routine to check the properties (\ref{Thm2.1(3)(1)}%
)--(\ref{Thm2.1(3)(5)}) in (\ref{Thm2.1(3)}).

To go from the object $\left(  \mathcal{K},\Omega,V_{1},\dots,V_{d}\right)  $
in (\ref{Thm2.1(3)}) to the state $\hat{\omega}$ on $\mathcal{O}_{d}$, we will
actually prove more:

(There is also a simple direct way of establishing this direction which will
be exhibited in Remark \ref{RemNew5.2}.)

\begin{lemma}
\label{Lem2.2}\cite{Pop89} Assume that $\mathcal{K},\Omega,V_{1},\dots,V_{d}$
satisfy the properties \textup{(\ref{Thm2.1(3)(1)})--(\ref{Thm2.1(3)(5)})}
under \textup{(\ref{Thm2.1(3)})}. It follows that there exists a unique linear
map $R\colon\mathcal{O}_{d}\rightarrow\mathcal{B}\left(  \mathcal{K}\right)  $
such that%
\[
R\left(  s_{I}^{{}}s_{J}^{\ast}\right)  =V_{I}^{{}}V_{J}^{\ast}%
\]
and this map is completely positive.
\end{lemma}

\begin{proof}
Let $\mathcal{T}_{d}$ be the Cuntz--Toeplitz algebra, realized on the
unrestricted Fock space $\mathcal{\hat{H}}=\bigoplus_{k=0}^{\infty}\left(
\mathbb{C}^{d}\right)  ^{\otimes\,k}$ in the usual way%
\[
L_{i}\colon\xi\longmapsto\dirack{i}\otimes\xi,
\]
where $L_{i}$, $i=1,\dots,d$, are the operators mapping into $s_{i}$ after
dividing out by the compact operators \cite{Eva,BEGJ}. Let $\lambda
\in\mathbb{C}$, $\left|  \lambda\right|  <1$, and define an operator%
\[
W_{\lambda}\colon\mathcal{K}\longrightarrow\mathcal{\hat{H}}\otimes\mathcal{K}%
\]
by%
\[
W_{\lambda}\varphi=\sqrt{1-\left|  \lambda\right|  ^{2}}\bigoplus
_{k=0}^{\infty}\lambda^{k}\sum_{I\in\mathcal{I}_{k}}\dirack{I}\otimes
V_{I}^{\ast}\varphi
\]
where $\mathcal{I}_{k}$ 
denotes 
all sequences $I=\left(  i_{1},\dots,i_{k}\right)
$ of length $k$ with $i_{j}\in\mathbb{Z}_{d}$, and%
\[
\dirack{I}=\dirack{i_{1}}\otimes\dirack{i_{2}}\otimes\dots\otimes\dirack
{i_{k}}.
\]
One checks that $W_{\lambda}$ is an isometry, and
\[
\left(  L_{i}^{\ast}\otimes\openone_{\mathcal{K}}^{{}}\right)  W_{\lambda}%
^{{}}=\lambda W_{\lambda}^{{}}V_{i}^{\ast}.
\]
{}From this intertwining 
relation, 
and its adjoint, it follows that%
\begin{align*}
R_{\lambda}^{{}}\left(  L_{I}^{{}}L_{J}^{\ast}\right)   &  \equiv W_{\lambda
}^{\ast}\left(  L_{I}^{{}}L_{J}^{\ast}\otimes\openone_{\mathcal{K}}^{{}%
}\right)  W_{\lambda}^{{}}\\
&  =\bar{\lambda}_{{}}^{n}\lambda_{{}}^{m}V_{I}^{{}}V_{J}^{\ast}%
\end{align*}
if $I\in\mathcal{I}_{n}$, $J\in\mathcal{I}_{m}$. It follows from this explicit
Stinespring representation that the linear map defined from%
\[
L_{I}^{{}}L_{J}^{\ast}\longmapsto\bar{\lambda}_{{}}^{n}\lambda_{{}}%
^{m}V_{I}^{{}}V_{J}^{\ast}%
\]
is 
then 
completely positive for all $\left|  \lambda\right|  <1$, and taking the
limit as $\lambda\rightarrow1$, it follows 
further
that $R$ is completely positive as
a map from $\mathcal{T}_{n}$ into $\mathcal{B}\left(  \mathcal{K}\right)  $.
To check that $R$, thus defined, defines a map from $\mathcal{O}_{d}$ into
$\mathcal{B}\left(  \mathcal{K}\right)  $, we have to show that $R$
annihilates the ideal generated by the one-dimensional projection
$p=\openone-\sum_{i}L_{i}^{{}}L_{i}^{\ast}$, i.e., that $R\left(  XpY\right)
=0$ for all polynomials $X,Y$ in the $L_{i}^{{}}$'s and the $L_{i}^{\ast}$'s.
We may take $X,Y$ to be Wick ordered monomials, i.e., of the form $L_{I}^{{}%
}L_{J}^{\ast}$. Since $pL_{i}=0$, we may assume that $Y$ contains no factor
$L_{I}$, and by the same token we may assume that $X$ contains no factor
$L_{J}^{\ast}$, and hence $XpY$ has the form%
\[
L_{I}^{{}}pL_{J}^{\ast}=L_{I}^{{}}L_{J}^{\ast}-\sum_{i\in\mathbb{Z}_{n}}%
L_{Ii}^{{}}L_{Ji}^{\ast}.
\]
Using the definition of $R_{\lambda}$, and the relation $\sum_{i=1}^{d}%
V_{i}^{{}}V_{i}^{\ast}=
\openone 
$, it now follows that $R_{\lambda}\left(  XpY\right)
=0$. Hence $R\left(  XpY\right)  =0$, and $R$ defines a completely positive map
from $\mathcal{O}_{d}$ into $\mathcal{B}\left(  \mathcal{K}\right)  $. This
ends the proof of Lemma \ref{Lem2.2}.
\end{proof}

\proofindent\proofheadfont{Proof of Theorem} \ref{Thm2.1}, 
\proofheadfont{continued:} To go from the object
$\left(  \mathcal{K},\Omega,V_{1},\dots,V_{d}\right)  $ to the state
$\hat{\omega}$ is now clear: put%
\[
\hat{\omega}\left(  X\right)  =\ip{\Omega}{R\left(  X \right) \Omega}%
\]
where $R\colon\mathcal{O}_{d}\rightarrow\mathcal{B}\left(  \mathcal{K}\right)
$ is the completely positive linear map defined in Lemma \ref{Lem2.2}. Then%
\[
\hat{\omega}\left(  s_{I}^{{}}s_{J}^{\ast}\right)  =\ip{\Omega}{V_{I}^{{}%
}V_{J}^{\ast}\Omega}%
\]
so (\ref{eq2.4}) is fulfilled.

This establishes the one-one correspondence stated in Theorem \ref{Thm2.1}.
Of course, the system $\left(  \mathcal{K},\Omega,V_{1},\dots,V_{d}\right)  $
is not unique, but determined only up to unitary equivalence. The argument for
why this is so is exactly the same as the standard argument from
representation theory \cite[Theorem 2.3.16]{BrRoI} to the effect that a state
on a $C^{\ast}$-algebra only determines 
a cyclic representation up to unitary equivalence.
\end{proof}

\begin{remark}
\label{Rem2.3}Note also that there is a simple direct way of going from the
state $\hat{\omega}$ in (\ref{Thm2.1(1)}) to the object $\left(
\mathcal{K},\Omega,V_{1},\dots,V_{d}\right)  $ in (\ref{Thm2.1(3)}). If
$\left(  \mathcal{H},\Omega,\pi\right)  $ is the cyclic representation of
$\mathcal{O}_{d}$ defined by $\Omega$, let $\mathcal{K}$ be the closure of the
linear span of all vectors $S_{I}^{\ast}\Omega$, where $S_{I}=\pi\left(
s_{I}\right)  $. Let $P$ be the projection from $\mathcal{H}$ onto
$\mathcal{K}$, and put%
\[
V_{i}^{\ast}=PS_{i}^{\ast}P=S_{i}^{\ast}P.
\]
The property $\sum_{i}V_{i}^{{}}V_{i}^{\ast}=P=\openone_{\mathcal{K}}$ follows
immediately from $\sum_{i}S_{i}^{{}}S_{i}^{\ast}=\openone_{\mathcal{H}}$.
\end{remark}

One can use Lemma \ref{Lem2.2}
to prove stronger versions of Popescu's dilation theorem:

\begin{corollary}
\label{Cor2.4}Let $\mathcal{K}$ be a Hilbert space, and $D\in\mathcal{B}%
\left(  \mathcal{K}\right)  $ a positive operator, and $V_{1},\dots,V_{d}%
\in\mathcal{B}\left(  \mathcal{K}\right)  $ operators such that%
\[
\sum_{i}V_{i}^{{}}DV_{i}^{\ast}=D.
\]
Then there exists a unique continuous linear map $R\colon\mathcal{O}%
_{d}\rightarrow\mathcal{B}\left(  \mathcal{K}\right)  $ such that%
\[
R\left(  s_{I}^{{}}s_{J}^{\ast}\right)  =V_{I}^{{}}DV_{J}^{\ast}%
\]
and this map is completely positive.
\end{corollary}

\begin{proof}
Roughly, if $R^{\prime}$ is the completely positive map defined in Lemma
\ref{Lem2.2} from the operators $V_{i}^{\prime}=D^{-\frac{1}{2}}V_{i}^{{}%
}D^{\frac{1}{2}}$, one verifies that%
\[
R^{\prime}\left(  s_{I}^{{}}s_{J}^{\ast}\right)  =V_{I}^{\prime}V_{J}%
^{\prime\ast}=D^{-\frac{1}{2}}V_{I}^{{}}DV_{J}^{\ast}D^{-\frac{1}{2}}.
\]
Putting%
\[
R\left(  \,\cdot\,\right)  =D^{\frac{1}{2}}R^{\prime}\left(  \,\cdot\,\right)
D^{\frac{1}{2}},
\]
we obtain Corollary \ref{Cor2.4}. A more careful argument is given in Remark
\ref{RemNew5.3}.
\end{proof}

\chapter{\label{Ergodic}Ergodic theory of completely positive maps on
$\mathcal{B}\left(  \mathcal{H}\right)  $}

In this \chaptername{} we prove 
some more or less known
results about completely 
positive unital normal maps $\varphi $
of  $\mathcal{B}\left( \mathcal{K}\right)  $,
and we analyze 
the
fixed-point set
\[
\mathcal{B}\left( \mathcal{K}\right)  ^{\varphi} := 
\left\{ X\in \mathcal{B}\left( \mathcal{K}\right) \mid
\varphi  \left( X\right) =X\right\} .
\]
We will need the arguments from the proofs here later in the paper.

Let $\mathcal{K}$ be a Hilbert space, and $\varphi\colon\mathcal{B}\left(
\mathcal{K}\right)  \rightarrow\mathcal{B}\left(  \mathcal{K}\right)  $ a
normal unital completely positive map. Then there exists a family of operators
$V_{i}\in\mathcal{B}\left(  \mathcal{K}\right)  $ such that%
\[
\sum_{i}V_{i}^{{}}V_{i}^{\ast}=\openone
\]
and%
\[
\varphi\left(  X\right)  =\sum_{i}V_{i}^{{}}X V_{i}^{\ast}%
\]
for all $X\in\mathcal{B}\left(  \mathcal{K}\right)  $, where the sum converges
in weak operator topology \cite{EvLe77}.

\begin{lemma}
\label{Lem3.1}Let $p$ be a projection in $\mathcal{B}\left(  \mathcal{K}%
\right)  $. Then the following conditions are equivalent.

\begin{enumerate}
\item \label{Lem3.1(1)}There is a $\lambda\geq0$ such that $\varphi\left(
p\right)  \leq\lambda p$.

\item \label{Lem3.1(2)}$V_{i}p=pV_{i}p$ for all $i$.

\item \label{Lem3.1(3)}$\varphi\left(  p\right)  \leq p$.
\end{enumerate}
\end{lemma}

\begin{remark}
\label{Rem3.1}Condition \textup{(\ref{Lem3.1(1)})} is of course equivalent to
the condition that the weakly closed hereditary subalgebras $p\mathcal{B}%
\left(  \mathcal{K}\right)  p$ of $\mathcal{B}\left(  \mathcal{K}\right)  $
are invariant under $\varphi$. The property that there are no nontrivial
weakly closed hereditary subalgebras of $\mathcal{B}\left(  \mathcal{H}%
\right)  $ invariant under $\varphi$ is called \emph{irreducibility} of
$\varphi$ in \cite{Far96}, and since any such subalgebra is of the form
$p\mathcal{B}\left(  \mathcal{H}\right)  p$, irreducibility of $\varphi$ is
equivalent to the nonexistence of projections $p$ with the property
\textup{(\ref{Lem3.1(1)})} or \textup{(\ref{Lem3.1(2)})}. The proof of Lemma
\ref{Lem3.1} is extracted from \cite{Far96}.
\end{remark}

\begin{proof}
\textup{(\ref{Lem3.1(1)})} $\Rightarrow$ \textup{(\ref{Lem3.1(2)})} Assume
$\varphi\left(  p\right)  \leq\lambda p$. Then
\begin{align*}
0  &  \leq\left(  \openone-p\right)  \varphi\left(  p\right)  \left(
\openone-p\right)  \leq\left(  \openone-p\right)  \lambda p\left(
\openone-p\right)  =0,\\%
\intertext{so}%
0  &  =\sum_{i}\left(  \openone-p\right)  V_{i}^{{}}pV_{i}^{\ast}\left(
\openone-p\right)  =\sum_{i}\left(  \left(  \openone-p\right)  V_{i}p\right)
\left(  \left(  \openone-p\right)  V_{i}p\right)  ^{\ast},\\%
\intertext{and hence}%
0  &  =\left(  \openone-p\right)  V_{i}p,
\end{align*}
which is \textup{(\ref{Lem3.1(2)})}.

\textup{(\ref{Lem3.1(2)})} $\Rightarrow$ \textup{(\ref{Lem3.1(3)})} Assume
that $V_{i}p=pV_{i}p$ for all $i$. Then%
\[
\varphi\left(  p\right)  =\sum_{i}V_{i}^{{}}pV_{i}^{\ast}=\sum_{i}pV_{i}^{{}%
}pV_{i}^{\ast}p=p\varphi\left(  p\right)  p\leq\left\|  \varphi\left(
p\right)  \right\|  p\leq p.
\]

\textup{(\ref{Lem3.1(3)})} $\Rightarrow$ \textup{(\ref{Lem3.1(1)})} is trivial.
\end{proof}

\begin{lemma}
\label{Lem3.3}Let $p$ be a projection in $\mathcal{B}\left(  \mathcal{K}%
\right)  $. Then the following conditions are equivalent.

\begin{enumerate}
\item \label{Lem3.3(1)}$\varphi\left(  p\right)  =p$.

\item \label{Lem3.3(2)}$p\in\left\{  V_{i}^{{}},V_{i}^{\ast}\right\}
_{i}^{\prime}$, i.e., $pV_{i}=V_{i}p$ for all $i$.
\end{enumerate}
\end{lemma}

\begin{proof}
Since $\sum_{i}V_{i}^{{}}V_{i}^{\ast}=\openone$, \textup{(\ref{Lem3.3(2)})}
$\Rightarrow$ \textup{(\ref{Lem3.3(1)})} is trivial. Conversely, assume%
\[
\varphi\left(  p\right)  =p.
\]
Applying Lemma \ref{Lem3.1}, \textup{(\ref{Lem3.1(1)})} $\Rightarrow$
\textup{(\ref{Lem3.1(2)})}, on $p$ and $\openone-p$, we obtain%
\begin{align*}
V_{i}p  &  =pV_{i}p\\%
\intertext{and}%
V_{i}\left(  \openone-p\right)   &  =\left(  \openone-p\right)  V_{i}\left(
\openone-p\right)  ,\\%
\intertext{i.e.,}%
pV_{i}  &  =pV_{i}p,
\end{align*}
so (\ref{Lem3.3(2)}) holds.
\end{proof}

If $\mathcal{B}\left(  \mathcal{K}\right)  ^{\varphi}=\left\{  X
\in\mathcal{B}\left(  \mathcal{K}\right)  \mid\varphi\left(  X\right)
=X\right\}  $ were an algebra (which necessarily is weakly closed and closed
under involution), it would follow from Lemma \ref{Lem3.3} that
\[
\mathcal{B}\left(  \mathcal{K}\right)  ^{\varphi}=\left\{  V_{i}^{{}}%
,V_{i}^{\ast}\right\}  ^{\prime}%
\]
(the inclusion $\supset$ is trivial, as mentioned before). There is one
important special case where $\mathcal{B}\left(  \mathcal{K}\right)
^{\varphi}$ is an algebra, 
namely when there is a faithful $\varphi$-invariant state:

\begin{lemma}
\label{Lem3.4}\cite{FNW94} Assume that there is a faithful $\varphi$-invariant
state $\omega$ on $\mathcal{B}\left(  \mathcal{K}\right)  $. Then
$\mathcal{B}\left(  \mathcal{K}\right)  ^{\varphi}$ is an algebra, and hence%
\[
\mathcal{B}\left(  \mathcal{K}\right)  ^{\varphi}=\left\{  V_{i}^{{}}%
,V_{i}^{\ast}\right\}  ^{\prime}.
\]
\end{lemma}

\begin{proof}
We follow \cite[proof of Proposition 2.2]{FNW94}. By \cite[Theorem
3.1]{Choi74}, if $\varphi$ is any $2$-positive map on a $C^{\ast}$-algebra
$\mathfrak{A}$, then%
\[
\left\{  a\in\mathfrak{A}\mid\varphi\left(  a^{\ast}a\right)  =\varphi\left(
a\right)  ^{\ast}\varphi\left(  a\right)  \right\}  =\left\{  a\in\mathfrak
{A}\mid\varphi\left(  xa\right)  =\varphi\left(  x\right)  \varphi\left(
a\right)  \text{ for all }x\in\mathfrak{A}\right\}  .
\]
Going back to our case, assume that $\varphi\left(  a\right)  =a$ for some
$a\in\mathcal{B}\left(  \mathcal{K}\right)  $. By the generalized Schwarz
inequality, we then have%
\[
\varphi\left(  a^{\ast}a\right)  -a^{\ast}a=\varphi\left(  a^{\ast}a\right)
-\varphi\left(  a^{\ast}\right)  \varphi\left(  a\right)  \geq0.
\]
But%
\[
\omega\left(  \varphi\left(  a^{\ast}a\right)  -a^{\ast}a\right)  =0
\]
by invariance of $\omega$, and as $\omega$ is faithful it follows that%
\[
\varphi\left(  a^{\ast}a\right)  =a^{\ast}a=\varphi\left(  a^{\ast}\right)
\varphi\left(  a\right)  .
\]
By Choi's theorem,%
\[
\varphi\left(  xa\right)  =\varphi\left(  x\right)  a
\]
for all $x\in\mathcal{B}\left(  \mathcal{H}\right)  $, and if in particular
$x\in\mathcal{B}\left(  \mathcal{H}\right)  ^{\varphi}$, then%
\[
\varphi\left(  xa\right)  =xa.
\]
Thus $xa\in\mathcal{B}\left(  \mathcal{H}\right)  ^{\varphi}$, and
$\mathcal{B}\left(  \mathcal{H}\right)  ^{\varphi}$ is an algebra.
\end{proof}

Note that the map $\omega\mapsto\omega\circ\varphi$ is obviously a
continuous map on the state space of $\mathcal{B}\left(  \mathcal{H}\right)
$, and this space is compact in the weak$^{\ast}$- topology from
$\mathcal{B}\left(  \mathcal{H}\right)  $. Hence it follows from the
Schauder--Tychonoff fixed point theorem that there exists a state $\omega$
such that $\omega\circ\varphi=\omega$ (\cite{T35c}, \cite[p.\ 456, \S V.10.5,
Theorem 5]{DS1}). Unfortunately the state $\omega$ is not necessarily faithful.
For example: let $\mathcal{K}=\mathbb{C}^{n}$, $e_{ij}$ a full set of matrix
units for $\mathcal{B}\left(  \mathbb{C}^{n}\right)  =M_{n}$, and put
$V_{i}=e_{i1}$ for $i=1,\dots,n$. Then 
\begin{align*}
\sum_{i}V_{i}^{{}}V_{i}^{\ast
} &=\openone ,  \\
\left\{  V_{i}^{{}},V_{i}^{\ast}\right\}  ^{\prime}%
&=\mathbb{C}%
\mkern2mu%
\openone ,
\end{align*}
but the unique invariant state for $\varphi$ is the pure
state $\omega\left(  \sum_{ij}X_{ij}e_{ij}\right)  =X_{11}$.

The states fixed by 
$\varphi$ need not 
in general
be normal either. We will discuss these
states further in the beginning of \Chaptername{} \ref{Restriction}.

Actually, there also exist examples where $\mathcal{B}\left(  \mathcal{K}%
\right)  ^{\varphi}$ is not an algebra. The following example is from
\cite{Arv69,Arv72}: $\mathcal{K}=\mathbb{C}^{3}$,%
\[
\varphi%
\begin{pmatrix}
X_{11} & X_{12} & X_{13}\\
X_{21} & X_{22} & X_{23}\\
X_{31} & X_{32} & X_{33}%
\end{pmatrix}
=%
\begin{pmatrix}
X_{11} & 0 & 0\\
0 & X_{22} & 0\\
0 & 0 & \frac{1}{2}\left(  X_{11}+X_{22}\right)
\end{pmatrix}
.
\]
Then $\varphi$ is completely positive, and one checks that%
\[
\mathcal{B}\left(  \mathcal{K}\right)  ^{\varphi}=\left\{
\begin{pmatrix}
a & 0 & 0\\
0 & b & 0\\
0 & 0 & \frac{1}{2}\left(  a+b\right)
\end{pmatrix}
\biggm|a,b\in\mathbb{C}\right\}  ,
\]
which contains no nontrivial subalgebras, and hence $\left\{  V_{i}^{{}}%
,V_{i}^{\ast}\right\}  ^{\prime}=\mathbb{C}%
\mkern2mu%
\openone$ whatever the choice of
$V_{i}$'s. One choice is
\[
V_{1}=
\begin{pmatrix}
1 & 0 & 0\\
0 & 0 & 0\\
0 & 0 & 0
\end{pmatrix}
,\;V_{2}=
\begin{pmatrix}
0 & 0 & 0\\
0 & 1 & 0\\
0 & 0 & 0
\end{pmatrix}
,\;V_{3}=
\begin{pmatrix}
0 & 0 & 0\\
0 & 0 & 0\\
\frac{1}{\sqrt{2}} & 0 & 0
\end{pmatrix}
,\;V_{4}=
\begin{pmatrix}
0 & 0 & 0\\
0 & 0 & 0\\
0 & \frac{1}{\sqrt{2}} & 0
\end{pmatrix}
.
\]
The invariant states are all 
the
convex combinations of the two states $\left(
X_{ij}\right)  \mapsto X_{11}$ or $\left(  X _{ij}\right)  \mapsto X
_{22}$. 
Thus, in the general situation, the following theorem is the best possible.

\begin{theorem}
\label{Thm3.5}Let $\varphi=\sum_{i}V_{i}^{{}}\,\cdot\,V_{i}^{\ast}$ be a
normal unital completely positive map of $\mathcal{B}\left(  \mathcal{K}%
\right)  $. Then%
\[
\left\{  V_{i}^{{}},V_{i}^{\ast}\right\}  ^{\prime}\subset\mathcal{B}\left(
\mathcal{H}\right)  ^{\varphi}.
\]
Furthermore, the space $\mathcal{B}\left(  \mathcal{H}\right)  ^{\varphi}$
contains a largest $\ast$-subalgebra, and this algebra is $\left\{  V_{i}^{{}%
},V_{i}^{\ast}\right\}  ^{\prime}$.
\end{theorem}

\begin{proof}
Since $\sum_{i}V_{i}^{{}}V_{i}^{\ast}=\openone$, the first assertion is
trivial. Next note that as $\varphi$ is normal, if $\mathfrak{A}$ is a $\ast
$-subalgebra of $\mathcal{B}\left(  \mathcal{H}\right)  ^{\varphi}$, then the
weak$^{\ast}$-closure $\overline{\mathfrak{A}}$ of $\mathfrak{A}$ is contained
in $\mathcal{B}\left(  \mathcal{H}\right)  ^{\varphi}$. But since
$\overline{\mathfrak
{A}}$ is the weak$^{\ast}$-closure of the linear span of its projections, it
follows from Lemma \ref{Lem3.3} that%
\[
\mathfrak{A}\subset\overline{\mathfrak{A}}\subset\left\{  V_{i}^{{}}%
,V_{i}^{\ast}\right\}  ^{\prime}.
\]
This proves Theorem \ref{Thm3.5}.
\end{proof}

\chapter{\label{Lifting}The commutant lifting theorem and pure states on
$\mathcal{O}_{d}$}

The main aim of this \chaptername{} is to decide which systems $\left(
\mathcal{K},\Omega,V_{1},\dots,V_{d}\right)  $ give rise to pure states
$\hat{\omega}$ on $\mathcal{O}_{d}$. To this end it will be convenient to
define a completely positive unital map $\bs$ of $\mathcal{B}\left(
\mathcal{K}\right)  $ by%
\begin{equation}
\bs\left(  A\right)  =\sum_{i=1}^{d}V_{i}^{{}}AV_{i}^{\ast}. \label{eq4.1}%
\end{equation}
We will actually establish an order isomorphism between the order interval
$\left[  0,\hat{\omega}\right]  $ in the set of positive functionals on
$\mathcal{O}_{d}$, and the set of operators $A\in\mathcal{B}\left(
\mathcal{K}\right)  $ such that $0\leq A\leq\openone$ and $\bs
\left(  A\right)  =A$. This is a natural generalization of the commutant
lifting theorem of \cite{Pop92}, and another version of this result is
Corollary 5.4 in \cite{BrJo97a}. The term ``commutant lifting'' is from
single-operator theory \cite{SzFo70,FoFr84,DMP68} where it refers to the
Sz.-Nagy lifting theorem, which for every contractive operator $V$ in a given
Hilbert space $\mathcal{K}$ yields a minimal coisometry, and in fact, by a
second step, also a unitary operator $U$, acting on a bigger Hilbert space
$\mathcal{H}$, and serving as a lifting of $V$. If $P$ denotes the projection
of $\mathcal{H}$ onto $\mathcal{K}$, i.e., $\mathcal{K}=P\mathcal{H}$, then
Sz.-Nagy's dilation theorem states the existence of $\left(  U,\mathcal{H}%
\right)  $ such that $UP=PUP$ and $V^{n}=PU^{n}P$ on $\mathcal{K}$ for all
$n\in\mathbb{N}$. ``Minimality'' here is the requirement that the subspace
$\mathcal{K}$ be cyclic for $\left\{  U^{n}\mid n\in\mathbb{Z}\right\}  $ in
$\mathcal{H}$. If we have two contractions $V_{i}\colon\mathcal{K}%
_{i}\rightarrow\mathcal{K}_{i}$, $i=1,2$, with corresponding minimal
coisometric (or unitary) dilations $\left(  U_{i},\mathcal{H}_{i}\right)  $
and projections $P_{i}\colon\mathcal{H}_{i}\rightarrow\mathcal{K}_{i}$,
$P_{i}\mathcal{H}_{i}=\mathcal{K}_{i}$, $i=1,2$, and if $Y\colon
\mathcal{K}_{1}\rightarrow\mathcal{K}_{2}$ is a bounded operator which is
given to intertwine the two contractions, i.e., $YV_{1}=V_{2}Y$, then $Y$
lifts, by \cite{DMP68}, to a bounded $X\colon\mathcal{H}_{1}\rightarrow
\mathcal{H}_{2}$ with the same operator norm, $\left\|  X\right\|  =\left\|
Y\right\|  $, and satisfying $XU_{1}=U_{2}X$, and $P_{2}XP_{1}=Y$ on
$\mathcal{K}_{1}$.

The analogy to the present setting refers to an operator $Y$ which intertwines
two given $V_{i}$-systems $\left\{  V_{i}\right\}  _{i=1}^{d}$ and $\left\{
W_{i}\right\}  _{i=1}^{d}$, say, and its canonical lifting to an operator
which intertwines the corresponding two representations of the Cuntz algebra
$\mathcal{O}_{d}$.

Proposition \ref{Pro4.1} and Theorem \ref{Thm4.4} below represent our
multivariable analogue of this lifting result, but only for the special case
when $V_{i}=W_{i}$, while Theorem \ref{Thm5.1} is our general multivariable
commutant lifting theorem.

\begin{proposition}
\label{Pro4.1}Adopt the notation in Remark\/ \textup{\ref{Rem2.3}}. Then the
selfadjoint part of the commutant $\pi\left(  \mathcal{O}_{d}\right)
^{\prime}$ is norm and order isomorphic to the space of selfadjoint fixed
points of the completely positive map $\bs$. This isomorphism takes $A\in
\pi\left(  \mathcal{O}_{d}\right)  ^{\prime}$ into $PAP\in\mathcal{B}\left(
\mathcal{K}\right)  ^{\bs}$.
\end{proposition}

\begin{proof}
Let $X\in\pi\left(  \mathcal{O}_{d}\right)  ^{\prime}$. Then $PXP$ is
determined by the matrix elements%
\[
\ip{V_{I}^{\ast}\Omega}{P X  PV_{J}^{\ast}\Omega}=\ip{S_{I}^{\ast}\Omega}{ X
S_{J}^{\ast}\Omega}.
\]
Writing the same expression for $I\mapsto Ii$ and $J\mapsto Ji$, and summing
over $i$, shows that $PXP\in\mathcal{B}\left(  \mathcal{K}\right)  ^{\bs}$.
(Why can we not conclude $PXP\in\left\{  V_{i}^{{}},V_{i}^{\ast}\right\}
^{\prime}$? We get
\begin{multline*}
\ip{V_{I}^{\ast}\Omega}{V_{i}^{}P X PV_{J}^{\ast}\Omega}=\ip{S_{i}^{\ast}%
S_{I}^{\ast}\Omega}{ X S_{J}^{\ast}\Omega}\\
=\ip{S_{I}^{\ast}\Omega}{ X S_{i}^{}S_{J}^{\ast}\Omega}=\ip{V_{I}^{\ast}%
\Omega}{P X S_{i}^{}PV_{J}^{\ast}\Omega}\\
\neq\ip{V_{I}^{\ast}\Omega}{P X PS_{i}^{}S_{J}^{\ast}\Omega},
\end{multline*}
because $S_{i}^{{}}S_{J}^{\ast}\Omega\notin\mathcal{K}$
in general. 
See Theorem
\ref{Thm3.5}.) Conversely, assume that $D\in\mathcal{B}\left(  \mathcal{K}%
\right)  ^{\bs}$, with $0\leq D\leq\openone$. That is, $D$ satisfies the
hypothesis of Corollary \ref{Cor2.4}. Hence the linear functional on
$\mathcal{O}_{d}$ defined by%
\[
\tilde{\omega}\left(  s_{I}^{{}}s_{J}^{\ast}\right)  =\ip{V_{I}^{\ast}\Omega
}{DV_{J}^{\ast}\Omega}%
\]
is positive, and, applying the same argument to $\openone-D$, we find that
$0\leq\tilde{\omega}\leq
\hat{\omega}$. 
Since $\Omega$ is cyclic for $\mathcal{O}%
_{d}$, there is an $X\in\pi\left(  \mathcal{O}_{d}\right)  ^{\prime}$ (with
$0\leq X\leq\openone$), which is uniquely determined by the equation%
\[
\tilde{\omega}\left(  s_{I}^{{}}s_{J}^{\ast}\right)  =\ip{\Omega}{ X  S_{I}%
^{}S_{J}^{\ast}\Omega}.
\]
But since $X\in\mathcal{O}_{d}^{\prime}$, we have $\ip{\Omega}{ X  S_{I}%
^{}S_{J}^{\ast}\Omega}=\ip{V_{I}^{\ast}\Omega}{ X  V_{J}^{\ast}\Omega}%
=\ip{V_{I}^{\ast}\Omega}{DV_{J}^{\ast}\Omega}$. That is to say $PXP=D$.

Since $\openone\in\mathcal{B}\left(  \mathcal{K}\right)  ^{\bs}$, the real
linear span of 
this 
positive cone in $\mathcal{B}\left(  \mathcal{K}\right)  ^{\bs
}$ is all of the selfadjoint part, and hence the map $X\mapsto PXP$ is onto,
and (by scaling with suitable positive factors) the above arguments show that
the map is an order isomorphism between 
the respective
selfadjoint parts of $\pi\left(
\mathcal{O}_{d}\right)  ^{\prime}$ and $\mathcal{B}\left(  \mathcal{K}\right)
^{\bs}$. The selfadjoint subspaces are also order unit spaces, i.e.,%
\[
\left\|  A\right\|  =\inf\left\{  \alpha\geq0\mid-\alpha\openone\leq
A\leq\alpha\openone\right\}  .
\]
(For $\mathcal{B}\left(  \mathcal{K}\right)  ^{\bs}$, this formula is inherited
from $\mathcal{B}\left(  \mathcal{K}\right)  $, using, of course, crucially
that $\openone\in\mathcal{B}\left(  \mathcal{K}\right)  ^{\bs}$.) {}From this
it is evident that the isomorphism is also isometric.
\end{proof}

Having now identified $\mathcal{B}\left(  \mathcal{K}\right)  ^{\bs}$ with
$P\pi\left(  \mathcal{O}_{d}\right)  ^{\prime}P$, let us return to the
question raised in Theorem \ref{Thm3.5} and the preceding remarks on when
$\mathcal{B}\left(  \mathcal{K}\right)  ^{\bs}$ is an algebra.

\begin{proposition}
\label{Pro4.2}Let $\mathcal{M}$ be a von Neumann algebra on a Hilbert space
$\mathcal{H}$, and let $P$ be a projection in $\mathcal{H}$ such that
$X\mapsto PXP$ is an isometry on the selfadjoint part of $\mathcal{M}$. Then
the following are equivalent.

\begin{enumerate}
\item \label{Pro4.2(1)}$P\mathcal{M}P$ is an algebra.

\item \label{Pro4.2(2)}$X\mapsto PXP$ is a homomorphism on $\mathcal{M}$.

\item \label{Pro4.2(3)}$P\in\mathcal{M}^{\prime}$.
\end{enumerate}
\end{proposition}

\begin{proof}
(\ref{Pro4.2(3)}) $\Rightarrow$ (\ref{Pro4.2(2)}) $\Rightarrow$
(\ref{Pro4.2(1)}) are trivial. (\ref{Pro4.2(2)}) $\Rightarrow$
(\ref{Pro4.2(3)}) follows from the observation that the homomorphism property
implies $PX^{\ast}\left(  \openone-P\right)  XP=0$, i.e., $\left(
\openone-P\right)  XP=0$, and $XP=PXP=PX$. Note that these steps do not even
depend on the isometry property.

The nontrivial bit, (\ref{Pro4.2(1)}) $\Rightarrow$ (\ref{Pro4.2(2)}), is
essentially 
contained
in the proof of Theorem \ref{Thm3.5}. Here is a slightly different
way of putting it: The isometry property means that the unit interval of
$\mathcal{M}$ is isometrically mapped onto that of the algebra $P\mathcal{M}%
P$. In particular, extremal points correspond to extremal points, which in a
von Neumann algebra means that projections go into projections, and
orthogonality of projections is preserved. By the spectral theorem, we find
that the compression map is a Jordan isomorphism, and hence the direct sum of
a homomorphism and an anti-homomorphism. Because it is completely positive, it
is a homomorphism.
\end{proof}

\begin{remark}
\label{Rem4.3}Note that it is not enough to require isometry on the whole
\textup{(}complex\textup{)} vector space $\mathcal{M}$, since the norms alone
do not give enough information. A counterexample can be made with a
two-dimensional abelian algebra. \textup{(}$\mathcal{M}=\mathbb{C}%
\mkern2mu%
%
Q+\mathbb{C}\left(  \openone-Q\right)  $, and $QPQ$ has both eigenvalues $0$
and $1$.\textup{)}

Note also that if Propositions \ref{Pro4.1} and \ref{Pro4.2} are applied to
the example of Arveson discussed prior to Theorem \ref{Thm3.5}, it follows
that the dilation of $\left\{  \mathbb{C}^{3},V_{1},\dots,V_{4}\right\}  $ to
a representation $\pi$ of $\mathcal{O}_{4}$ decomposes into two disjoint
irreducible representations, and $P$ is not contained in $\pi\left(
\mathcal{O}_{4}\right)  ^{\prime\prime}$. Note that the nontrivial projections
in $\pi\left(  \mathcal{O}_{4}\right)  ^{\prime}$ then cannot be Popescu
dilations of anything in $\mathcal{B}\left(  \mathbb{C}^{3}\right)  $.
\end{remark}

We are now ready to state the characterization of pure states $\hat{\omega}$
on $\mathcal{O}_{d}$. If $\hat{\omega}$ is any state on $\mathcal{O}_{d}$, let
again $\left(  \mathcal{H},\Omega,S_{1},\dots,S_{d}\right)  $ be the
corresponding representation, and $\left(  \mathcal{K},\Omega,V_{1}%
,\dots,V_{d}\right)  $ the corresponding Popescu system, and define the
corresponding endomorphism $\sigma$ of $\mathcal{B}\left(  \mathcal{H}\right)
$ by $\sigma\left(  \,\cdot\,\right)  =\sum_{i=1}^{d}S_{i}^{{}}\,\cdot
\,S_{i}^{\ast}$, and the unital completely positive map $\bs$ of $\mathcal{B}%
\left(  \mathcal{K}\right)  $ by $\bs\left(  \,\cdot\,\right)  =\sum_{i=1}%
^{d}V_{i}^{{}}\,\cdot\,V_{i}^{\ast}$.

\begin{theorem}
\label{Thm4.4}If $\hat{\omega}$ is a state on $\mathcal{O}_{d}$, the following
conditions are equivalent.

\begin{enumerate}
\item \label{Thm4.4(1)}$\hat{\omega}$ is pure.

\item \label{Thm4.4(2)}$\sigma\left(  X\right)  =X$ implies $X\in
\mathbb{C}%
\mkern2mu%
\openone_{\mathcal{H}}$, $X\in\mathcal{B}\left(  \mathcal{H}\right)
$.

\item \label{Thm4.4(3)}$\bs\left(  Y\right)  =Y$ implies $Y\in\mathbb{C}%
\mkern2mu%
%
\openone_{\mathcal{K}}$, $Y\in\mathcal{B}\left(  \mathcal{K}\right)  $.

\item \label{Thm4.4(4)}$\left\{  V_{i}^{{}},V_{i}^{\ast}\right\}  $ acts
irreducibly on $\mathcal{K}$, and $P\in\pi\left(  \mathcal{O}_{d}\right)
^{\prime\prime}$.
\end{enumerate}
\end{theorem}

\begin{proof}
In general the fixed point algebra for $\sigma$ is $\pi\left(  \mathcal{O}%
_{d}\right)  ^{\prime}$ (see, e.g., \cite[formula (3.5)]{BJP96} or
\cite[Proposition 3.1]{Lac93a}) and hence (\ref{Thm4.4(1)}) $\Leftrightarrow$
(\ref{Thm4.4(2)}), and (\ref{Thm4.4(1)}) $\Leftrightarrow$ (\ref{Thm4.4(3)})
follows from Proposition \ref{Pro4.1}. If $\hat{\omega}$ is pure, then
$P\in\pi\left(  \mathcal{O}_{d}\right)  ^{\prime\prime}$ and $\left\{
V_{i}^{{}},V_{i}^{\ast}\right\}  $ acts irreducibly on $\mathcal{K}$ because
$\left\{  V_{i}^{{}},V_{i}^{\ast}\right\}  ^{\prime}\subset\mathcal{B}\left(
\mathcal{K}\right)  ^{\bs}=\mathbb{C}%
\mkern2mu%
\openone_{\mathcal{H}}$ by
(\ref{Thm4.4(3)}), hence (\ref{Thm4.4(1)}) $\Rightarrow$ (\ref{Thm4.4(4)}).
Conversely if $P\in\pi\left(  \mathcal{O}_{d}\right)  ^{\prime\prime}$ it
follows, by applying Proposition \ref{Pro4.2} on $\mathcal{M}=\pi\left(
\mathcal{O}_{d}\right)  ^{\prime}$, that $P\pi\left(  \mathcal{O}_{d}\right)
^{\prime}P$ is an algebra. But this algebra is $\mathcal{B}\left(
\mathcal{K}\right)  ^{\bs}$ by Proposition \ref{Pro4.1}, and if $\left\{
V_{i}^{{}},V_{i}^{\ast}\right\}  $ acts irreducibly, it follows from Theorem
\ref{Thm3.5} that $\mathcal{B}\left(  \mathcal{K}\right)  ^{\bs}%
=\mathbb{C}%
\mkern2mu%
\openone_{\mathcal{K}}$. Thus (\ref{Thm4.4(4)}) $\Rightarrow$
(\ref{Thm4.4(3)}), and Theorem \ref{Thm4.4} is proved.
\end{proof}

\chapter{\label{Representations}Representations of $\mathcal{O}_{d}$}

For the wavelet applications 
described in \Chaptername{} \ref{Intro},
we will need versions of Theorem \ref{Thm2.1},
Proposition \ref{Pro4.1} and Theorem \ref{Thm4.4} where the state $\hat
{\omega}$ is replaced merely by the system $\left(  \mathcal{K},V_{1}%
,\dots,V_{d}\right)  $.

\begin{theorem}
\label{Thm5.1}Let $\mathcal{K}$ be a Hilbert space, and let $V_{1},\dots
,V_{d}\in\mathcal{B}\left(  \mathcal{K}\right)  $ be operators satisfying%
\[
\sum_{i\in\mathbb{Z}_{d}}V_{i}^{{}}V_{i}^{\ast}=\openone.
\]
Then $\mathcal{K}$ can be embedded into a larger Hilbert space $\mathcal{H}%
=\mathcal{H}_{V}$ carrying a representation $S_{1},\dots,S_{d}$ of the Cuntz
algebra $\mathcal{O}_{d}$ such that if $P\colon\mathcal{H}\rightarrow
\mathcal{K}$ is the projection onto $\mathcal{K}$ we have%
\[
V_{i}^{\ast}=S_{i}^{\ast}P
\]
\textup{(}i.e., $S_{i}^{\ast}\mathcal{K}\subset\mathcal{K}$ and $S_{i}^{\ast
}P=PS_{i}^{\ast}P=V_{i}^{\ast}$\textup{)} and $\mathcal{K}$ is cyclic for the
representation. The system $\left(  \mathcal{H},S_{1},\dots,S_{d},P\right)  $
is unique up to unitary equivalence, and if $\bs\colon\mathcal{B}\left(
\mathcal{K}\right)  \rightarrow\mathcal{B}\left(  \mathcal{K}\right)  $ is
defined by%
\[
\bs\left(  A\right)  =\sum_{i}V_{i}^{{}}AV_{i}^{\ast},
\]
then the commutant of the representation $\left\{  S_{1},\dots,S_{d}\right\}
^{\prime}$ is isometrically order isomorphic to the fixed point set
$\mathcal{B}\left(  \mathcal{K}\right)  ^{\bs}=\left\{  A\in\mathcal{B}\left(
\mathcal{K}\right)  \mid\bs\left(  A\right)  =A\right\}  $ by the map
$A^{\prime}\mapsto PA^{\prime}P$. More generally, if $W_{1},\dots,W_{d}%
\in\mathcal{B}\left(  \mathcal{K}\right)  $ is another set of operators
satisfying%
\[
\sum_{i\in\mathbb{Z}_{d}}W_{i}^{{}}W_{i}^{\ast}=\openone
\]
and $T_{1},\dots,T_{d}$ are the corresponding representatives of $s_{1}%
,\dots,s_{d}$, then there is an isometric linear isomorphism between
intertwiners $U\colon\mathcal{H}_{V}\rightarrow\mathcal{H}_{W}$, i.e.,
operators satisfying%
\[
US_{i}=T_{i}U,
\]
and operators $V\in\mathcal{B}\left(  \mathcal{K}\right)  $ such that%
\begin{equation}
\sum_{i\in\mathbb{Z}_{d}}W_{i}^{{}}VV_{i}^{\ast}=V, \label{eq5.1}%
\end{equation}
given by the map $U\mapsto V=PUP$.
\end{theorem}

\begin{proof}
Inspecting the proof of Lemma \ref{Lem2.2}, we see that the vector $\Omega$
plays no role in the proof, so the map $R\colon\mathcal{O}_{d}\rightarrow
\mathcal{B}\left(  \mathcal{K}\right)  $ defined by%
\[
R\left(  s_{I}^{{}}s_{J}^{\ast}\right)  =V_{I}^{{}}V_{J}^{\ast}%
\]
is well defined and completely positive. The representation $S_{1},\dots
,S_{d}$ of $\mathcal{O}_{d}$ on $\mathcal{H}$ thus may be taken to be the
Stinespring dilation of $R$ \cite{Arv69}, \cite[p.\ 229, Notes and Remarks to
Chapter 5]{BrRoII}, \cite{Sti55}, and uniqueness up to unitary equivalence
follows from uniqueness of the Stinespring representation.

The commutant lifting property is established as in Proposition \ref{Pro4.1},
using Corollary \ref{Cor2.4}.

To establish the final intertwiner lifting property, one considers the
direct sum
representation of $\mathcal{O}_{d}$ on $\mathcal{H}_{V}\oplus\mathcal{H}_{W}$
given by%
\[
s_{i}\longmapsto S_{i}\oplus T_{i}.
\]
Note that
some operator
$U\colon\mathcal{H}_{V}\rightarrow\mathcal{H}_{W}$ is an intertwiner
if and only if $\left(
\begin{smallmatrix}
0 & 0\\
U & 0
\end{smallmatrix}
\right)  $ is in the commutant of this
sum
representation. But the operators
corresponding to $V_{i}$ of this latter representation, relative to the
subspace $\mathcal{K}\oplus\mathcal{K}\subset\mathcal{H}_{V}\oplus
\mathcal{H}_{W}$, are
\[%
\begin{pmatrix}
V_{i} & 0\\
0 & W_{i}%
\end{pmatrix}
,\quad i\in\mathbb{Z}_{d},
\]
so, using the commutant lifting property of the direct sum representation, one
verifies that $U$ intertwines the $S_{i}$'s and the $T_{i}$'s if and only if
$V=PUP$ is fixed under the map $\sum_{i}W_{i}^{{}}\,\cdot\,V_{i}^{\ast}$.
Specifically, if
$\beta \left( \,\cdot \,\right) :=\sum_{i}T_{i}^{{}}\cdot S_{i}^{*}$,
then we have the identity
\[P\beta \left( X\right) P=\sum_{i}W_{i}^{{}}PXPV_{i}^{*},
\]
valid for all operators
$X\colon \mathcal{H}_{V}\rightarrow \mathcal{H}_{W}$.
Now note that $U$ intertwines the two
$\mathcal{O}_{d}$-representations, if and only if
$\beta \left( U\right) =U$,
and the assertion follows from this.
\end{proof}

\begin{remark}
\label{RemNew5.2}Another more direct way of constructing the representation of
$\mathcal{O}_{d}$ in Theorem \textup{\ref{Thm5.1}} is the following: Let
$\mathcal{I}_{n}$ be the set of finite sequences $I=\left(  i_{1},\dots
,i_{m}\right)  $ where $m\leq n$ and $i_{k}\in\mathbb{Z}_{d}$ for all $k$
\textup{(}including the empty sequence\textup{)}, and let $\mathbb{C}%
\mkern2mu%
\mathcal{I}_{n}$ be the complex linear space of formal linear combinations of
elements in $\mathcal{I}_{n}$. Put $\mathcal{I}=\bigcup_{n}\mathcal{I}_{n}$ as
in \Chaptername{} \textup{\ref{General}}, and define%
\[
H_{n}=\mathbb{C}%
\mkern2mu%
\mathcal{I}_{n}\otimes\mathcal{K}%
\]
\textup{(}algebraic tensor product\textup{)}. For each $I\in\mathcal{I}$,
define a linear operator $S_{I}$ on $H=\bigcup_{n}H_{n}$ by%
\[
S_{I}\left(  J\otimes\xi\right)  =IJ\otimes\xi
\]
and linearity. Define a semi-inner product on $H$ by requiring%
\begin{align*}
\ip{I\otimes\xi}{IJ\otimes\eta} &=\ip{\xi}{V_{J}\eta},  \\
\ip{IJ\otimes\xi}{I\otimes\eta} &=\ip{V_{J}\xi}{\eta}  \\
\intertext{for all $I,J\in\mathcal{I}$, $\xi,\eta\in\mathcal{K}$, and}%
\ip{I\otimes\xi}{J\otimes\eta} &=0
\end{align*}%
if the pair $I,J$ does not have one of the forms above. To show that this
sesquilinear form is indeed positive and well defined, we proceed by
induction: This is true for $H_{0}=\mathcal{K}$. Suppose this is proved for
$H_{n-1}$ and let $\zeta\in H_{n}$. We express $\zeta$ as%
\[
\zeta=\sum_{j\in\mathbb{Z}_{d}}S_{j}\zeta_{j}+\zeta_{0}%
\]
where $\zeta_{j}\in H_{n-1}$ and $\zeta_{0}\in H_{0}=\mathcal{K}$. Then%
\begin{align*}
\ip{\zeta}{\zeta} &=\vphantom{\sum_{jk}}
\ip{\vphantom{\sum}\smash{\sum_{j}}S_{j}\zeta_{j}+\zeta_{0}}
{\smash{\sum_{k}}S_{k}\zeta_{k}+\zeta_{0}}  \\
&=\sum_{j}\ip{\zeta_{j}}{\zeta_{j}}+\sum_{j}\ip{\zeta_{j}^{{}}}{V_{j}^{\ast
}\zeta_{0}^{{}}}+\sum_{k}\ip{V_{k}^{\ast}\zeta_{0}^{{}}}
{\zeta_{k}^{{}}}+\ip{\zeta_{0}}{\zeta_{0}}  \\
&=\sum_{j}\left\|  \zeta_{j}^{{}}+V_{j}^{\ast}\zeta_{0}^{{}}\right\|  ^{2}%
\geq0.
\end{align*}%
Let $\mathcal{H}$ be the completion of $H$ modulo zero-vectors and
$\Lambda\colon H\rightarrow\mathcal{H}$ the canonical map. We define a bounded
operator $S_{i}$ on $\mathcal{H}$ by%
\[
S_{i}\Lambda\left(  \zeta\right)  =\Lambda\left(  S_{i}\zeta\right)  ,
\]
and, using $\sum_{j}V_{j}^{{}}V_{j}^{\ast}=\openone$, one easily verifies that
$s_{i}\mapsto S_{i}$ is a representation of $\mathcal{O}_{d}$ satisfying the
required properties.
\end{remark}

\begin{remark}
\label{RemNew5.3}Note that Corollary \textup{\ref{Cor2.4}} can also be proved
along the lines in \textup{Remark \ref{RemNew5.2}}, but now one defines the
semi-inner product $\ip{\,\cdot
\,}{\,\cdot\,}_{D}$ on $\mathcal{H}$ by requiring%
\begin{align*}
\ip{I\otimes\xi}{IJ\otimes\eta}_{D}  &  =\ip{V_{J}^{*}\xi}{D\eta},\\
\ip{IJ\otimes\xi}{I\otimes\eta}_{D}  &  =\ip{\xi}{DV_{J}^{*}\eta},\text{\quad
etc.}%
\end{align*}
\end{remark}

\begin{remark}
\label{RemNew5.4}In comparing the single-operator commutant lifting
\cite{DMP68} with our Theorem \textup{\ref{Thm5.1}}, we note that the naive
\textup{(}or natural\textup{)} multivariable generalization of the
intertwining property for an operator $Y$ on the ``small'' Hilbert space
$\mathcal{K}$ would be $YV_{i}=W_{i}Y$. But this property is slightly
different from the present one, $\sum_{i}W_{i}^{{}}YV_{i}^{\ast}=Y$, i.e.,
\textup{(\ref{eq5.1})}, used in Theorem \ref{Thm5.1}.

There is naturally a variety of
ways of generalizing the classical
single-operator commutant lifting
theorem to several variables, each
serving different purposes. In
addition to ours and the others
mentioned above, there are
related, but different, approaches
\textup{(}to the multivariable theory\textup{)} in
recent papers by Arveson \cite{Arv97}
and Bhat \cite{Bha97}.
\end{remark}

\chapter{\label{Restriction}Irreducible representations of $\mathcal{O}_{d}$
and their restriction to $\operatorname*{UHF}\nolimits_{d}$}

Consider an irreducible representation $S_{1},\dots,S_{d}$ of $\mathcal{O}%
_{d}$ on a Hilbert space $\mathcal{H}$, and let $\mathcal{K}$ be a cyclic
subspace of $\mathcal{H}$ invariant under $S_{1}^{\ast},\dots,S_{d}^{\ast}$.
Define again $V_{1},\dots,V_{d}\in\mathcal{B}\left(  \mathcal{K}\right)  $ by%
\[
V_{i}^{\ast}P=S_{i}^{\ast}P=PS_{i}^{\ast}P
\]
where $P\colon\mathcal{H}\rightarrow\mathcal{K}$ is the projection onto
$\mathcal{K}$. By Theorem \ref{Thm5.1}, irreducibility
on $\mathcal{H}$
is equivalent to
ergodicity of the completely positive map $\bs$ defined on $\mathcal{B}\left(
\mathcal{K}\right)  $ by $\bs\left(  \,\cdot\,\right)  =\sum_{i}V_{i}^{{}%
}\,\cdot\,V_{i}^{\ast}$. Since $\bs\left(  \openone\right)  =\openone$, $\bs$
maps the state space of $\mathcal{B}\left(  \mathcal{K}\right)  $ into itself,
and hence there is a $\bs$-invariant state $\varphi$. If $\mathcal{K}$ is
finite-dimensional, we will show that $\varphi$ is unique. The state $\varphi$
is automatically normal since $\mathcal{K}$ is finite-dimensional. Let $E$ be
the support of $\varphi$.

\begin{lemma}
\label{Lem6.1}$S_{i}^{\ast}E\mathcal{K}\subset E\mathcal{K}$ for all
$i\in\mathbb{Z}_{d}$.
\end{lemma}

\begin{proof}
Since%
\[
\varphi\left(  \bs\left(  E\right)  \right)  =\varphi\left(  E\right)
=\openone
\]
and $0\leq\bs\left(  E\right)  \leq\openone$, it follows that%
\[
\bs\left(  E\right)  \geq E.
\]
Applying Lemma \ref{Lem3.1}, (\ref{Lem3.1(1)}) $\Rightarrow$ (\ref{Lem3.1(2)}%
), on $p=\openone-E$ gives%
\[
V_{i}^{\ast}E=EV_{i}^{\ast}E.
\]
This proves Lemma \ref{Lem6.1}.
\end{proof}

But from \cite[Lemma 6.3]{BrJo97a}, it follows that there is only one $\bs
$-invariant state with support inside $E$. So, if $\varphi_{1}$, $\varphi_{2}$
are two $\bs$-invariant states with respective support projections $E_{1}$,
$E_{2}$, then $\frac{1}{2}\left(  \varphi_{1}+\varphi_{2}\right)  $ is a $\bs
$-invariant state with support $E_{1}\vee E_{2}$, and hence $\varphi_{1}%
=\frac{1}{2}\left(  \varphi_{1}+\varphi_{2}\right)  =\varphi_{2}$ by the
argument above. We have proved:

\begin{lemma}
\label{Lem6.2}If $\mathcal{K}$ is finite-dimensional, then $\mathcal{B}\left(
\mathcal{K}\right)  $ has a unique $\bs$-invariant state when $\bs$ is ergodic.
\end{lemma}

The example after the proof of Lemma \ref{Lem3.4} shows that this $\bs
$-invariant state need not be faithful. However, replacing $P$ by the support
$E$ of $\varphi$, and using Lemma \ref{Lem6.1}, the following theorem is
applicable to general irreducible representations when $\mathcal{K}$ is
finite-dimensional, replacing $P$ by $E$.

\begin{theorem}
\label{Thm6.3}\raggedright
Consider an irreducible representation of $\mathcal{O}_{d}$ on
$\mathcal{H}$, and let $\mathcal{K},V_{1},\dots,V_{d},P,\bs$ be as in the
introduction to this \chaptername{}. 
Assume that there exists a normal faithful $\bs
$-invariant state $\varphi$ on $\mathcal{B}\left(  \mathcal{K}\right)  $. Let
$\psi$ be the state of $\mathcal{O}_{d}$ defined by%
\[
\psi\left(  s_{I}^{{}}s_{J}^{\ast}\right)  =\varphi\left(  V_{I}^{{}}%
V_{J}^{\ast}\right)  .
\]
The following three subsets of the circle group $\mathbb{T}$ are equal.

\begin{enumerate}
\item \label{Thm6.3(1)}$\left\{  t\in\mathbb{T}\mid\psi\circ\tau_{t}%
=\psi\right\}  $, where $\tau$ is the gauge action.

\item \label{Thm6.3(2)}$\left\{  t\in\mathbb{T}\mid\psi\circ\tau_{t}\text{ is
quasi-equivalent to }\psi\right\}  $.

\item \label{Thm6.3(3)}$\operatorname{PSp}\left(  \bs\right)  \cap\mathbb{T}$,
where $\operatorname{PSp}\left(  \bs\right)  $ is the set of eigenvalues of
$\bs$.
\end{enumerate}

Furthermore, this set is a finite subgroup of $\mathbb{T}$. If $k$ is the
order of this subgroup, the restriction of the representation to
$\operatorname*{UHF}\nolimits_{d}$ decomposes into $k$ mutually disjoint
irreducible representations, and these are mapped cyclically into each other
by the one-sided shift $\lambda\left(  \,\cdot\,\right)  =\sum_{i}s_{i}^{{}%
}\,\cdot\,s_{i}^{\ast}$.
\end{theorem}

\begin{remark}
\label{Rem6.4}Since the normal states on $\mathcal{B}\left(  \mathcal{K}%
\right)  $ are given by density matrices, it follows from \cite[Lemma
6.3]{BrJo97a} \textup{(}as in the proof of Lemma \ref{Lem6.2} above\textup{)}
that if there is a faithful $\bs$-invariant normal state $\varphi$, then this
is the unique $\bs$-invariant normal state. Note that the state $\psi$ defined
in Theorem \ref{Thm6.3} is well defined by Lemma \ref{Lem2.2}, and
\[
\psi\circ\lambda=\psi
\]
since $\varphi\circ\bs=\varphi$.
\end{remark}

During the proof of Theorem \ref{Thm6.3} we will establish that%
\[
\operatorname{PSp}\left(  \bs\right)  \cap\mathbb{T}=\operatorname{PSp}\left(
\lambda\right)  \cap\mathbb{T}%
\]
and that each of the corresponding eigenspaces is spanned by a unitary
operator (in $\mathcal{B}\left(  \mathcal{K}\right)  $, $\mathcal{B}\left(
\mathcal{H}\right)  $, respectively), and this unitary operator in
$\mathcal{B}\left(  \mathcal{H}\right)  $ implements $\tau_{t}$ if
$\bar{t}$ is the eigenvalue. In fact, if $U\in\mathcal{B}\left(
\mathcal{H}\right)  \diagdown\left\{  0\right\}  $ and $\lambda\left(
U\right)  =\bar{t}U$, then $\lambda\left(  U^{\ast}U\right)  =\lambda\left(
U\right)  ^{\ast}\lambda\left(  U\right)  =U^{\ast}U$, hence $U^{\ast}U$, and
likewise $UU^{\ast}$, is a scalar multiple of $\openone$. Thus, renormalizing
$U$, we may take $U$ to be unitary. But $\lambda\left(  U\right)  =\sum
_{i}S_{i}^{{}}US_{i}^{\ast}=\bar{t}U$, so $S_{i}U=\bar{t}US_{i}$ and hence
\[
US_{i}U^{\ast}=tS_{i}=\tau_{t}\left(  S_{i}\right)  .
\]
Conversely, if $U$
implements $\tau_{t}$, then $\lambda\left(  U\right)  =\bar{t}U$, and we have
shown%
\[
\operatorname{PSp}\left(  \lambda\right)  \cap\mathbb{T}=\left\{  t\mid
\tau_{t}\text{ is inner}\right\}  .
\]
In particular, this shows that $\operatorname{PSp}\left(  \bs
\right)  \cap\mathbb{T}$ is independent of the particular state $\varphi$
chosen (with the required properties).

\begin{proof}
[Proof of Theorem \textup{\ref{Thm6.3}}]
We first prove the inclusion

\begin{lemma}
\label{Lem6.5}$\left\{  t\in\mathbb{T}\mid\psi\circ\tau_{t}=\psi\right\}
\subset\operatorname{PSp}\left(  \bs\right)  \cap\mathbb{T}$.
\end{lemma}

\begin{proof}
Since $\varphi$ is normal on $\mathcal{B}\left(  \mathcal{K}\right)  $,
$\varphi$ is a (possibly infinite) convex combination of vector states, and
thus $\psi$ is a convex combination of vector states. Since the given
representation of $\mathcal{O}_{d}$ on $\mathcal{H}$ is irreducible, $\psi$ is
a type $\mathrm{I}$ factor state. If $\psi\circ\tau_{t}=\psi$, it follows that
there is a unitary $\bar{U}_{t}\in\mathcal{B}\left(  \mathcal{H}\right)  $
such that $\tau_{t}=\operatorname{Ad}\left(  \bar{U}_{t}\right)  $. But if
$E=\operatorname*{supp}\psi$, the invariance implies $\tau_{t}\left(
E\right)  =E$, and hence $E\in\bar{U}_{t}^{\prime}$. Thus $U_{t}=E\bar{U}%
_{t}E=\bar{U}_{t}E=E\bar{U}_{t}$ is unitary. But $\bar{U}_{t}^{{}}S_{i}^{{}%
}\bar{U}_{t}^{\ast}=tS_{i}^{{}}$, so multiplying to the 
left 
with $E$,
we get
\[
U_{t}^{{}}V_{i}^{{}}U_{t}^{\ast}=tV_{i}^{{}}.
\]
Multiplying to the right with $U_{t}^{{}}V_{i}^{\ast}$, and summing over $i$, we
then
obtain%
\[
U_{t}=t\bs\left(  U_{t}\right)  ,
\]
i.e., $U_{t}$ is an eigenvector of $\bs$ with eigenvalue $\bar{t}$. Thus
$U_{t}^{\ast}$ is an eigenvector with eigenvalue $t$, and Lemma \ref{Lem6.5}
is proved.
\end{proof}

We next establish the converse inclusion.

\begin{lemma}
\label{Lem6.6}$\operatorname{PSp}\left(  \bs\right)  \cap\mathbb{T}%
\subset\left\{  t\in\mathbb{T}\mid\psi\circ\tau_{t}=\psi\right\}  $.
\end{lemma}

\begin{proof}
If $t\in\operatorname{PSp}\left(  \bs\right)  \cap\mathbb{T}$, let $U_{{}%
}^{\ast}=U_{t}^{\ast}$ be a corresponding eigenvector, and assume that
$\left\|  U^{\ast}\right\|  =1$. We argue that $U$ is unitary by using the
argument employed in the proof of Lemma \ref{Lem3.4}: By the generalized
Schwarz inequality,%
\[
\bs\left(  U^{\ast}U\right)  \geq\bs\left(  U^{\ast}\right)  \bs\left(
U\right)  =t\bar{t}U^{\ast}U=U^{\ast}U,
\]
so%
\[
\bs\left(  U^{\ast}U\right)  -U^{\ast}U\geq0.
\]
But by $\bs$-invariance of $\varphi$,%
\[
\varphi\left(  \bs\left(  U^{\ast}U\right)  -U^{\ast}U\right)  =0,
\]
and as $\varphi$ is faithful,%
\[
\bs\left(  U^{\ast}U\right)  =U^{\ast}U.
\]
Since $\bs$ is ergodic and $\left\|  U^{\ast}\right\|  =1$, it follows that%
\[
U^{\ast}U=\openone.
\]
In the same way, one shows
that
$UU^{\ast}=\openone$, so $U=U_{t}$ is unitary. But
we have%
\[
U_{t}=t\bs\left(  U_{t}\right)  =t\sum_{i}V_{i}^{{}}U_{t}^{{}}V_{i}^{\ast}.
\]

Before continuing the proof of Lemma \ref{Lem6.6}, we now prove

\begin{lemma}
\label{Lem6.7}If $U$ is a unitary operator in $\mathcal{B}\left(
\mathcal{K}\right)  $ with $\bs\left(  U\right)  =\bar{t}U$, where
$t\in\mathbb{T}\subset\mathbb{C}$, then%
\begin{equation}
UV_{i}U^{\ast}=tV_{i}
\label{eqNew6.1}
\end{equation}
for $i\in\mathbb{Z}_{d}$.
\end{lemma}

\begin{proof}
By \cite[Theorem 3.1]{Choi74}, we have%
\[
\bs\left(  XU\right)  =\bs\left(  X\right)  \bs\left(  U\right)  =\bs\left(
X\right)  \bar{t}U
\]
for all $X\in\mathcal{B}\left(  \mathcal{K}\right)  $. Define%
\[
X_{i}=UV_{i}U^{\ast}-tV_{i}.
\]
Then%
\begin{align*}
\sum_{i}X_{i}^{{}}X_{i}^{\ast}  &  =U\left(  \vphantom{\sum}\smash{\sum_{i}%
}V_{i}^{{}}V_{i}^{\ast}\right)  U_{{}}^{\ast}-t\sum_{i}V_{i}^{{}}UV_{i}^{\ast
}U_{{}}^{\ast}-\bar{t}\sum_{i}UV_{i}^{{}}U_{{}}^{\ast}V_{i}^{\ast}+\sum
_{i}V_{i}^{{}}V_{i}^{\ast}\\
&  =\openone-\bs\left(  U\right)  tU^{\ast}-\bar{t}U\bs\left(  U^{\ast
}\right)  +\openone\\
&  =\openone-\bs\left(  UU^{\ast}\right)  -\bs\left(  UU^{\ast}\right)
+\openone\\
&  =0.
\end{align*}
It follows that%
\[
X_{i}=UV_{i}U^{\ast}-tV_{i}=0,
\]
and Lemma \ref{Lem6.7} is proved.
\end{proof}

\proofindent\proofheadfont{Continuation of the proof of Lemma} \ref{Lem6.6}. 
We may now finalize
the proof of Lemma \ref{Lem6.6} by extending the unitary $U_{t}$ on
$\mathcal{K}$, to a unitary $\bar{U}_{t}$ on $\mathcal{H}$, through the
definition%
\[
\bar{U}_{t}\left(  \vphantom{\sum}\smash{\sum_{I}}\alpha_{I}S_{I}\xi
_{I}\right)  =\sum_{I}\alpha_{I}t^{\left|  I\right|  }S_{I}U_{t}%
\xi_{I}%
\]
where $I$ is a finite multi-index with elements from $\mathbb{Z}_{d}$,
$\alpha_{I}\in\mathbb{C}$ and $\xi_{I}\in\mathcal{K}$. $\bar{U}_{t}$ is well
defined and unitary by the following computation, where $J$, $I$ are
multi-indices related by $J=IJ^{\prime}$, where $J^{\prime}$ is another
multi-index. Lemma \ref{Lem6.7} is used in the computation.%
\begin{align*}
\ip{S_{I}U_{t}\xi_{I}}{S_{J}U_{t}\xi_{J}}
& =\ip{U_{t}\xi_{I}}{S_{J^{\prime}}U_{t}\xi_{J}}\\
& =\ip{U_{t}\xi_{I}}{V_{J^{\prime}}U_{t}\xi_{J}}\\
& =t^{\left|  J^{\prime}\right|  }\ip{U_{t}\xi_{I}}{U_{t}V_{J^{\prime}}\xi
_{J}}\\
& =t^{\left|  J^{\prime}\right|  }\ip{\xi_{I}}{V_{J^{\prime}}\xi_{J}}\\
& =\ip{t^{\vphantom{I}\smash{\left|  I\right| }}S_{I}\xi_{I}}{t^{\vphantom
{J}\smash{\left|  J\right| }}S_{J}\xi_{J}}.
\end{align*}%
But, from the definition of $\bar{U}_{t}$, it follows that%
\begin{equation}
\bar{U}_{t}^{{}}S_{i}^{{}}\bar{U}_{t}^{\ast}=tS_{i}^{{}},
\label{eqNew6.2}
\end{equation}
so $\bar{U}_{t}$ implements $\tau_{t}$.
(In passing from (\ref{eqNew6.1}) to (\ref{eqNew6.2})
with the lifting $U\mapsto \bar{U}_{t}$, we
note that this is a ``scaled''
version of the commutant
lifting in \chaptername{} \ref{Representations}.)
Use now the same symbol $\tau_{t}$ to denote also the
normal extension of $\tau_{t}$ to $\mathcal{B}\left(  \mathcal{H}\right)  $.
By construction of $\bar{U}_{t}$, we have $\bar{U}_{t}P=P\bar{U}_{t}$, so%
\[
\tau_{t}\left(  P\right)  =P.
\]
We now argue that $\psi\circ\tau_{t}=\psi$. Put $\psi_{t}=\psi\circ\tau_{t}$.
Since $\tau_{t}$ is unitarily implemented, $\psi_{t}$ is normal in the given
representation and extends to $\mathcal{B}\left(  \mathcal{H}\right)  $. Since
$\psi_{t}\left(  P\right)  =\psi\left(  \tau_{t}\left(  P\right)  \right)
=\psi\left(  P\right)  =\openone$, we have $\operatorname*{supp}\left(
\psi_{t}\right)  \leq P$, and we may define a state $\varphi_{t}$ on
$\mathcal{B}\left(  \mathcal{K}\right)  $ by%
\[
\varphi_{t}\left(  PXP\right)  =\psi_{t}\left(  X\right)
\]
for $X\in\mathcal{B}\left(  \mathcal{H}\right)  $. But%
\[
\lambda\tau_{t}\left(  X\right)  =\sum_{i}S_{i}^{{}}\tau_{t}^{{}}\left(
X\right)  S_{i}^{\ast}=\tau_{t}\lambda\left(  X\right)
\]
for $X\in\mathcal{B}\left(  \mathcal{H}\right)  $, and, as $\psi\circ
\lambda=\psi$, we deduce that%
\begin{align*}
\varphi_{t}\circ\bs\left(  PXP\right)   &  =\psi_{t}\circ\lambda\left(
X\right) \\
&  =\psi\tau_{t}\lambda\left(  X\right) \\
&  =\psi\lambda\tau_{t}\left(  X\right) \\
&  =\psi\tau_{t}\left(  X\right) \\
&  =\varphi\tau_{t}\left(  PXP\right) \\
&  =\varphi_{t}\left(  PXP\right)  ,
\end{align*}
so $\varphi_{t}\circ\bs=\varphi_{t}$.
Using the fact that
$\mathcal{B}\left(
\mathcal{K}\right)  $ has a unique $\bs$-invariant normal state by assumption,
we conclude that $\varphi_{t}=\varphi$, and hence%
\[
\psi\circ\tau_{t}=\psi_{t}=\psi.
\]
This ends the proof of Lemma \ref{Lem6.6}.
\end{proof}

We have now established that the sets (\ref{Thm6.3(1)}) and (\ref{Thm6.3(3)})
in Theorem \ref{Thm6.3} are equal. Clearly set (\ref{Thm6.3(1)}) is contained
in set (\ref{Thm6.3(2)}), and to establish the converse, we have to show that,
if $\psi$ is $\tau_{t}$-covariant for some $t\in\mathbb{T}$, then $\psi$ is
actually $\tau_{t}$-invariant. To this end, note that, as%
\[
\psi\circ\tau_{t}\circ\lambda=\psi\circ\lambda\circ\tau_{t}=\psi\circ\tau
_{t},
\]
this will follow once we can show the following lemma:

\begin{lemma}
\label{Lem6.8}Adopt the assumptions of Theorem \textup{\ref{Thm6.3}}. Then
$\psi$ is a unique $\lambda$-invariant normal state on $\mathcal{B}\left(
\mathcal{H}\right)  $.
\end{lemma}

\begin{proof}
If $X\in\mathcal{B}\left(  \mathcal{H}\right)  $, then%
\[
\wstarlim_{N\rightarrow\infty}\frac{1}{N+1}\sum_{k=0}^{N}\lambda^{k}\left(
X\right)  =\psi\left(  X\right)  \openone
\]
by the following reasoning: Putting%
\[
X_{N}=\frac{1}{N+1}\sum_{k=0}^{N}\lambda^{k}\left(  X\right) ,
\]
we have $\lambda\left(  X_{N}\right)  -X_{N}=\left(  \lambda^{N+1}\left(
X\right)  -X\right)  \diagup\left(  N+1\right)  $, and hence%
\[
\left\|  \lambda\left(  X_{N}\right)  -X_{N}\right\|  \leq\frac{2\left\|
X\right\|  }{N+1}.
\]
It follows that any weak*-limit point of the sequence $X_{N}$ is $\lambda
$-invariant. But, as the representation is irreducible, the only $\lambda
$-invariant elements in $\mathcal{B}\left(  \mathcal{H}\right)  $ are the
scalar multiples of $\openone$ (see, e.g., (3.5) in \cite{BJP96}). Moreover,
as $\psi\circ\lambda=\psi$, we have $\psi\left(  X_{N}\right)  =\psi\left(
X\right)  $, and the claim follows. Finally, if $\omega$ is a $\lambda
$-invariant normal state and $X\in\mathcal{B}\left(  \mathcal{H}\right)  $, it
follows that%
\[
\omega\left(  X\right)  =\omega\left(  X_{N}\right)  ,
\]
and therefore
\[
\omega\left(  X\right)  =\lim_{N\rightarrow\infty}\omega\left(  X_{N}\right)
=\omega\left(  \psi\left(  X\right)  \openone\right)  =\psi\left(  X\right)  .%
\settowidth{\qedskip}{$\displaystyle
\omega\left(   X \right)  =\lim_{N\rightarrow\infty}\omega\left(   X
_{N}\right)  =\omega\left(  \psi\left(   X \right)  \openone\right)
=\psi\left(   X \right)  .$}
\addtolength{\qedskip}{-\textwidth}
\rlap{\hbox to-0.5\qedskip{\hfil\qed}}%
\]%
\renewcommand{\qed}{}%
\end{proof}

Note that Lemma \ref{Lem6.8} could be used to simplify the last part of the
proof of Lemma \ref{Lem6.6}.

Next, we establish the finiteness of the three sets in Theorem \ref{Thm6.3}:

\begin{lemma}
\label{Lem6.9}$\left\{  t\in\mathbb{T}\mid\psi\circ\tau_{t}=\psi\right\}  $ is
a finite subgroup of $\mathbb{T}$.
\end{lemma}

\begin{proof}
The set is clearly a closed subgroup of $\mathbb{T}$, so if it is not finite
it is equal to $\mathbb{T}$. But in that case the automorphism group
$t\mapsto\tau_{t}$ extends to the weak closure $\pi_{\psi}\left(
\mathcal{O}_{d}\right)  ^{\prime\prime}$ of $\pi_{\psi}\left(  \mathcal{O}%
_{d}\right)  $ in the GNS representation defined by $\psi$. Since the
original representation of $\mathcal{O}_{d}$ on $\mathcal{H}$ is irreducible,
and $\psi$ is a normal state in this representation, $\pi_{\psi}\left(
\mathcal{O}_{d}\right)  ^{\prime\prime}$ is a type $\mathrm{I}$ factor and
$\pi_{\psi}$ extends canonically to a $\ast$-isomorphism from $\mathcal{B}%
\left(  \mathcal{H}\right)  =\mathcal{O}_{d}^{\prime\prime}$ to this factor.
Transporting $\tau_{t}$ back by this isomorphism, it follows that there exists
a unitary representation $t\mapsto U_{t}$ of $\mathbb{T}$ on $\mathcal{H}$
such that
\[
\tau_{t}^{{}}\left(  x\right)  =U_{t}^{{}}xU_{t}^{\ast}%
\]
for $x\in\mathcal{O}_{d}$.
For this covariant representation,
let%
\[
U_{t}=\sum_{n\in\mathbb{Z}}t^{n}E_{n}%
\]
be the Stone--Naimark--Ambrose--Godement (SNAG) decomposition \cite{Mac49} of
$U$. As%
\[
U_{t}S_{i}=\tau_{t}\left(  S_{i}\right)  U_{t}=tS_{i}U_{t},
\]
we obtain%
\[
E_{n}S_{i}=S_{i}E_{n-1},
\]
and thus%
\[
\lambda\left(  E_{n}\right)  =E_{n+1}.
\]
But%
\[
\psi\left(  \,\cdot\,\right)  =\psi\circ\int_{\mathbb{T}}\operatorname{Ad}%
U_{t}\left(  \,\cdot\,\right)  \,dt=\sum_{n}\psi\left(  E_{n}\,\cdot
\,E_{n}\right)  ,
\]
where we identify $\psi$ with the vector state it defines on the bounded
operators on the representation Hilbert space. Therefore, if%
\begin{align*}
\psi_{n}\left(  X\right)   &  =\psi\left(  E_{n}XE_{n}\right)  ,\\%
\intertext{then}%
\psi_{n}\left(  \lambda\left(  X\right)  \right)   &  =\psi\left(
E_{n}\lambda\left(  X\right)  E_{n}\right) \\
&  =\psi\left(  \lambda\left(  E_{n-1}XE_{n-1}\right)  \right) \\
&  =\psi_{n-1}\left(  X\right)  ,
\end{align*}
since $\psi\circ\lambda=\psi$. (Here we use implicitly the facts that both
$\psi$ and $\lambda$ extend by weak*-continuity to $\mathcal{B}\left(
\mathcal{H}\right)  =\mathcal{O}_{d}^{\prime\prime}$, and 
that
the invariance
$\psi\circ\lambda=\psi$ is preserved in the extension. If the original
representation were not irreducible, this point would be problematic.) But
then%
\[
\psi\left(  \openone\right)  =\sum_{n}\psi_{n}\left(  \openone\right)
=\sum_{n}\psi_{0}\left(  \openone\right)  =\infty,
\]
since $\lambda\left(  \openone\right)  =\openone$. This is impossible, so
Lemma \ref{Lem6.9} is established.
\end{proof}

It remains to prove the last statements of Theorem \ref{Thm6.3}. To this end,
define $k\in\mathbb{N}$ such that the finite group in Lemma \ref{Lem6.9} is%
\[
\left\{  \exp\left(  2\pi il\diagup k\right)  \mid l=0,1,\dots,k-1\right\}  .
\]

\begin{lemma}
\label{Lem6.10}With $k$ defined as above and $\psi$ as in Theorem
\textup{\ref{Thm6.3}}, we have%
\[
\pi_{\psi}\left(  \mathcal{O}_{d}^{\tau_{1/k}}\right)  ^{\prime\prime}%
=\pi_{\psi}\left(  \operatorname*{UHF}\nolimits_{d}\right)  ^{\prime\prime}%
\]
and%
\[
\pi_{\psi}\left(  \mathcal{O}_{d}\right)  ^{\prime\prime}\cap\pi_{\psi}\left(
\mathcal{O}_{d}^{\tau_{1/k}}\right)  ^{\prime}\cong\mathbb{C}^{k},
\]
where $\mathcal{O}_{d}^{\tau_{1/k}}$ denotes the fixed point algebra in
$\mathcal{O}_{d}$ under the gauge automorphism $\tau_{\exp\left(  2\pi
i\diagup k\right)  }$.
\end{lemma}

\begin{proof}
Since $\pi_{\psi}$ is merely a multiple of the given irreducible
representation on $\mathcal{H}$ (by normality of $\psi$), we only need to show%
\begin{align*}
\left(  \mathcal{O}_{d}^{\tau_{1/k}}\right)  ^{\prime\prime}  &  =\left(
\operatorname*{UHF}\nolimits_{d}\right)  ^{\prime\prime}\\%
\intertext{and}%
\left(  \mathcal{O}_{d}^{\tau_{1/k}}\right)  ^{\prime}  &  \cong\mathbb{C}%
^{k}.
\end{align*}
But, if $U$ is the unitary on $\mathcal{H}$ implementing $\tau_{\frac{1}{k}}$,
we have shown that $U$ is an eigenunitary of $\lambda$ with eigenvalue
$e^{-\frac{2\pi i}{k}}$. As $\lambda\left(  U^{k}\right)  =e^{-\frac{2\pi
ik}{k}}U^{k}=U^{k}$, we have $U^{k}\in\mathbb{C}%
\mkern2mu%
\openone$; and we may assume
$U^{k}=\openone$ by changing $U$ by a phase factor. Thus $U$
will
have a spectral
decomposition%
\[
U=\sum_{l\in\mathbb{Z}_{k}}e^{\frac{i2\pi l}{k}}E_{l}%
\]
where $E_{l}$, $l\in\mathbb{Z}_{k}$, are mutually orthogonal projections
summing up to $\openone$. Moreover, as%
\[
\lambda\left(  U\right)  =\sum_{l\in\mathbb{Z}_{k}}e^{\frac{i2\pi l}{k}%
}\lambda\left(  E_{l}\right)  =e^{-\frac{i2\pi}{k}}U=\sum_{l\in\mathbb{Z}_{k}%
}e^{\frac{i2\pi\left(  l-1\right)  }{k}}E_{l},
\]
we see that%
\[
\lambda\left(  E_{l}\right)  =E_{l+1}%
\]
for $l\in\mathbb{Z}_{k}$. It follows that all the
projections
$E_{l}$ are nonzero. Thus%
\begin{align*}
\left(  \mathcal{O}_{d}^{\tau_{1/k}}\right)  ^{\prime\prime}  &  =\left(
\mathcal{O}_{d}\right)  ^{\prime\prime}\cap\left\{  U,U^{\ast}\right\}
^{\prime}\\
&  =\left\{  E_{l}\mid l\in\mathbb{Z}_{k}\right\}  ^{\prime}\\
&  =\bigoplus_{l\in\mathbb{Z}_{k}}E_{l}\mathcal{B}\left(  \mathcal{H}\right)
E_{l},\\%
\intertext{so}%
\left(  \mathcal{O}^{\tau_{1/k}}\right)  ^{\prime}  &  =\sum_{l\in
\mathbb{Z}_{k}}\mathbb{C}%
\mkern2mu%
E_{l}\cong\mathbb{C}^{k}%
\end{align*}
as asserted.

To prove that%
\[
\left(  \operatorname*{UHF}\nolimits_{d}\right)  ^{\prime\prime}=\left(
\mathcal{O}^{\tau_{1/k}}\right)  ^{\prime\prime}%
\]
we first note that if $\tau_{t}$ is restricted to $t\in\left[  0,\frac{1}%
{k}\right\rangle $, then $\tau$ defines a representation of $\mathbb{T}$ in
$\operatorname{Aut}\left(  \mathcal{O}^{\tau_{1/k}}\right)  $.
Now consider
the
direct integral
representation%
\[
\pi=\int_{\left[  0,\frac{1}{k}\right\rangle }^{\oplus}dt\,\pi_{\psi
|_{\mathcal{O}_{d}^{\tau_{1/k}}}}\circ\tau_{t}%
\]
of $\mathcal{O}_{d}^{\tau_{1/k}}$ on $\mathcal{H}_{\psi|_{\mathcal{O}%
_{d}^{\tau_{1/k}}}}\otimes L^{2}\left(  \left[  0,\frac{1}{k}\right\rangle
\right)  $.
(See \cite{BEEK}.)

We establish the following observation
concerning this representation
before finalizing the proof of Lemma
\ref{Lem6.10}:

\begin{lemma}
\label{Lem6.11}$\openone\otimes L^{\infty}\left(  \left[  0,\frac{1}%
{k}\right\rangle \right)  \subset\pi\left(  \mathcal{O}_{d}^{\tau_{1/k}%
}\right)  ^{\prime\prime}$.
\end{lemma}

\begin{proof}
Note that%
\[
\lambda\left(  \mathcal{O}_{d}^{\tau_{1/k}}\right)  \subset\mathcal{O}%
_{d}^{\tau_{1/k}},
\]
since $\lambda\tau_{t}=\tau_{t}\lambda$. If $t_{1},t_{2}\in\left[  0,\frac
{1}{k}\right\rangle $ and $t_{1}\neq t_{2}$, it follows from the already
proved part of Theorem \ref{Thm6.3} that there exists an $x\in\mathcal{O}_{d}$
with%
\[
\psi\left(  \tau_{t_{1}}\left(  x\right)  \right)  \neq\psi\left(  \tau
_{t_{2}}\left(  x\right)  \right)  .
\]
Replacing $x$ with its mean over $\mathbb{Z}_{k}$,%
\[
\frac{1}{k}\sum_{l\in\mathbb{Z}_{k}}\tau_{\frac{l}{k}}\left(  x\right)  ,
\]
we may assume that $x\in\mathcal{O}_{d}^{\tau_{1/k}}$. 
Since
\[
\wlim_{N\rightarrow\infty}\frac{1}{N+1}\sum_{n=0}^{N}\lambda^{n}\left(
\tau_{t}\left(  x\right)  \right)  =\psi\left(  \tau_{t}\left(  x\right)
\right)  \openone
,
\]
by the reasoning in the proof of Lemma \ref{Lem6.8}, 
it follows that
\[
\wlim_{N\rightarrow\infty}\frac{1}{N+1}\sum_{n=0}^{N}\pi\left(  \lambda
^{n}\left(  x\right)  \right)  =\openone\otimes f,
\]
where%
\[
f\left(  t\right)  =\psi\left(  \tau_{t}\left(  x\right)  \right)  .
\]
Lemma \ref{Lem6.11} follows, as these $f$'s separate points.
\end{proof}

\proofindent\proofheadfont{Continuation of the proof of Lemma} \ref{Lem6.10}. 
It follows from
Lemma \ref{Lem6.11} that%
\[
\pi\left(  \mathcal{O}_{d}^{\tau_{1/k}}\right)  ^{\prime\prime}=\pi_{\psi
}\left(  \mathcal{O}_{d}^{\tau_{1/k}}\right)  ^{\prime\prime}\otimes
L^{\infty}\left(  \textstyle\left[  0,\frac{1}{k}\right\rangle \right)  .
\]
But $\pi$ is clearly $\tau$-covariant, $\mathbb{T}$ acting by translation, and
therefore%
\[
\pi\left(  \operatorname*{UHF}\nolimits_{d}\right)  ^{\prime\prime}=\pi\left(
\mathcal{O}_{d}^{\tau_{1/k}}\right)  ^{\prime\prime\,\tau}=\pi_{\psi}\left(
\mathcal{O}_{d}^{\tau_{1/k}}\right)  ^{\prime\prime}\otimes\openone.
\]
This equality then also holds on fibers, so%
\[
\pi_{\psi}\left(  \operatorname*{UHF}\nolimits_{d}\right)  ^{\prime\prime}%
=\pi_{\psi}\left(  \mathcal{O}_{d}^{\tau_{1/k}}\right)  ^{\prime\prime},
\]
and this ends the proof of Lemma \ref{Lem6.10}.
\end{proof}

\proofindent\proofheadfont{End of proof of Theorem} \ref{Thm6.3}. 
We finally observe from the
proof that this means that $\operatorname*{UHF}\nolimits_{d}$ acts irreducibly
on each of the subspaces $E_{l}\mathcal{H}$, that these representations are
mutually disjoint, and since $\lambda\left(  E_{l}\right)  =E_{l+1}$, the
endomorphism $\lambda$ maps these representations of $\operatorname*{UHF}%
\nolimits_{d}$ cyclically
one into another, i.e.,
\[
\pi_0\rightarrow \pi_1\rightarrow \dots \rightarrow \pi_{k-1} \rightarrow \pi_0,
\]
where $\pi_l$ is the ``cut down'' of $\pi_\psi$ by $E_l$,
$\pi_l\left(\,\cdot \,\right) =\pi_\psi\left(\,\cdot \,\right) E_l$,
$l\in \mathbb{Z}_k$.
\end{proof}

\chapter{\label{Quantum}Translationally invariant states on the two-sided
quantum chain}

Let us recall the definition of finitely correlated pure states from
\cite{FNW92,FNW94}. These are translationally invariant states defined on the
one-dimensional quantum chain $\bigotimes_{\mathbb{Z}}M_{d}$ as follows: Let
$\mathcal{K}$ be a finite-dimensional Hilbert space and let $V\colon
\mathcal{K}\rightarrow\mathcal{K}\otimes\mathbb{C}^{d}$ be an isometry. Define%
\[
\mathbb{E}\colon\mathcal{B}\left(  \mathcal{K}\right)  \otimes M_{d}%
\longrightarrow\mathcal{B}\left(  \mathcal{K}\right)
\]
by%
\[
\mathbb{E}\left(  X\right)  =V^{\ast}XV.
\]
Let $\varphi$ be a state on $\mathcal{B}\left(  \mathcal{K}\right)  $ such
that%
\[
\varphi\left(  \mathbb{E}\left(  B\otimes\openone\right)  \right)
=\varphi\left(  B\right)
\]
for all $B\in\mathcal{B}\left(  \mathcal{K}\right)  $. Define%
\[
\mathbb{E}_{A}\colon\mathcal{B}\left(  \mathcal{K}\right)  \longrightarrow
\mathcal{B}\left(  \mathcal{K}\right)
\]
by%
\[
B\longmapsto\mathbb{E}\left(  B\otimes A\right)
\]
for $A\in M_{d}$. Then%
\[
\omega\left(  A_{1}\otimes A_{2}\otimes\dots\otimes A_{m}\right)
=\varphi\left(  \mathbb{E}_{A_{1}}\circ\mathbb{E}_{A_{2}}\circ\dots
\circ\mathbb{E}_{A_{m}}\left(  \openone_{\mathcal{K}}\right)  \right)
\]
defines a translation-invariant state on $\bigotimes_{\mathbb{Z}}M_{d}$. It is
proved in \cite[Theorem 1.5]{FNW94} that this state is pure if the completely
positive map $\bs=\mathbb{E}_{\openone}$ has trivial peripheral spectrum,
i.e., the only eigenvectors of $\mathbb{E}_{\openone}$ with eigenvalue of
modulus one are the scalar multiples of $\openone$. Conversely, if $\omega$ is
pure, there does exist a realization of $\omega$ as above such that $\bs$ has
trivial peripheral spectrum (but it might not be the given one; see the
remarks at the end of \Chaptername{} \ref{Intro}). Now $V\colon\mathcal{K}%
\rightarrow\mathcal{K}\otimes\mathbb{C}^{d}=\bigoplus_{1}^{d}\mathcal{K}$ has
the matrix form $V=\left[
\begin{smallmatrix}
V_{1}^{\ast}\\%
\raisebox{0pt}[10pt]{$\vdots$}%
\\
V_{d}^{\ast}%
\end{smallmatrix}
\right]  $ and the property that $V$ is an isometry translates into%
\[
V^{\ast}V=\sum_{k\in\mathbb{Z}_{d}}V_{k}^{{}}V_{k}^{\ast}=\openone.
\]
We check that%
\[
\bs\left(  X\right)  =\mathbb{E}\left(  X\otimes\openone\right)  =\sum
_{k\in\mathbb{Z}_{d}}V_{k}^{{}}X V_{k}^{\ast}.
\]
Using Theorem \ref{Thm5.1} we can thus associate a representation $\pi_{V}$ of
$\mathcal{O}_{d}$ to $V$, and since $\varphi$ is normal by
finite-dimensionality of $\mathcal{K}$, we can associate a state $\psi$ on
$\mathcal{O}_{d}$ to $\left(  V,\varphi\right)  $ which is normal in the given
representation. We next verify that the restriction of $\psi$ to the one-sided
tensor product $\operatorname*{UHF}\nolimits_{d}=\bigotimes_{\mathbb{N}}M_{d}$
is equal to the restriction of $\omega$ to $\bigotimes_{\mathbb{N}}M_{d}$:%
\begin{align*}
\omega\left(  e_{i_{1}j_{1}}\otimes\dots\otimes e_{i_{m}j_{m}}\right)   &
=\varphi\left(  \mathbb{E}_{e_{i_{1}j_{1}}}\circ\mathbb{E}_{e_{i_{2}j_{2}}%
}\circ\dots\circ\mathbb{E}_{e_{i_{m}j_{m}}}\left(  \openone_{\mathcal{K}%
}\right)  \right) \\
&  =\varphi\left(  \mathbb{E}_{e_{i_{1}j_{1}}}\circ\mathbb{E}_{e_{i_{2}j_{2}}%
}\circ\dots\circ\mathbb{E}_{e_{i_{m-1}j_{m-1}}}\left(  V_{i_{m}}^{{}}V_{j_{m}%
}^{\ast}\right)  \right) \\
&  =\varphi\left(  V_{i_{1}}^{{}}\cdots V_{i_{m}}^{{}}V_{j_{m}}^{\ast}\cdots
V_{j_{1}}^{\ast}\right) \\
&  =\psi\left(  s_{i_{1}}^{{}}\cdots s_{i_{m}}^{{}}s_{j_{m}}^{\ast}\cdots
s_{j_{1}}^{\ast}\right) \\
&  =\psi\left(  e_{i_{1}j_{1}}^{\left(  1\right)  }\otimes\dots\otimes
e_{i_{m}j_{m}}^{\left(  m\right)  }\right)  .
\end{align*}
Note that finite-dimensionality of $\mathcal{K}$ and normality of $\varphi$ do
not play any role in the computation above.

The main theorem in this \chaptername{} is the following.

\begin{theorem}
\label{Thm7.1}Let $\mathcal{M}$ be a factor, $\varphi$ a faithful normal state
on $\mathcal{M}$,
$V_{1},\dots,V_{d}$ operators in $\mathcal{M}$ satisfying
\[
\sum_{k\in\mathbb{Z}_{d}}V_{k}^{{}}V_{k}^{\ast}=\openone,
\]
and $\bs$ the completely positive unital normal map of $\mathcal{B}\left(
\mathcal{K}\right)  $ defined by%
\[
\bs\left(  X\right)  =\sum_{k\in\mathbb{Z}_{d}}V_{k}^{{}}XV_{k}^{\ast}%
\]
for $X\in\mathcal{B}\left(  \mathcal{K}\right)  $, and assume that%
\[
\mathcal{B}\left(  \mathcal{K}\right)  ^{\bs}=\mathcal{M}^{\prime}.
\]
If $\mathcal{M}$ is type $\mathrm{I}$, 
the following two conditions are equivalent.

\begin{enumerate}
\item \label{Thm7.1(1)}The translationally invariant state $\omega$ defined by
$\left\{  \varphi,V_{1},\dots,V_{d}\right\}  $ on $\bigotimes_{\mathbb{Z}%
}M_{d}$ is pure.

\item \label{Thm7.1(2)}$\operatorname{PSp}\left(  \bs|_{\mathcal{M}}\right)
\cap\mathbb{T}=\left\{  1\right\}  $.
\end{enumerate}

\noindent If $\mathcal{M}$ is not assumed to be type $\mathrm{I}$, the
condition \textup{(\ref{Thm7.1(1)})} is nevertheless equivalent to each of the
following two conditions.

\begin{enumerate}
\setcounter{enumi}{2}

\item \label{Thm7.1(3)}$\omega$ is a factor state, i.e.,%
\[
\lim_{\left|  n\right|  \rightarrow\infty}\omega\left(  x\lambda^{n}\left(
y\right)  \right)  =\omega\left(  x\right)  \omega\left(  y\right)
\]
for all $x,y\in\bigotimes_{\mathbb{Z}}M_{d}$, where $\lambda$ is the shift.

\item \label{Thm7.1(4)}The Connes spectrum of $\tau|_{H}$ on $\pi_{\psi
}\left(  \mathcal{O}_{d}\right)  ^{\prime\prime}$ is $\hat{H}$, where $H$ is
the subgroup of $t\in\mathbb{T}$ such that $\tau_{t}$ extends from
$\mathcal{O}_{d}$ to the weak closure $\pi_{\psi}\left(  \mathcal{O}%
_{d}\right)  ^{\prime\prime}$ and $\psi$ is defined in Theorem
\textup{\ref{Thm6.3}}.
\end{enumerate}
\end{theorem}

\begin{remark}
\label{Rem7.2}The condition $\mathcal{B}\left(  \mathcal{K}\right)  ^{\bs
}=\mathcal{M}^{\prime}$ implies in particular that the von\ Neumann algebra
generated by $V_{1},\dots,V_{d}$ is $\mathcal{M}$. 
See \Chaptername{} \ref{Ergodic}
for a further discussion.

Note that the condition $\mathcal{B}\left(  \mathcal{K}\right)  ^{\bs
}=\mathcal{M}^{\prime}$ does not depend on the particular normal
representation of $\mathcal{M}$ \textup{(}when $\bs$ is defined by the
representatives for $V_{i}$\textup{)}. 
The reason for this is that
any normal representation of $\mathcal{M}$
is a product of a spatial isomor%
phism, an induction $\mathcal{M}\ni X\mapsto XP$
where $P$ is a projection in $\mathcal{M}^{\prime }$, and
an amplification $\mathcal{M}\ni X\mapsto X\otimes \openone$,
and applying these three types of maps
on the $V_{i}$'s, one verifies that
the condition remains the same;
see
\cite[Th\'eor\`eme I.4.3]{Dix96}
for details on normal representations.
When developing a duality theory later,
we will use the representation where $\varphi$ is defined by a separating and
cyclic vector.
\end{remark}

The rest of this \chaptername{} 
will be devoted to a proof of Theorem \ref{Thm7.1}.
To this end we have to develop a certain duality theory for the objects
$\left(  \mathcal{M}\text{,}\varphi,V_{1},\dots,V_{d}\right)  $. But
before that we will mention
some more pedestrian results on translationally invariant states.

Recall from \cite[Example 4.3.24]{BrRoI} that any translationally invariant
factor state of $\bigotimes_{\mathbb{Z}}M_{d}$ is extremal among the invariant
states, i.e., is ergodic. Conversely, an ergodic state need not be a factor
state: If, for example, $\omega_{1}$, $\omega_{2}$ are distinct pure states on
$M_{d}$, the mean of the pure product state
\[
\cdots\otimes\omega_{1}\otimes\omega_{2}\otimes\omega_{1}\otimes\omega
_{2}\otimes\cdots\text{\quad on }\bigotimes_{\mathbb{Z}}M_{d}%
\]
and its shift is extremally invariant, but not a factor state (see
\cite[Example 4.3.26]{BrRoI}). 

The difference between factor states
and ergodic states is reflected in the
fact that if $\omega$ is a translationally in%
variant state on $\bigotimes_\mathbb{Z}M_d$, then
$\omega$ is a factor state if and only
if it is strongly clustering,
\[
\lim_{\left| n\right| \rightarrow \infty }
\omega \left( x\lambda^{n}\left( y\right) \right) 
=\omega \left( x\right) \omega \left( y\right) 
\]
(see \cite{Pow67}),
while $\omega$ is
ergodic if and only if it is clustering
in the mean
\[
\lim_{N\rightarrow \infty }
\frac{1}{N+1}\sum_{n=0}^{N}
\omega \left( x\lambda^{n}\left( y\right) \right) 
=\omega \left( x\right) \omega \left( y\right) 
\]
(see \cite[Example 4.3.5 and Theorem 4.3.17]{BrRoI}).
However, the following is true:

\begin{proposition}
\label{Pro7.3}There is a canonical one-one correspondence between the
following three sets.

\begin{enumerate}
\item \label{Pro7.3(1)}The set of extremal translationally invariant states on
$\bigotimes_{\mathbb{Z}}M_{d}$.

\item \label{Pro7.3(2)}The set of states on $\bigotimes_{\mathbb{N}}M_{d}$
which are extremal among the states invariant under the one-sided shift
$\lambda$.

\item \label{Pro7.3(3)}The set of orbits under the gauge action $\tau$ in the
sets of states $\psi$ on $\mathcal{O}_{d}$ such that%
\[
\psi\circ\lambda=\psi
\]
and%
\[
\psi\text{ is a factor state with an ergodic restriction to }%
\operatorname*{UHF}\nolimits_{d}.
\]
\end{enumerate}

The maps giving the correspondence are defined by the restriction maps from
Set \textup{(\ref{Thm7.1(1)})} and Set \textup{(\ref{Thm7.1(3)})} to Set
\textup{(\ref{Thm7.1(2)})}, using the inclusions%
\begin{align*}
\bigotimes_{\mathbb{N}}M_{d}  &  \subset\bigotimes_{\mathbb{Z}}M_{d}\\%
\intertext{and}%
\bigotimes_{\mathbb{N}}M_{d}=\operatorname*{UHF}\nolimits_{d}  &
\subset\mathcal{O}_{d}.
\end{align*}
\end{proposition}

\begin{proof}
\renewcommand{\qed}{} If we also use $\lambda$ to denote the two-sided shift
on $\bigotimes_{\mathbb{Z}}M_{d}$, the new $\lambda$ extends the old, and%
\[
\bigotimes_{\mathbb{Z}}M_{d}=\overline{\bigcup_{n=1}^{\infty}\lambda
^{-n}\left(
\vphantom{\bigotimes}\smash{\bigotimes
_{\mathbb{N}}}%
M_{d}\right)  }^{\left\|  \,\cdot\,\right\|  },
\]
so the one-one correspondence between Set \textup{(\ref{Thm7.1(1)})} and Set
\textup{(\ref{Thm7.1(2)})} is trivial.

It is clear that the map from Set \textup{(\ref{Thm7.1(3)})} to Set
\textup{(\ref{Thm7.1(2)})} is well defined. To prove that it is injective, let
$\omega^{\prime}$ be an extremal invariant state on $\bigotimes_{\mathbb{N}%
}M_{d}$, and consider the set%
\[
K=\left\{  \psi\mid\psi\text{ is a state of }\mathcal{O}_{d}\text{ such that
}\psi\circ\lambda=\psi\text{ and }\psi|_{\operatorname*{UHF}\nolimits_{d}%
}=\omega^{\prime}\right\}  .
\]
By applying an invariant mean on an extension of $\omega^{\prime}$ to
$\mathcal{O}_{d}$ it is clear that $K$ is nonempty, and $K$ is clearly convex
and compact, and a face in the set of $\lambda$-invariant states since
$\omega^{\prime}$ is extremal. We finish the proof of Proposition \ref{Pro7.3}
by proving:
\end{proof}

\begin{lemma}
\label{Lem7.4}$\psi\in K$ is an extremal point in $K$ if and only if $\psi$ is
a factor state, and then all other extremal points have the form $\psi
\circ\tau_{t}$ for some $t\in\mathbb{T}$.
\end{lemma}

\begin{proof}
If $\psi$ is not factorial, there is a nontrivial projection $E\in\pi_{\psi
}\left(  \mathcal{O}_{d}\right)  ^{\prime\prime}\cap\pi_{\psi}\left(
\mathcal{O}_{d}\right)  ^{\prime}$. But%
\[
\psi_{E}\left(  x\right)  =\ip{E\Omega_{\psi}}{\pi_{\psi}\left( x\right
) E\Omega_{\psi}}.
\]
Then $\psi_{E}\leq\psi$ and as $\lambda\left(  E\right)  =E$ we have $\psi
_{E}\circ\lambda=\psi_{E}$. But%
\[
\psi_{E}|_{\operatorname*{UHF}\nolimits_{d}}\leq\omega^{\prime}%
\]
and it follows from extremality of $\omega^{\prime}$ that there is a scalar
$c$ such that%
\[
\psi_{E}|_{\operatorname*{UHF}\nolimits_{d}}=c%
\mkern2mu%
\omega^{\prime}.
\]
But as $\psi=\psi_{E}+\psi_{\openone-E}$, and $\psi_{E}$ and $\psi
_{\openone-E}$ are disjoint, this contradicts the extremality of $\psi$.
Thus the extremal points in $K$ are factor states.

Conversely, if $\psi\in K$ is a factor state, it follows as in the proof of
Lemma \ref{Lem6.8} that%
\[
\wlim_{N\rightarrow\infty}\frac{1}{N+1}\sum_{k=0}^{N}\pi_{\psi}\left(
\lambda^{k}\left(  x\right)  \right)  =\psi\left(  x\right)  \openone,
\]
and hence $\psi$ is ergodic by \cite[Theorems 4.3.17 and 4.3.23]{BrRoI}.
(Strictly speaking, these theorems are proved under the assumption that
$\lambda$ is an automorphism, but extending $\lambda$ to an automorphism of
the inductive
limit%
\[
\mathcal{O}_{d}\overset{\lambda}{\longrightarrow}\mathcal{O}_{d}%
\overset{\lambda}{\longrightarrow}\mathcal{O}_{d}\overset{\lambda
}{\longrightarrow}\cdots
\]
and extending $\psi$ by requiring $\lambda$-invariance, one still has the
clustering%
\[
\lim_{N\rightarrow\infty}\psi\left(  y\left(  \frac{1}{N+1}%
\vphantom{\sum}\smash{\sum_{k=0}^{N}}%
\lambda^{k}\left(  x\right)  \right)  z\right)  =\psi\left(  yz\right)
\psi\left(  x\right)  ,
\]
so the extended $\psi$ is ergodic, and thus the original $\psi$ is so,
since there
is a one-one correspondence between the $\lambda$-invariant states on
$\mathcal{O}_{d}$ and those on the inductive
limit.)

Finally, let $\psi$ be a given extremal point in the face $K$ in the invariant
states. It follows from \cite[Theorem 4.3.19]{BrRoI} and the previous
paragraph that any two translates $\psi\circ\tau_{t_{1}}$, $\psi\circ
\tau_{t_{2}}$ of $\psi$ are either equal or disjoint. Put%
\[
G=\left\{  t\in\mathbb{T}\mid\psi\circ\tau_{t}=\psi\right\}
\]
and define
\[
\psi_{0}=\omega^{\prime}\circ\int_{\mathbb{T}}\tau_{t}\,dt=\int_{\mathbb{T}%
}\left(  \psi\circ\tau_{t}\right)  \,dt.
\]
Then%
\[
\pi_{\psi_{0}}=\int_{\mathbb{T}\diagup G}^{\oplus}\left(  \pi_{\psi}\circ
\tau_{t^{\prime}}\right)  \,dt^{\prime}%
\]
is the central decomposition of $\pi_{\psi_{0}}$
by Lemma \ref{Lem6.11} and its proof%
. If now $\psi^{\prime}$ is an
extremal point in $K$, i.e., $\psi^{\prime}$ is a factorial $\lambda
$-invariant state with $\psi^{\prime}|_{\operatorname*{UHF}\nolimits_{d}%
}=\omega^{\prime}$, then%
\[
\frac{1}{\varepsilon}\int_{0}^{\varepsilon}\left(  \psi^{\prime}\circ\tau
_{t}\right)  \,dt\leq\frac{1}{\varepsilon}\int_{\mathbb{T}}\left(
\psi^{\prime}\circ\tau_{t}\right)  \,dt=\frac{1}{\varepsilon}\psi_{0}%
\]
and hence, since the left-hand side of the above inequality is $\lambda
$-invariant, by Segal's Radon--Nikodym theorem \cite[Theorem 2.3.19]{BrRoI},
there is a function $g_{\varepsilon}\in L^{\infty}\left(  \mathbb{T}\diagup
G\right)  $ such that%
\[
\frac{1}{\varepsilon}\int_{0}^{\varepsilon}\left(  \psi^{\prime}\circ\tau
_{t}\right)  \,dt=\int_{\mathbb{T}\diagup G}\left(  \psi\circ\tau_{t}\right)
g_{\varepsilon}\left(  t\right)  \,dt.
\]
Letting $\varepsilon\rightarrow0$, we find a measure $\mu$ on $\mathbb{T}%
\diagup G$ such that%
\[
\psi^{\prime}=\int_{\mathbb{T}\diagup G}\left(  \psi\circ\tau_{t}\right)
\,d\mu\left(  t\right)  .
\]
But as $\psi^{\prime}$ is extremal in $K$, this must be a Dirac measure, and
$\psi^{\prime}=\psi\circ\tau_{t}$ for some $t$.

This ends the proof of Lemma \ref{Lem7.4} and Proposition \ref{Pro7.3}.
\end{proof}

In order to prove Theorem \ref{Thm7.1}, we need to develop a duality theory
for the objects $\left(  \mathcal{M},\varphi,V_{1},\dots,V_{d},\bs\right)
$ somewhat different from the duality theory in \cite{Jor96a}. The starting
point is more restrictive in that the normal state $\varphi$ is assumed to be
faithful. We assume that $V_{1},\dots,V_{d}\in\mathcal{M}$ and%
\[
\sum_{j=1}^{d}V_{j}^{{}}V_{j}^{\ast}=\openone
,
\]
and assume invariance%
\[
\varphi\circ\bs=\varphi ,
\]
where $\bs$ is the unital completely positive map on $\mathcal{M}$ defined
by%
\[
\bs\left(  X\right)  =\sum_{j=1}^{d}V_{j}^{{}}X V_{j}^{\ast}.
\]
We will construct a dual object $\vphantom{\widetilde{\mathcal{M}},
\tilde{\varphi},\tilde{V}_{1},\dots,\tilde{V}_{d},\tilde{\bs}}\left(
\vphantom{\mathcal{M}V_{1}V_{d}}\smash{\widetilde{\mathcal{M}},
\tilde{\varphi},\tilde{V}_{1},\dots,\tilde{V}_{d},\tilde{\bs}}\right)  $
satisfying the same axioms. To this end we assume $\mathcal{M}$ is acting on a
Hilbert space $\mathcal{K}$ with a cyclic vector $\Phi$ such that%
\[
\varphi\left(  X\right)  =\ip{\Phi}{ X \Phi}.
\]
Note that $\Phi$ is then separating for $\mathcal{M}$ by faithfulness of
$\varphi$. In the application to Theorem \ref{Thm7.1}, the system $\left(
\mathcal{M},\varphi,V_{1},\dots,V_{d},\bs\right)  $ will roughly correspond
to a state on a Cuntz algebra $\mathcal{O}_{d}$ with associated
$\operatorname*{UHF}$ algebra $\bigotimes_{1}^{\infty}M_{d}$, and the dual
object to a state on an isomorphic Cuntz algebra $\widetilde{\mathcal{O}}_{d}$
with associated $\operatorname*{UHF}$ algebra $\bigotimes_{-\infty}^{0}M_{d}$,
and $\bigotimes_{1}^{\infty}M_{d}$ and $\bigotimes_{-\infty}^{0}M_{d}$ will be
embedded into $\bigotimes_{-\infty}^{\infty}M_{d}$ in the obvious manner.
This statement will be made more precise in
Lemma \ref{LemNew7.15}.

Then to the definitions: We let $\widetilde{\mathcal{M}}=\mathcal{M}^{\prime}%
$. The cyclic and separating vector $\Phi$ for $\mathcal{M}$ defines 
the associated Tomita
modular
conjugation $J$, and modular operator $\Delta$;
see, e.g.,
\cite[Theorem 2.5.14]{BrRoI}. Let $\sigma_{t}$
be the modular automorphism group of $\mathcal{M}$:%
\[
\sigma_{t}\left(  X\right)  =\Delta^{it}X\Delta^{-it}%
\]
for $X\in\mathcal{M}$. We put%
\begin{align*}
\tilde{V}_{j}^{{}}  &  =\left(  J\sigma_{\frac{i}{2}}\left(  V_{j}^{\ast
}\right)  J\right)  \overline{\quad\mathstrut}\\
&  =\left(  J\Delta^{-\frac{1}{2}}V_{j}^{\ast}\Delta^{\frac{1}{2}}J\right)
\overline{\quad\mathstrut}%
\end{align*}
where $\overline{\quad\mathstrut}$ denotes closure of the respective
operators. To show that this is a well-defined operator in $\mathcal{M}%
^{\prime}$, define a positive sesquilinear form $Q_{j}$ on $J\mathfrak{A}\Phi
$, where $\mathfrak{A}$ is the $\ast$-algebra of $\sigma_{t}$-entire elements
in $\mathcal{M}$, by%
\begin{align*}
Q_{j}\left(  JX\Phi,JY\Phi\right)   &  =\vphantom{\tilde{V}_{j}^{\ast}}%
\ip{\vphantom{V_{j}^{\ast}}\smash{\tilde{V}_{j}^{\ast}}J X \Phi}%
{\vphantom{V_{j}^{\ast}}\smash{\tilde{V}_{j}^{\ast}}J Y \Phi}\\
&  =\ip{J\sigma_{-\frac{i}{2}}\left( V_{j}\right) X \Phi}{J\sigma_{-\frac
{i}{2}}\left( V_{j}\right) Y \Phi}\\
&  =\ip{J\Delta^{\frac{1}{2}}V_{j}\sigma_{\frac{i}{2}}\left( X \right) \Phi
}{J\Delta^{\frac{1}{2}}V_{j}\sigma_{\frac{i}{2}}\left( Y \right) \Phi}\\
&  =\ip{\sigma_{-\frac{i}{2}}\left( X _{}^{\ast}\right) V_{j}^{\ast}\Phi
}{\sigma_{-\frac{i}{2}}\left( Y _{}^{\ast}\right) V_{j}^{\ast}\Phi}.
\end{align*}
Hence%
\begin{align*}
\sum_{j\in\mathbb{Z}_{d}}Q_{j}\left(  JX\Phi,JY\Phi\right)   &  =\varphi
\left(  \bs\left(  \sigma_{\frac{i}{2}}\left(  X\right)  \sigma_{-\frac
{i}{2}}\left(  Y^{\ast}\right)  \right)  \right) \\
&  =\varphi\left(  \sigma_{\frac{i}{2}}\left(  X\right)  \sigma_{-\frac{i}{2}%
}\left(  Y^{\ast}\right)  \right) \\
&  =\ip{\sigma_{-\frac{i}{2}}\left( X _{}^{\ast}\right) \Phi}{\sigma
_{-\frac{i}{2}}\left( Y _{}^{\ast}\right) \Phi}\\
&  =\ip{\Delta^{\frac{1}{2}} X ^{\ast}\Phi}{\Delta^{\frac{1}{2}} Y ^{\ast}%
\Phi}\\
&  =\ip{J\Delta^{\frac{1}{2}} Y ^{\ast}\Phi}{J\Delta^{\frac{1}{2}} X ^{\ast
}\Phi}\\
&  =\ip{ Y \Phi}{ X \Phi}\\
&  =\ip{J X \Phi}{J Y \Phi}.
\end{align*}
It follows both that%
\[
\vphantom{\tilde{V}_{j}^{\ast}}\left\|  \vphantom{V_{j}^{\ast}}\smash
{\tilde{V}_{j}^{\ast}}JX\Phi\right\|  \leq\left\|  JX\Phi\right\|  ,
\]
i.e., $\tilde{V}_{j}^{\ast}$ is bounded, and that
\[
\sum_{j=1}^{d}\tilde{V}_{j}^{{}}\tilde{V}_{j}^{\ast}=\openone.
\]
We now naturally define a completely positive map $\tilde{\bs}$ on
$\widetilde{\mathcal{M}}=\mathcal{M}^{\prime}$ by%
\[
\tilde{\bs}\left(  X\right)  =\sum_{j\in\mathbb{Z}_{d}}\tilde{V}_{j}^{{}}X
\tilde{V}_{j}^{\ast}%
\]
for $X\in\widetilde{\mathcal{M}}$, and a faithful normal state $\tilde
{\varphi}$ on $\widetilde{\mathcal{M}}$ by%
\[
\tilde{\varphi}\left(  X\right)  =\ip{\Phi}{ X \Phi}%
\]
for $X\in\mathcal{M}^{\prime}$. We introduce the terminology $\vphantom
{\tilde{\varphi},\tilde{V}_{1},\dots,\tilde{V}_{d},\tilde{\bs}}\left(
\mathcal{M}^{\prime},\vphantom{V_{1}V_{d}}\smash{\tilde{\varphi},\tilde{V}%
_{1},\dots,\tilde{V}_{d},\tilde{\bs}}\right)  $ for the \emph{dual system}
of $\left(  \mathcal{M},\varphi,V_{1},\dots,V_{d},\bs\right)  $. Note that%
\[
\tilde{\varphi}\circ\tilde{\bs}=\tilde{\varphi},
\]
since, for $X\in\mathcal{M}^{\prime}$,%
\begin{align*}
\tilde{\varphi}\tilde{\bs}\left(  X\right)   &  =\sum_{j}\ip{\Phi}%
{J\sigma_{\frac{i}{2}}\left( V_{j}^{\ast}\right) J X J\sigma_{-\frac{i}{2}%
}\left( V_{j}^{{}}\right) \Phi}\\
&  =\sum_{j}\ip{V_{j}^{\ast}\Phi}{ X V_{j}^{\ast}\Phi}\\
&  =\vphantom{\sum_{j}} \ip{\Phi}{ X \left( \vphantom{\sum}\smash{\sum_{j}%
}V_{j}^{}V_{j}^{\ast}\right) \Phi}\\
&  =\tilde{\varphi}\left(  X\right)  .
\end{align*}
The term ``dual system'' is justified by the fact that the dual system
of
\linebreak 
$\vphantom{\tilde{\varphi},\tilde{V}_{1},\dots,\tilde
{V}_{d},\tilde{\bs}}\left(  \mathcal{M}^{\prime},\vphantom{V_{1}V_{d}%
}\smash{\tilde{\varphi},\tilde{V}_{1},\dots,\tilde
{V}_{d},\tilde{\bs}}\right)  $ is $\left(  \mathcal{M},\varphi,V_{1}%
,\dots,V_{d},\bs\right)  $ again. For this, we just need to check
$\Tilde{\Tilde{V}}_{j}=V_{j}$. But this follows from the computation%
\begin{align*}
\Tilde{\Tilde{V}}_{j}  
&  =J\Delta^{\frac{1}{2}}\tilde{V}_{j}^{\ast}\Delta^{-\frac
{1}{2}}J\\
&  =J\Delta^{\frac{1}{2}}\left(  J\Delta^{\frac{1}{2}}V_{j}^{{}}\Delta
^{-\frac{1}{2}}J\right)  \Delta^{-\frac{1}{2}}J\\
&  =V_{j}%
\end{align*}
where we used that $J$ and $\Delta^{-1}$ are the modular conjugation and
modular operator associated to the pair $\left(  \mathcal{M}^{\prime}%
,\Phi\right)  $, $J\Delta=\Delta^{-1}J$ and $J^{2}=\openone$.

This duality has several nice properties. For example $\bs$ is ergodic if
and only if $\tilde{\bs}$ is, and $\operatorname{PSp}\left(  \bs\right)
\cap\mathbb{T}=\operatorname{PSp}\left(  \tilde{\bs}\right)  \cap
\mathbb{T}$. These properties will be discussed in \Chaptername{} 
\ref{Hilbert}. For
the moment we return to the proof of Theorem \ref{Thm7.1}. So let $\left(
\mathcal{M},\varphi,V_{1},\dots,V_{d},\bs\right)  $ be as in the hypothesis
of the theorem, put $\mathcal{K}=\mathcal{H}_{\varphi}$ and identify $V_{i}$
with its representative $\pi_{\omega}\left(  V_{i}\right)  $ on $\mathcal{K}$.
If $\vphantom{\widetilde{\mathcal{M}},\tilde{\varphi},\tilde{V}_{1}%
,\dots,\tilde
{V}_{d},\tilde{\bs}}\left(  \vphantom{\mathcal{M}V_{1}V_{d}}\smash
{\widetilde{\mathcal{M}},\tilde{\varphi},\tilde{V}_{1},\dots,\tilde
{V}_{d},\tilde{\bs}}\right)  $ is the dual system, we have the canonical
identification $\mathcal{H}_{\tilde{\varphi}}=\mathcal{K}$, and $\tilde
{\varphi}$ is the vector state on $\widetilde{\mathcal{M}}$ defined by the
same vector $\Phi$ as $\varphi$. By Theorem \ref{Thm5.1} there are Hilbert
spaces $\mathcal{H}_{0}$, $\widetilde{\mathcal{H}}_{0}$ containing
$\mathcal{K}$, with projectors $P_{0}\colon\mathcal{H}_{0}\rightarrow
\mathcal{K}$, $\tilde{P}_{0}\colon\widetilde{\mathcal{H}}_{0}\rightarrow
\mathcal{K}$ and representations $S_{i}$, $\tilde{S}_{i}$ of the Cuntz
relations on $\mathcal{H}_{0}$, $\widetilde{\mathcal{H}}_{0}$, respectively
such that $\mathcal{K}$ is cyclic for both representations and%
\begin{align*}
P_{0}^{{}}S_{I}^{{}}S_{J}^{\ast}P_{0}^{{}}  &  =V_{I}^{{}}V_{J}^{\ast},\\
\tilde{P}_{0}^{{}}\tilde{S}_{I}^{{}}\tilde{S}_{J}^{\ast}\tilde{P}_{0}^{{}}  &
=\tilde{V}_{I}^{{}}\tilde{V}_{J}^{\ast}.
\end{align*}
We will now form a sort of amalgamated tensor product of $\mathcal{H}_{0}$ and
$\widetilde{\mathcal{H}}_{0}$ over the joint subspace $\mathcal{K}$ and thus
obtain a Hilbert space $\mathcal{H}$ carrying two commuting representations of
$\mathcal{O}_{d}$. To this end we generalize the construction in Remark
\ref{RemNew5.2}. $\mathcal{H}$ is the completion of the quotient of
\[
H=\mathbb{C}%
\mkern2mu%
\mathcal{I}\otimes\mathbb{C}%
\mkern2mu%
\widetilde{\mathcal{I}}\otimes\mathcal{K},
\]
where $\mathcal{I}$, $\widetilde{\mathcal{I}}$ both consist of all finite
sequences in $\mathbb{Z}_{d}$, by the equivalence relation defined by a
semi-inner product defined on $H$ by requiring%
\begin{align*}
\vphantom{\tilde{I}\tilde{J}\tilde{V}_{\tilde{J}}}\ip{I\otimes\smash{\tilde
{I}}\otimes\xi}{IJ\otimes\smash{\tilde{I}\tilde{J}}\otimes\eta}  &  =\ip{\xi
}{V_{J}\smash{\tilde{V}_{\tilde{J}}}\eta},\\
\vphantom{\tilde{I}\tilde{J}\tilde{V}_{\tilde{J}}}\ip{I\otimes\smash{\tilde
{I}\tilde{J}}\otimes\xi}{IJ\otimes\smash{\tilde{I}}\otimes\eta}  &
=\ip{\smash{\tilde{V}_{\tilde{J}}}\xi}{V_{J}\eta},
\end{align*}
etc., all inner products that cannot be put in these forms being zero. Since
the $V_{J}$'s and $\tilde{V}_{\tilde{J}}$'s commute along with all
combinations of their adjoints, we see that this gives rise to two commuting
representations of $\mathcal{O}_{d}$ on $\mathcal{H}_{0}$ as follows:%
\begin{align}
\vphantom{\tilde{J}}S_{I}\Lambda\left(  J\otimes\smash{\tilde{J}}\otimes
\xi\right)   &  =\Lambda\left(  IJ\otimes\smash{\tilde{J}}\otimes\xi\right)
,\label{eqNew7.1}  \\
\vphantom{\tilde{I}\tilde{J}}\tilde{S}_{\tilde{I}}\Lambda\left(
J\otimes\smash{\tilde{J}}\otimes\xi\right)   &  =\Lambda\left(  J\otimes
\smash{\tilde{I}\tilde{J}}\otimes\xi\right)  ,\label{eqNew7.2}  
\end{align}
where $\Lambda\colon H\rightarrow\mathcal{H}$ is the quotient map. This is a
slight abuse of notation as the earlier $S_{I}$, $\tilde{S}_{\tilde{I}}$
identify with the restriction of the present $S_{I}$, $\tilde{S}_{\tilde{I}}$
to the subspaces $\mathcal{H}_{0}$, $\widetilde{\mathcal{H}}_{0}$ of
$\mathcal{H}$ spanned by vectors $\Lambda\left(  I\otimes\left\{
\varnothing\right\}  \otimes\xi\right)  $ and $\vphantom{\tilde{I}}%
\Lambda\left(  \left\{  \varnothing\right\}  \otimes\vphantom{I}\smash
{\tilde{I}}\otimes\xi\right)  $, respectively.

All the previous statements are easy to check. For example the positivity of
the sesquilinear form on $H\times H$ is checked by induction as follows, where
the operators $S_{i}$ and $\tilde{S}_{i}$ on $H$ are defined in the obvious
manner: If%
\[
\zeta=\sum_{ij}S_{i}\tilde{S}_{j}\zeta_{ij}+\sum_{i}S_{i}\zeta_{i\,0}+\sum
_{j}\tilde{S}_{j}\zeta_{0\,j}+\zeta_{0}%
\]
is a general element in%
\[
H_{n}=\mathbb{C}%
\mkern2mu%
\mathcal{I}_{n}\otimes\mathbb{C}%
\mkern2mu%
\widetilde{\mathcal{I}}_{n}\otimes\mathcal{K}%
\]
where%
\begin{align*}
\zeta_{ij}  &  \in H_{n-1},\\
\zeta_{i\,0}  &  \in\mathbb{C}%
\mkern2mu%
\mathcal{I}_{n-1}\otimes\left\{  \varnothing\right\}  \otimes\mathcal{K},\\
\zeta_{0\,j}  &  \in\left\{  \varnothing\right\}  \otimes\mathbb{C}%
\mkern2mu%
\widetilde{\mathcal{I}}_{n-1}\otimes\mathcal{K},\\
\zeta_{0}  &  \in\mathcal{K},
\end{align*}
and we assume the form is positive on $H_{n-1}\times H_{n-1}$, we compute%
\begin{multline*}
\ip{\zeta}{\zeta}=\sum_{ij}\left\|  \zeta_{ij}\right\|  ^{2}+\sum_{i}\left\|
\zeta_{i\,0}\right\|  ^{2}+\sum_{j}\left\|  \zeta_{0\,j}\right\|  ^{2}+\left\|
\zeta_{0}\right\|  ^{2}\\
+\sum_{ij}\left\{\ip{\vphantom{V}\smash{\tilde{V}}_{j}\zeta_{ij}}{\zeta_{i\,0}}
+\ip{\zeta_{i\,0}}{\vphantom{V}\smash{\tilde{V}}_{j}\zeta_{ij}}\right\}
+\sum_{ij}\left\{  \ip{V_{i}\zeta_{ij}}{\zeta_{0\,j}}
+\ip{\zeta_{0\,j}}{V_{i}\zeta_{ij}}\right\}  \\
+\sum_{ij}\left\{  
\ip{V_{i}\vphantom{V}\smash{\tilde{V}}_{j}\zeta_{ij}}{\zeta_{0}}
+\ip{\zeta_{0}}{V_{i}\vphantom{V}\smash{\tilde{V}}_{j}\zeta_{ij}}\right\}
+\sum_{i}\left\{  \ip{V_{i}\zeta_{i\,0}}{\zeta_{0}}
+\ip{\zeta_{0}}{V_{i}\zeta_{i\,0}}\right\}  \\
+\sum_{j}\left\{  \ip{\vphantom{V}\smash{\tilde{V}}_{j}\zeta_{0\,j}}{\zeta_{0}}
+\ip{\zeta_{0}}{\vphantom{V}\smash{\tilde{V}}_{j}\zeta_{0\,j}}\right\}
+\sum_{ij}\left\{  
\ip{V_{i}\zeta_{i\,0}}{\vphantom{V}\smash{\tilde{V}}_{j}\zeta_{0\,j}}
+\ip{\vphantom{V}\smash{\tilde{V}}_{j}\zeta_{0\,j}}
{V_{i}\zeta_{i\,0}}\right\}  \\
=\sum_{ij}\left\|  \zeta_{ij}^{{}}
+\vphantom{V}\smash{\tilde{V}}_{j}^{\ast}\zeta_{i\,0}^{{}}
+V_{i}^{\ast}\zeta_{0\,j}^{{}}
+V_{i}^{\ast}\vphantom{V}\smash{\tilde{V}}_{j}^{{}}\zeta_{0}^{{}}\right\|  ^{2}
\geq0.
\end{multline*}%
Note that $\mathcal{K}=\mathcal{H}_{\varphi}$ identifies with a subspace of
$\mathcal{H}$ through the map%
\[
\mathcal{K}\ni\xi\longmapsto\Lambda\left\{  \left\{  \varnothing\right\}
\otimes\left\{  \varnothing\right\}  \otimes\xi\right\}  .
\]
Then $\mathcal{K}=\mathcal{H}_{0}\cap\widetilde{\mathcal{H}}_{0}$, 
so
$\mathcal{H}$ 
may be viewed as an
amalgamated tensor product%
\[
\mathcal{H}=\mathcal{H}_{0}\otimes_{\mathcal{K}}\widetilde{\mathcal{H}}_{0}.%
\]
Let $P$ be the projection from $\mathcal{H}$
onto $\mathcal{K}$. Then%
\begin{align*}
S_{i}^{\ast}P  &  =PS_{i}^{\ast}P=V_{i}^{\ast},\\
\tilde{S}_{i}^{\ast}P  &  =P\tilde{S}_{i}^{\ast}P=\tilde{V}_{i}^{\ast}.
\end{align*}
We can thus define states $\psi$, $\tilde{\psi}$ on $\mathcal{O}_{d}$ through
the requirement%
\begin{align*}
\psi\left(  s_{I}^{{}}s_{J}^{\ast}\right)   &  =\varphi\left(  V_{I}^{{}}%
V_{J}^{\ast}\right)  ,\\
\vphantom{\tilde{V}_{I}^{{}}\tilde{V}_{J}^{\ast}}\tilde{\psi}\left(  s_{I}%
^{{}}s_{J}^{\ast}\right)   &  =\tilde{\varphi}\left(  \vphantom{V_{I}^{{}%
}V_{J}^{\ast}}\smash{\tilde{V}_{I}^{{}}\tilde{V}_{J}^{\ast}}\right)  .
\end{align*}
Let $E$ be the support projection of $\psi$ as a state on $\mathcal{O}%
_{d}^{\prime\prime}$ on the amalgamated tensor product, and similarly let
$\tilde{E}$ be the support projection of $\tilde{\psi}$. Here $\mathcal{O}%
_{d}^{\prime\prime}$ and $\widetilde{\mathcal{O}}_{d}^{\prime\prime}$ denote
the von Neumann algebras generated by $\left\{  S_{1},\dots,S_{d}\right\}  $
and $\vphantom{\tilde{S}_{1},\dots,\tilde{S}_{d}}\left\{  \vphantom{S_{1}%
S_{d}}\smash{\tilde{S}_{1},\dots,\tilde{S}_{d}}\right\}  $
of (\ref{eqNew7.1})--(\ref{eqNew7.2}),
respectively. (The
amalgamated tensor product thus carries a representation of $\mathcal{O}%
_{d}\otimes\widetilde{\mathcal{O}}_{d}$, where $\widetilde{\mathcal{O}}%
_{d}\cong\mathcal{O}_{d}$, and the states $\psi$, $\tilde{\psi}$ identify with
the restriction of the vector state $\ip{\Phi}{\,\cdot\,\Phi}$ to each of the
two tensor factors.)

\begin{lemma}
\label{Lem7.5}$\mathcal{H}_{0}$ is an invariant subspace for $\mathcal{O}_{d}$
and $P|_{\mathcal{H}_{0}}=E|_{\mathcal{H}_{0}}$.
\end{lemma}

\begin{proof}
$\mathcal{H}_{0}$ is obviously an invariant subspace for $\mathcal{O}_{d}$ by
construction. But the map%
\[
\left(  \mathcal{O}_{d}|_{\mathcal{H}_{0}}\right)  ^{\prime}\longrightarrow
\mathcal{B}\left(  \mathcal{K}\right)  ^{\bs}=\mathcal{M}^{\prime
}\colon Q\longmapsto PQP
\]
is an order isomorphism onto $\mathcal{M}^{\prime}$ by Proposition
\ref{Pro4.1} and the assumptions of Theorem \ref{Thm7.1}. Hence, as
$\mathcal{M}^{\prime}$ is a factor by assumption and $\left(  \mathcal{O}%
_{d}|_{\mathcal{H}_{0}}\right)  ^{\prime}$ is a von Neumann algebra, the map
$Q\mapsto PQP$ is either an isomorphism or anti-isomorphism by
\cite[Proposition 3.22]{BrRoI}, or \cite{JaRi50,Kad51,Kad65}. But as the map
$Q\mapsto PQP$ is clearly completely positive, it is an isomorphism. Hence
$P\in\left(  \left(  \mathcal{O}_{d}|_{\mathcal{H}_{0}}\right)  ^{\prime
}\right)  ^{\prime}=\mathcal{O}_{d}^{\prime\prime}|_{\mathcal{H}_{0}}$ by
Proposition \ref{Pro4.2}. But as $P$ is the support projection of the normal
state $\psi$ on $\mathcal{O}_{d}^{\prime\prime}|_{\mathcal{H}_{0}}$,
it follows
that $P$ is the image of $E$ under the map $\mathcal{O}_{d}^{\prime\prime
}\ni A\mapsto A|_{\mathcal{H}_{0}}$, i.e.,%
\[
P|_{\mathcal{H}_{0}}=E|_{\mathcal{H}_{0}}.%
\settowidth{\qedskip}{$\displaystyle
P|_{\mathcal{H}_{0}}=E|_{\mathcal{H}_{0}}.$}
\addtolength{\qedskip}{-\textwidth}
\rlap{\hbox to-0.5\qedskip{\hfil\qed}}%
\]%
\renewcommand{\qed}{}%
\end{proof}

\begin{lemma}
\label{Lem7.6}$E\tilde{E}=P$.
\end{lemma}

\begin{proof}
Clearly $E\geq P$, $\tilde{E}\geq P$, so $E\tilde{E}\geq P$. 
The converse inequality follows by
using Lemma \ref{Lem7.5}, $E\tilde{E}=\tilde{E}E$,
and $E\tilde{S}_{\tilde{I}}=\tilde{S}_{\tilde{I}}E$:
\begin{align*}
\tilde{E}E
\vphantom{\tilde{I}}\Lambda\left(  I\otimes\smash{\tilde{I}}\otimes
\xi\right)   
&  =\tilde{E}\tilde{S}_{\tilde{I}}E
\Lambda\left(  I\otimes
\left\{  \varnothing\right\} 
\otimes\xi\right)  \\
&  =\tilde{E}\tilde{S}_{\tilde{I}}P
\Lambda\left(  I\otimes
\left\{  \varnothing\right\} 
\otimes\xi\right)  \\
&  =\tilde{E}\tilde{S}_{\tilde{I}}PS_{I}\xi  \\
&  =\tilde{E}\tilde{S}_{\tilde{I}}V_{I}\xi  \\
&  =\tilde{E}\vphantom{\tilde{I}}\Lambda\left(
\left\{  \varnothing\right\} 
\otimes\smash{\tilde{I}}\otimes V_{I}\xi\right)   \\
&  =P\vphantom{\tilde{I}}\Lambda\left(
\left\{  \varnothing\right\} 
\otimes\smash{\tilde{I}}\otimes V_{I}\xi\right)   \\
&  =\tilde{V}_{\tilde{I}}V_{I}\xi \in \mathcal{K}, \\%
\intertext{so}%
\tilde{E}E &\leq P.
\settowidth{\qedskip}{$\displaystyle 
\tilde{E}E
\vphantom{\tilde{I}}\Lambda\left(  I\otimes\smash{\tilde{I}}\otimes
\xi\right)   
\leq P.$}
\settowidth{\qedadjust}{$\displaystyle 
\tilde{E}E
\vphantom{\tilde{I}}\Lambda\left(  I\otimes\smash{\tilde{I}}\otimes
\xi\right)   
=\tilde{E}\tilde{S}_{\tilde{I}}E
\Lambda\left(  I\otimes
\left\{  \varnothing\right\} 
\otimes\xi\right)  $}
\addtolength{\qedadjust}{\textwidth}
\addtolength{\qedskip}{-0.5\qedadjust}
\rlap{\hbox to-\qedskip{\hfil $\qedsymbol $}\hss }
\end{align*}
\renewcommand{\qed}{}%
\end{proof}

\begin{lemma}
\label{Lem7.7}$\mathcal{O}_{d}^{\prime\prime}\vee\widetilde{\mathcal{O}}%
_{d}^{\prime\prime}=\left(  \mathcal{O}_{d}^{{}}\cup\smash{\widetilde
{\mathcal{O}}_{d}^{{}}}\right)  ^{\prime\prime}=\mathcal{B}\left(
\mathcal{H}\right)  $.
\end{lemma}

\begin{proof}
We have $P=E\tilde{E}\in\mathcal{O}_{d}^{\prime\prime}\vee\widetilde
{\mathcal{O}}_{d}^{\prime\prime}$. But by Lemma \ref{Lem7.5} (applied both to
$\mathcal{O}_{d}$ and $\widetilde{\mathcal{O}}_{d}$), $P\vphantom
{\widetilde{\mathcal{O}}_{d}^{\prime\prime}}\left(  \mathcal{O}_{d}%
^{\prime\prime}\vee\smash{\widetilde{\mathcal{O}}_{d}^{\prime\prime}}\right) P
$ contains both $\mathcal{M}$ and $\mathcal{M}^{\prime}$, and as $\mathcal{M}$
is a factor, we have%
\[
P\vphantom{\widetilde{\mathcal{O}}_{d}^{\prime\prime}}\left(  \mathcal{O}%
_{d}^{\prime\prime}\vee\smash{\widetilde{\mathcal{O}}_{d}^{\prime\prime}%
}\right)  P=\mathcal{B}\left(  \mathcal{K}\right)  .
\]
Since $\mathcal{K}$ is cyclic for $\mathcal{O}_{d}^{\prime\prime}%
\vee\widetilde{\mathcal{O}}_{d}^{\prime\prime}$, Lemma \ref{Lem7.7} follows.
\end{proof}

\begin{lemma}
\label{Lem7.8}$\mathcal{O}_{d}^{\prime}=\widetilde{\mathcal{O}}_{d}%
^{\prime\prime}$.
\end{lemma}

\begin{proof}
Since $\mathcal{O}_{d}$ and $\widetilde{\mathcal{O}}_{d}$ mutually commute, it
follows from Lemma \ref{Lem7.7} that $\mathcal{O}_{d}^{\prime}$ and
$\widetilde{\mathcal{O}}_{d}^{\prime\prime}$ are factors, and $\mathcal{O}%
_{d}^{\prime}\supset\widetilde{\mathcal{O}}_{d}^{\prime\prime}$. The
projection $E\in\mathcal{O}_{d}^{\prime\prime}$ commutes with $\mathcal{O}%
_{d}^{\prime}$ and thus also with $\widetilde{\mathcal{O}}_{d}^{\prime\prime}%
$, and hence it suffices to show that $\mathcal{O}_{d}^{\prime}E=\widetilde
{\mathcal{O}}_{d}^{\prime\prime}E$. Introduce $\mathcal{N}_{1}=\mathcal{O}%
_{d}^{\prime}E$ and $\mathcal{N}_{2}=\widetilde{\mathcal{O}}_{d}^{\prime
\prime}E$. Then $P\in\mathcal{N}_{2}\subset\mathcal{N}_{1}$, and
\begin{align*}
P\mathcal{N}_{1}^{\prime}P  &  =\mathcal{M},\\
P\mathcal{N}_{2}P  &  =\mathcal{M}^{\prime},
\end{align*}
so%
\[
P\mathcal{N}_{1}P=P\mathcal{N}_{2}P=\mathcal{M}^{\prime}.
\]
The
relations $P\in\mathcal{N}_{2}\subset\mathcal{N}_{1}$, and $P\mathcal{N}%
_{1}P=P\mathcal{N}_{2}P$, imply $\mathcal{N}_{1}=\mathcal{N}_{2}$. 
(We can find a type
$\mathrm{I}$ subfactor $M$ 
of $\mathcal{N}_{2}$
such that $P
$ dominates a minimal projection in $
M$, and the above conditions imply
that
\[
\mathcal{N}_{1}\cap M^{\prime}=\mathcal{N}_{2}\cap M^{\prime},
\]
which again implies $\mathcal{N}_{1}=\mathcal{N}_{2}$.)
\end{proof}

Denote the gauge action of $\mathbb{T}$ on $\mathcal{O}_{d}$, respectively
$\widetilde{\mathcal{O}}_{d}$, by $\tau$, respectively $\tilde{\tau}$. Define%
\[
H=\left\{  z\in\mathbb{T}\mid\tau_{z}\text{ extends to an automorphism of
}\mathcal{O}_{d}^{\prime\prime}\right\}  .
\]
As in Theorem \ref{Thm6.3}, it follows from $\psi\circ\lambda=\psi$ that%
\[
H=\left\{  z\in\mathbb{T}\mid\psi\circ\tau_{z}=\psi\right\}  ,
\]
and hence $H$ is a closed subgroup of $\mathbb{T}$. Define a subgroup
$\tilde{H}$ in the same way as $H$ by using $\tilde{\tau}$ instead of $\tau$.

As mentioned before, the algebra $\mathcal{O}_{d}\otimes\widetilde
{\mathcal{O}}_{d}$ is naturally represented on $\mathcal{H}$. Define
\[
G=\vphantom{\widetilde{\mathcal{O}}_{d}^{\prime\prime}}\left\{  \left(
z_{1},z_{2}\right)  \in\mathbb{T}^{2}\mid\tau_{z_{1}}\otimes\tilde{\tau
}_{z_{2}}\text{ extends to an automorphism of }\mathcal{O}_{d}^{\prime\prime
}\vee\smash{\widetilde{\mathcal{O}}_{d}^{\prime\prime}}=\mathcal{B}\left(
\mathcal{H}\right)  \right\}  .
\]

\begin{lemma}
\label{Lem7.9}$\tilde{H}=H$ and $\left\{  \left(  z,z\right)  \mid z\in
H\right\}  \subset G\subset H\times H$.
\end{lemma}

\begin{proof}
Once we can show $\left\{  \left(  z,z\right)  \mid z\in H\right\}  \subset G$
it follows that $H\subset\tilde{H}$, and then it follows by symmetry that
$\tilde{H}\subset H$, so $H=\tilde{H}$. But then $G\subset H\times H$ is
obvious. So it remains to show%
\[
\left\{  \left(  z,z\right)  \mid z\in H\right\}  \subset G.
\]
For this, let $z\in H$, i.e., $\psi\circ\tau_{z}=\psi$. Thus $\tau_{z}\left(
E\right)  =E$ (where still $E=\operatorname*{supp}\psi$), and one can define a
unitary operator $U_{0}$ on $P\mathcal{H}=\mathcal{K}$ by%
\[
U_{0}Q\Phi=\tau_{z}\left(  Q\right)  \Phi
\]
for $Q\in\mathcal{M}$. (We are now working in the cyclic representation
defined by $\varphi$, and $\tau_{z}$ also denotes the extension of $\tau_{z}$
to $\mathcal{O}_{d}^{\prime\prime}$.) If $J$, $\Delta$ are the modular
conjugation and modular operator associated to $\left(  \mathcal{M}%
,\Phi\right)  $, it follows from $\tau_{z}$-invariance that%
\begin{align*}
U_{0}^{{}}JU_{0}^{\ast}  &  =J,\\
U_{0}^{{}}\Delta U_{0}^{\ast}  &  =\Delta.
\end{align*}
Thus%
\begin{align*}
\operatorname{Ad}U_{0}^{{}}\vphantom{\tilde{V}_{j}^{{}}}\left(  \vphantom
{V_{j}^{{}}}\smash{\tilde{V}_{j}^{{}}}\right)   &  =U_{0}^{{}}\left(
J\Delta^{-\frac{1}{2}}V_{j}^{\ast}\Delta^{\frac{1}{2}}J\right)  U_{0}^{\ast}\\
&  =z\tilde{V}_{j}.
\end{align*}
Using this, we can extend $U_{0}$ to a unitary operator $U$ on $\mathcal{H}$
by the definition%
\[
U\Lambda\left(  I\otimes\smash{\tilde{I}}\otimes\xi\right)  =z^{\left|
I\right|  +\left|  \smash{\tilde{I}}\right|  }\Lambda\left(  I\otimes
\smash{\tilde{I}}\otimes U_{0}\xi\right)  .
\]
This operator $U$ is indeed well defined and unitary because%
\[
\ip{\Lambda\left( I\otimes\smash{\tilde{I}}\otimes U_{0}\xi\right) }
{\Lambda\left( J\otimes\smash{\tilde{J}}\otimes U_{0}\eta\right) }=z^{\left|
I\right|  +\left|  \smash{\tilde{I}}\right|  -\left|  J\right|  -\left|
\smash{\tilde{J}}\right|  }\ip{\Lambda\left( I\otimes\smash{\tilde{I}}
\otimes\xi\right) }{\Lambda\left( J\otimes\smash{\tilde{J}}\otimes\eta
\right) }.
\]%
We have
\begin{align*}
US_{i}\Lambda\left(  I\otimes\smash{\tilde{I}}\otimes\xi\right)   &
=z^{\left|  I\right|  +1+\left|  \smash{\tilde{I}}\right|  }\Lambda\left(
iI\otimes\smash{\tilde{I}}\otimes U_{0}\xi\right) \\
&  =zS_{i}\Lambda\left(  I\otimes\smash{\tilde{I}}\otimes\xi\right)  ,
\end{align*}
and similarly $U\tilde{S}_{i}=z\tilde{S}_{i}U$. Hence, $\operatorname{Ad}%
U|_{\mathcal{O}_{d}\otimes\widetilde{\mathcal{O}}_{d}}=\tau_{z}\otimes
\tilde{\tau}_{z}$, so $\left(  z,z\right)  \in G$.
\end{proof}

We next prove an analogue of Theorem \ref{Thm6.3} in this situation.

\begin{lemma}
\label{Lem7.10}If $z\in\mathbb{T}$, the following conditions are equivalent.

\begin{enumerate}
\item \label{Lem7.10(1)}$\left(  z,1\right)  \in G$.

\item \label{Lem7.10(2)}$\tau_{z}$ extends to an inner automorphism of
$\mathcal{O}_{d}^{\prime\prime}$.

\item \label{Lem7.10(3)}$z\in\operatorname{PSp}\left(  \bs|_{\mathcal{M}%
}\right)  $.
\end{enumerate}
\end{lemma}

\begin{proof}
(\ref{Lem7.10(1)}) $\Rightarrow$ (\ref{Lem7.10(2)}): If $\left(  z,1\right)
\in G$, there is a unitary $U$ on $\mathcal{H}$ such that $\tau_{z}%
\otimes\operatorname{id}=\operatorname{Ad}U$. But then $\operatorname{Ad}%
U|_{\widetilde{\mathcal{O}}_{d}}=\operatorname{id}$, i.e., $U\in
\widetilde{\mathcal{O}}_{d}^{\prime}=\mathcal{O}_{d}^{\prime\prime}$ (by Lemma
\ref{Lem7.8}) and hence $\tau_{z}=\operatorname{Ad}U|_{\mathcal{O}_{d}}$
extends to an inner automorphism of $\mathcal{O}_{d}^{\prime\prime}$.

(\ref{Lem7.10(2)}) $\Rightarrow$ (\ref{Lem7.10(3)}): If (\ref{Lem7.10(2)})
holds, then $\psi\circ\tau_{z}=\psi$ by the comment prior to Lemma
\ref{Lem7.9}, and hence $UEU^{\ast}=E$ for a unitary $U\in\mathcal{O}%
_{d}^{\prime\prime}$ with $\tau_{z}=\operatorname{Ad}U$. Thus $U_{0}%
=UE\tilde{E}=UP$ is a unitary in $\mathcal{M}$ with%
\begin{align*}
\bs\left(  U_{0}\right)   &  =\bar{z}U_{0}\\%
\intertext{and thus}%
\bs\left(  U_{0}^{\ast}\right)   &  =zU_{0}^{\ast}.
\end{align*}
Thus $z\in\operatorname{PSp}\left(  \bs|_{\mathcal{M}}\right)  $.

(\ref{Lem7.10(3)}) $\Rightarrow$ (\ref{Lem7.10(1)}): Since $\mathcal{M}^{\bs
}=\mathbb{C}%
\mkern2mu%
\openone$, it follows from (\ref{Lem7.10(3)}) and the beginning of
the proof of Lemma \ref{Lem6.6} that there exists a unitary operator $U_{0}%
\in\mathcal{M}$ with $\bs\left(  U_{0}\right)  =zU_{0}$, and from Lemma
\ref{Lem6.7} it follows that%
\[
\operatorname{Ad}U_{0}\left(  V_{i}\right)  =\bar{z}V_{i}.
\]
Proceeding as in the final parts of the proofs of Lemma \ref{Lem6.6} and Lemma
\ref{Lem7.9}, we extend $U_{0}$ to a unitary $U$ on $\mathcal{H}$ by%
\[
U\Lambda\left(  I\otimes\smash{\tilde{I}}\otimes\xi\right)  =\bar
{z}^{\,\left|  I\right|  }\Lambda\left(  I\otimes\smash{\tilde{I}}\otimes
U_{0}\xi\right)  .
\]
We check as there that $U$ is a well-defined unitary, and that
$\operatorname{Ad}\left(  U^{\ast}\right)  =\tau_{z}\otimes\operatorname{id}$.
Thus $\left(  z,1\right)  \in G$.
\end{proof}

We now show that the representation of $\operatorname*{UHF}\nolimits_{d}%
\subset\mathcal{O}_{d}$ on $\mathcal{H}$ is quasi-equivalent to the
subrepresentation on $\left[  \operatorname*{UHF}\nolimits_{d}\Phi\right]  $:

\begin{lemma}
\label{Lem7.11}The representation of $\operatorname*{UHF}\nolimits_{d}$ on
$\mathcal{H}$ is quasi-equivalent to $\pi_{\psi|_{\operatorname*{UHF}%
\nolimits_{d}}}$.
\end{lemma}

\begin{proof}
Since $\Phi$ is cyclic for the representation of $\mathcal{O}_{d}%
\otimes\widetilde{\mathcal{O}}_{d}$ on $\mathcal{H}$ by Lemma \ref{Lem7.7},
the vectors%
\[
\xi=S_{I}^{{}}S_{J}^{\ast}\tilde{S}_{\tilde{I}}^{{}}\tilde{S}_{\tilde{J}%
}^{\ast}\Phi
\]
span a dense subspace of $\mathcal{H}$, and thus it suffices to show that
$\omega_{\xi}$ is normal in $\pi_{\psi}|_{\operatorname*{UHF}\nolimits_{d}}$.
For this, if $n\geq\left|  I\right|  $ and $x\in\operatorname*{UHF}%
\nolimits_{d}$, we have, using the Cuntz relations,
\begin{align*}
\ip{\xi}{\lambda^{n}\left( x\right) \xi}  &  =\vphantom{S_{I}^{{}}\tilde
{S}_{\tilde{I}}^{{}} S_{J}^{\ast}\tilde{S}_{\tilde{J}}^{\ast}\Phi\sum
_{\left| I^{\prime}\right| =n} S_{I^{\prime}}^{{}}xS_{I^{\prime}}^{\ast}%
S_{I}^{{}}\tilde{S}_{\tilde{I}}^{{}} S_{J}^{\ast}\tilde{S}_{\tilde{J}}^{\ast
}\Phi}\ip{S_{I}^{{}}\tilde{S}_{\tilde{I}}^{{}} S_{J}^{\ast}\tilde{S}%
_{\tilde{J}}^{\ast}\Phi} {\vphantom{\sum}\smash{\sum_{\left| I^{\prime}%
\right| =n}} S_{I^{\prime}}^{{}}xS_{I^{\prime}}^{\ast} S_{I}^{{}}\tilde
{S}_{\tilde{I}}^{{}} S_{J}^{\ast}\tilde{S}_{\tilde{J}}^{\ast}\Phi}\\
&  =\ip{\tilde{S}_{\tilde{I}}^{{}}S_{J}^{\ast} \tilde{S}_{\tilde{J}}^{\ast
}\Phi}{\lambda^{n-\left| I\right| }\left( x\right) \tilde{S}_{\tilde{I}}^{{}%
}S_{J}^{\ast} \tilde{S}_{\tilde{J}}^{\ast}\Phi}\\
&  =\ip{S_{J}^{\ast}\tilde{S}_{\tilde{J}}^{\ast}\Phi} {\lambda^{n-\left
| I\right| }\left( x\right) S_{J}^{\ast}\tilde{S}_{\tilde{J}}^{\ast}\Phi}.
\end{align*}
Since $\tilde{S}_{\tilde{J}}^{\ast}$ commute with the other
factors in this inner product and 
$\tilde{S}_{\tilde{J}}^{{}}\tilde{S}_{\tilde{J}}^{\ast}$
is a projection, it follows
for positive $x$ that
\begin{align*}
\ip{\xi}{\lambda^{n}\left( x\right) \xi}  
&  \leq\ip{\Phi
}{S_{J}^{{}}\lambda^{n-\left| I\right| }%
\left( x\right) S_{J}^{\ast}\Phi}\\
&  \leq\sum_{
\vphantom{J^{\prime}}\left|
\vphantom{J}\smash{J^{\prime}}\right|  =\left|  J\right|  }\ip{\Phi
}{S_{J^{\prime}}^{{}}\lambda^{n-\left| I\right| }%
\left( x\right) S_{J^{\prime}}^{\ast}\Phi}\\
&  =\ip{\Phi} {\lambda^{n-\left| I\right| +\left| J\right| 
}\left( x\right) \Phi}\\
&  =\psi\circ\lambda^{n-\left|  I\right|  +\left|  J\right|  
}\left(  x\right) \\
&  =\psi\left(  x\right)  ,
\end{align*}
where the last identity follows from Remark \ref{Rem6.4}. Hence, $\omega_{\xi
}\circ\lambda^{n}|_{\operatorname*{UHF}\nolimits_{d}}$ is a vector state in
the $\psi|_{\operatorname*{UHF}\nolimits_{d}}$-representation, and since%
\[
\operatorname*{UHF}\nolimits_{d}\cong M_{d^{n}}\otimes\lambda^{n}\left(
\operatorname*{UHF}\nolimits_{d}\right)
\]
in a canonical fashion, it follows that $\omega_{\xi}|_{\operatorname*{UHF}%
\nolimits_{d}}$ is normal in $\psi|_{\operatorname*{UHF}\nolimits_{d}}$. See
\cite[proof of Lemma 5.2]{BJP96}.
\end{proof}

\begin{lemma}
\label{Lem7.12}Let $\mathcal{N}$ be a factor and $\alpha$ an action of a group
$G$ on $\mathcal{N}$. Assume that $G$ is the circle group or a finite cyclic
group. If the fixed point algebra $\mathcal{N}^{\alpha}$ is a factor, then
$\mathcal{N}\cap\left(  \mathcal{N}^{\alpha}\right)  ^{\prime}$ is the abelian
von Neumann algebra generated by a unitary operator $V$ such that%
\[
\alpha_{g}\left(  V\right)  =\ip{g}{\gamma_{0}}V
\]
for some $\gamma_{0}\in\hat{G}$, and%
\[
\operatorname*{Ad}V|_{\mathcal{N}}=\alpha_{h}%
\]
for some $h\in G$.
\end{lemma}

\begin{proof}
Since $\mathcal{N}^{\alpha}$ is a factor, $\alpha|_{\mathcal{N}\cap\left(
\mathcal{N}^{\alpha}\right)  ^{\prime}}$ is ergodic, and therefore
$\operatorname*{Sp}\left(  \alpha|_{\mathcal{N}\cap\left(  \mathcal{N}%
^{\alpha}\right)  ^{\prime}}\right)  $ is a subgroup of $\hat{G}$, and thus
$\operatorname*{Sp}\left(  \alpha|_{\mathcal{N}\cap\left(  \mathcal{N}%
^{\alpha}\right)  ^{\prime}}\right)  $ is either a finite cyclic group or
$\mathbb{Z}$. In any case, $\operatorname*{Sp}\left(  \alpha|_{\mathcal{N}%
\cap\left(  \mathcal{N}^{\alpha}\right)  ^{\prime}}\right)  $ has a generator
$\gamma_{0}$, and by ergodicity the corresponding eigen-subspace is the linear
span of a unitary operator $V$ \cite[Section 3.2.3]{BrRoI}, \cite{OPS77},
\cite{Bra86}. Also $V$ generates $\mathcal{N}\cap\left(  \mathcal{N}^{\alpha
}\right)  ^{\prime}$ as a von Neumann algebra, so this algebra is abelian. Let%
\[
\beta=\operatorname*{Ad}V|_{\mathcal{N}}.
\]
Since $\alpha_{g}\left(  V\right)  =\ip{g}{\gamma_{0}}V$ for all $g\in G$, we
have%
\[
\beta\alpha_{g}=\alpha_{g}\beta
\]
for all $g\in G$, and hence $\beta$ fixes each spectral subspace of $\alpha$
in $\mathcal{N}$. As $V\in\left(  \mathcal{N}^{\alpha}\right)  ^{\prime}$, we
have%
\[
\beta|_{\mathcal{N}^{\alpha}}=\operatorname*{id}.
\]
Since $\mathcal{N}^{\alpha}$ is a factor, each spectral subspace
$\mathcal{N}^{\alpha}\left(  \gamma\right)  $ for $\gamma\in\hat{G}$ either is
$0$, or has the form%
\[
\mathcal{N}^{\alpha}\left(  \gamma\right)  =
\mathcal{N}^{\alpha}V\left(  \gamma\right)
\]
for an isometry $V\left(  \gamma\right)  \in\mathcal{N}^{\alpha}\left(
\gamma\right)  $,
or the form%
\[
\mathcal{N}^{\alpha}\left(  \gamma\right)  =
V\left(  \gamma\right) \mathcal{N}^{\alpha}
\]
for a coisometry $V\left(  \gamma\right)  \in\mathcal{N}^{\alpha}\left(
\gamma\right)  $.
We may assume that $G$ acts faithfully, and this excludes
the case $\mathcal{N}^{\alpha}\left(  \gamma\right)  =0$.

Now
consider the case that $V\left(  \gamma\right) $ is an isometry.
Since $\beta\left(  V\left(  \gamma\right)  \right)  \in\mathcal{N}%
^{\alpha}\left(  \gamma\right)  $, there is an operator $U\left(
\gamma\right)  \in\mathcal{N}^{\alpha}$ such that%
\[
\beta\left(  V\left(  \gamma\right)  \right)  =U\left(  \gamma\right)
V\left(  \gamma\right)  .
\]
Since the projection
$Q=V\left( \gamma\right)  V\left(  \gamma\right)  ^{\ast}$
is in
$\mathcal{N}^{\alpha}$,
we may replace
$U\left( \gamma\right) $ by
$U\left(  \gamma\right) Q$
without changing the equation above,
and then
\[
U\left(  \gamma\right) Q
=U\left(  \gamma\right) .
\]
Then since
$U\left(  \gamma\right) 
=\beta \left(  V\left(  \gamma\right) \right) 
V\left(  \gamma\right)  ^{\ast}$,
we have likewise
\begin{align*}
QU\left(  \gamma\right) 
&=V\left(  \gamma\right) 
V\left(  \gamma\right)  ^{\ast}
\beta \left(  V\left(  \gamma\right) \right) 
V\left(  \gamma\right)  ^{\ast}  \\
&=\beta \left(  V\left(  \gamma\right) 
V\left(  \gamma\right)  ^{\ast}
V\left(  \gamma\right) \right) 
V\left(  \gamma\right)  ^{\ast}  \\
&=\beta \left(  V\left(  \gamma\right) \right) 
V\left(  \gamma\right)  ^{\ast}  \\
&=U\left(  \gamma\right) ,  
\end{align*}
so
\[
U\left(  \gamma\right) \in Q\mathcal{N}^{\alpha}Q.
\]
If $A\in Q\mathcal{N}^{\alpha}Q\subset \mathcal{N}^{\alpha}$, 
then $AV\left(  \gamma\right)  \in
\mathcal{N}^{\alpha}\left(  \gamma\right)  $, and we have on one side%
\begin{align*}
\beta\left(  AV\left(  \gamma\right)  \right)   &  =\beta\left(  A\right)
\beta\left(  V\left(  \gamma\right)  \right) \\
&  =AU\left(  \gamma\right)  V\left(  \gamma\right)
\end{align*}
and on the other side, since $V\left(  \gamma\right)  ^{\ast}AV\left(
\gamma\right)  \in\mathcal{N}^{\alpha}$,%
\begin{align*}
\beta\left(  AV\left(  \gamma\right)  \right)   &  =\beta\left(  V\left(
\gamma\right)  V\left(  \gamma\right)  ^{\ast}AV\left(  \gamma\right)  \right)
\\
&  =\beta\left(  V\left(  \gamma\right)  \right)  \beta\left(  V\left(
\gamma\right)  ^{\ast}AV\left(  \gamma\right)  \right) \\
&  =U\left(  \gamma\right)  V\left(  \gamma\right)  \left(  V\left(
\gamma\right)  ^{\ast}AV\left(  \gamma\right)  \right) \\
&  =U\left(  \gamma\right)  AV\left(  \gamma\right)  .
\end{align*}
Comparing the last two expressions, we see that 
\[
AU\left(  \gamma\right)
=U\left( \gamma\right) A
\]
for all $A\in Q\mathcal{N}^{\alpha}Q$. Hence%
\[
U\left(  \gamma\right)  \in
Q\left( 
\left(  \mathcal{N}^{\alpha}\right)  ^{\prime}%
\cap\mathcal{N}^{\alpha}
\right) Q
=\mathbb{C}%
\mkern2mu%
Q,
\]
so $U\left(  \gamma\right)  $ is a scalar multiple $f\left(  \gamma\right)  
\in\mathbb{T}$ of 
$Q$, and%
\[
\beta\left(  B\right)  =f\left(  \gamma\right)  B
\]
for all $B\in\mathcal{N}^{\alpha}\left(  \gamma\right) 
=\mathcal{N}^{\alpha}V\left(  \gamma\right)  $. 
If $V\left( \gamma\right) $ is a coisometry, the
verification of this relation is analogous.
Since $\beta$ is an
automorphism, one verifies that $f\in\Hat{\Hat{G}}=G$, so there exists an
$h\in G$ with%
\[
\beta\left(  B\right)  =\ip{h}{\gamma}B=\alpha_{h}\left(  B\right)
\]
for $B\in\mathcal{N}^{\alpha}\left(  \gamma\right)  $ and $\gamma\in\hat{G}$;
and hence%
\[
\beta=\alpha_{h}.
\settowidth{\qedskip}{$\displaystyle
\beta=\alpha_{h}.$}
\addtolength{\qedskip}{-\textwidth}
\rlap{\hbox to-0.5\qedskip{\hfil\qed}}%
\]%
\renewcommand{\qed}{}%
\end{proof}

\begin{lemma}
\label{Lem7.13}If $\omega$ is a translationally invariant state on
$\mathfrak{A}=\bigotimes_{\mathbb{Z}}M_{d}$ such that $\pi_{\omega}\left(
\mathfrak{A}\right)  ^{\prime\prime}$ is a type $\mathrm{I}$ factor, then
$\omega$ is pure.
\end{lemma}

\begin{proof}
Let $\pi$ be an irreducible representation quasi-equivalent to $\pi_{\omega}$.
There is a density matrix $\rho$ on $\mathcal{H}_{\pi}$ such that%
\[
\omega\left(  x\right)  =\operatorname*{Tr}\left(  \pi\left(  x\right)
\rho\right)
\]
for $x\in\mathfrak{A}$. Since $\pi_{\omega}$ and thus $\pi$ are translationally
covariant, there is a unitary operator $U$ on $\mathcal{H}_{\pi}$ such that%
\[
U\pi\left(  x\right)  U^{\ast}=\pi\left(  \lambda\left(  x\right)  \right)
\]
for all $x\in\mathfrak{A}$, where $\lambda$ is the translation automorphism.
Since $\omega\circ\lambda=\omega$, we obtain that%
\[
U^{\ast}\rho U=\rho.
\]
Assume \emph{ad absurdum} that $\rho$ is not a one-dimensional projection.
Then there are at least two orthogonal eigenvectors $\xi_{1}$, $\xi_{2}$ of
$U$. Thus, for any $x\in A$,%
\begin{align*}
\ip{\xi_{i}}{\pi\left( \lambda^{n}\left( x\right) \right) \xi_{i}}  &
=\ip{U^{\ast\,n}\xi_{i}}{\pi\left( x\right) U^{\ast\,n}\xi_{i}}\\
&  =\ip{\xi_{i}}{\pi\left( x\right) \xi_{i}}.
\end{align*}
But any weak*-limit point of $\pi\left(  \lambda^{n}\left(  x\right)  \right)
$, as $n\rightarrow\infty$, is in $\pi\left(  \mathfrak{A}\right)  ^{\prime
}=\mathbb{C}%
\mkern2mu%
\openone$, and it follows that%
\[
\ip{\xi_{1}}{\pi\left( x\right) \xi_{1}}=\ip{\xi_{2}}{\pi\left( x\right
) \xi_{2}}%
\]
for all $x\in\mathfrak{A}$. But as $\pi$ is irreducible, this is a
contradiction. Thus $\rho$
must be
a one-dimensional projection, and $\omega$ is pure.
\end{proof}

\begin{lemma}
\label{Lem7.14}Let $\mathcal{N}$ be a type $\mathrm{I}$ von Neumann algebra
and $\alpha$ an action of a group $G$ on $\mathcal{N}$. Assume that $G$ is the
circle group or a finite cyclic group. Then the fixed point subalgebra
$\mathcal{N}^{\alpha}$ is of type $\mathrm{I}$.
\end{lemma}

\begin{proof}
By considering the action $\alpha$ on the center $\mathcal{N}\cap
\mathcal{N}^{\prime}$, the von Neumann algebra decomposes into algebras of the
form $\mathcal{M}\otimes L^{\infty}\left(  G\diagup H\right)  $, where
$\mathcal{M}$ is a type $\mathrm{I}$ factor, $H$ is a closed subgroup of $G$,
and $G$ acts on $L^{\infty}\left(  G\diagup H\right)  $ by translation, until
reaching the end of the orbit, in which case the action may be modified by an
automorphism of $\mathcal{M}$ (if $H\neq\left\{  0\right\}  $). The latter
automorphism is inner and of finite order, except in the case $H=G=\mathbb{T}%
$, in which case we have an inner action of $\mathbb{T}$ on the type
$\mathrm{I}$ factor $\mathcal{M}$. In any case, it is clear that the fixed
point subalgebra of $\mathcal{M}\otimes L^{\infty}\left(  G\diagup H\right)  $
under the action is a type $\mathrm{I}$ von Neumann algebra, 
and the lemma follows.
\end{proof}

\begin{lemma}
\label{LemNew7.15}The vector state
$\omega _{\Phi}$ on
\[
\widetilde{\operatorname*{UHF}}
_{d}\otimes\operatorname*{UHF}%
\nolimits_{d}\cong
\vphantom{\bigotimes_{-\infty}^{0}M_{d}\otimes \bigotimes_{1}^{\infty}M_{d}}
\left( \vphantom{\bigotimes}\smash{\bigotimes_{-\infty}^{0}}M_{d}\right) 
\otimes \left( \vphantom{\bigotimes}\smash{\bigotimes_{1}^{\infty}}M_{d}\right) 
\cong\bigotimes_{\mathbb{Z}}^{{}}M_{d}
\]
is equal to $\omega $.
\end{lemma}

\begin{proof}
Since
\[
\tilde{S}_{j}^{\ast}\Phi
=\tilde{V}_{j}^{\ast}\Phi 
=J\Delta^{\frac{1}{2}}V_{j}^{{}}
\Delta^{-\frac{1}{2}}J\Phi
=V_{j}^{\ast}\Phi 
=S_{j}^{\ast}\Phi ,
\]
it follows that
\[
\ip{\Phi}{\tilde{S}_{i\tilde{I}}^{{}}\tilde{S}_{j\tilde{J}}^{\ast}
S_{I}^{{}}S_{J}^{\ast}\Phi}
=\ip{\tilde{S}_{i}^{\ast}\Phi}
{\tilde{S}_{\tilde{I}}^{{}}\tilde{S}_{\tilde{J}}^{\ast}
S_{I}^{{}}S_{J}^{\ast}
\tilde{S}_{j}^{\ast}\Phi}
=\ip{\Phi}{\tilde{S}_{\tilde{I}}^{{}}\tilde{S}_{\tilde{J}}^{\ast}
S_{iI}^{{}}S_{jJ}^{\ast}\Phi}.
\]
This proves Lemma \ref{LemNew7.15} since the
vector state $\omega _{\Phi}$ on 
$\operatorname*{UHF}\nolimits_{d}$ is the
restriction of $\omega $ to 
$\operatorname*{UHF}\nolimits_{d}=\bigotimes_{\mathbb{N}}M_{d}$.
\end{proof}

\begin{proof}
[Proof of Theorem \textup{\ref{Thm7.1}}]
We first merely assume that $\mathcal{M}$ is a factor.
\medskip

(\ref{Thm7.1(1)}) $\Rightarrow$ (\ref{Thm7.1(3)}): This is trivial.
\medskip

(\ref{Thm7.1(3)}) $\Rightarrow$ (\ref{Thm7.1(1)}): By Lemma \ref{Lem7.7},
$\mathcal{O}_{d}^{\prime\prime}$ is a factor, and by assumption
(\ref{Thm7.1(3)}), $\pi_{\omega}\left(  \bigotimes_{\mathbb{Z}}M_{d}\right)
^{\prime\prime}$ is a factor. By viewing $\operatorname*{UHF}\nolimits_{d}%
=\bigotimes_{\mathbb{N}}M_{d}$ as a subalgebra of both $\bigotimes
_{\mathbb{Z}}M_{d}$ and $\mathcal{O}_{d}$, we have%
\[
\omega|_{\bigotimes_{\mathbb{N}}M_{d}}=\psi|_{\operatorname*{UHF}%
\nolimits_{d}}.
\]
Hence, by Lemma \ref{Lem7.11}, $\operatorname*{UHF}\nolimits_{d}^{\prime
\prime}$ is a factor. Recall from the remark prior to Lemma \ref{Lem7.9} that
the set
\[
H=\left\{  z\in\mathbb{T}\mid\tau_{z}\text{ extends to an automorphism of
}\mathcal{O}_{d}^{\prime\prime}\right\}
\]
is a closed subgroup of $\mathbb{T}$, and, as in the proof of Lemma
\ref{Lem6.11}, we have%
\[
\mathcal{O}_{d}^{\prime\prime\,\tau_{H}}=\operatorname*{UHF}\nolimits_{d}%
^{\prime\prime}.
\]
By Lemma \ref{Lem7.12}, the algebra $\mathcal{O}_{d}^{\prime\prime}%
\cap\operatorname*{UHF}\nolimits_{d}^{\prime}$ is abelian, and thus%
\begin{align*}
\left(  \mathcal{O}_{d}^{\prime\prime}\cap\operatorname*{UHF}\nolimits_{d}%
^{\prime}\right)  ^{\prime}  &  =\mathcal{O}_{d}^{\prime}\vee
\operatorname*{UHF}\nolimits_{d}^{\prime\prime}\\
&  =\widetilde{\mathcal{O}}_{d}^{\prime\prime}\vee\operatorname*{UHF}%
\nolimits_{d}^{\prime\prime}%
\end{align*}
is of type $\mathrm{I}$, where we used Lemma \ref{Lem7.8} for the last
identity. For any $z\in\mathbb{T}$,
\[
\tilde{\tau}_{z}\otimes\operatorname*{id}|_{\widetilde{\mathcal{O}}_{d}%
\otimes\operatorname*{UHF}\nolimits_{d}}=\tilde{\tau}_{z}\otimes\tau
_{z}|_{\widetilde{\mathcal{O}}_{d}\otimes\operatorname*{UHF}\nolimits_{d}}%
\]
extends to an automorphism of $\widetilde{\mathcal{O}}_{d}^{\prime\prime}%
\vee\operatorname*{UHF}\nolimits_{d}^{\prime\prime}$ if and only if $z\in H$.
Then, again invoking the argument in the proof of Lemma \ref{Lem6.11}, we
obtain%
\[
\vphantom{\widetilde{\mathcal{O}}}\left(  \vphantom{\mathcal{O}}%
\smash{\widetilde{\mathcal{O}}}_{d}^{\prime\prime}\vee\operatorname*{UHF}%
\nolimits_{d}^{\prime\prime}\right)  ^{\tilde{\tau}_{H}\otimes
\operatorname*{id}}=\widetilde{\operatorname*{UHF}}
\vphantom{\operatorname*{UHF}}_{d}^{\prime\prime}%
\vee\operatorname*{UHF}\nolimits_{d}^{\prime\prime}.
\]
Since $H$ is either $\mathbb{T}$ or a closed subgroup of $\mathbb{T}$, it
follows from Lemma \ref{Lem7.14} that $\displaystyle
\widetilde{\operatorname*{UHF}}%
\vphantom{\operatorname*{UHF}}
_{d}^{\prime\prime}\vee\operatorname*{UHF}\nolimits_{d}^{\prime\prime}$ is of
type $\mathrm{I}$. 
By Lemma \ref{LemNew7.15},
$\omega=\omega_{\Phi
|_{\widetilde{\operatorname*{UHF}}
_{d}\otimes\operatorname*{UHF}%
\nolimits_{d}}}$, 
and
by assumption (\ref{Thm7.1(3)}), the restriction of
$\displaystyle \widetilde{\operatorname*{UHF}}
\vphantom{\operatorname*{UHF}}
_{d}^{\prime\prime}\vee\operatorname*{UHF}%
\nolimits_{d}^{\prime\prime}$ to the closed subspace $\displaystyle
\vphantom{\widetilde
{\operatorname*{UHF}}}\left[  \left(  \smash{\widetilde{\operatorname*{UHF}}}
\vphantom{\operatorname*{UHF}}
_{d}^{\prime\prime}\vee\operatorname*{UHF}\nolimits_{d}^{\prime\prime
}\right)  \Phi\right]  $ is a factor. 
Thus
it follows from Lemma
\ref{Lem7.13} that $\omega$ is pure. This ends the proof of (\ref{Thm7.1(3)})
$\Rightarrow$ (\ref{Thm7.1(1)}).
\medskip

{}From now on, assume that $\mathcal{M}$ is a type $\mathrm{I}$ factor.
\medskip

(\ref{Thm7.1(2)}) $\Rightarrow$ (\ref{Thm7.1(1)}): Since $\mathcal{M}\simeq
P\mathcal{O}_{d}^{\prime\prime}P\simeq E\mathcal{O}_{d}^{\prime\prime}E$ and
$\mathcal{O}_{d}^{\prime\prime}$ is a factor, it follows that $\mathcal{O}%
_{d}^{\prime\prime}$ is a type $\mathrm{I}$ factor. Since any automorphism of
a type $\mathrm{I}$ factor is inner, it follows from Lemma \ref{Lem7.10} that
$H=\left\{  0\right\}  $. Thus $\mathcal{O}_{d}^{\prime\prime}%
=\operatorname*{UHF}\nolimits_{d}^{\prime\prime}$ as in Lemma \ref{Lem6.10}.
Similarly $\displaystyle
\widetilde{\mathcal{O}}_{d}^{\prime\prime}=\widetilde
{\operatorname*{UHF}}
\vphantom{\operatorname*{UHF}}
_{d}^{\prime\prime}$, and it follows from Lemma
\ref{Lem7.7} that%
\[
\widetilde{\operatorname*{UHF}}
\vphantom{\operatorname*{UHF}}
_{d}^{\prime\prime}\vee\operatorname*{UHF}%
\nolimits_{d}^{\prime\prime}=\mathcal{B}\left(  \mathcal{H}\right)  ,
\]
so $\omega$ is pure.
\medskip

(\ref{Thm7.1(1)}) $\Rightarrow$ (\ref{Thm7.1(2)}): If $\mathcal{M}$ is of type
$\mathrm{I}$, it follows again that $\mathcal{O}_{d}^{\prime\prime}$ and
$\widetilde{\mathcal{O}}_{d}^{\prime\prime}$ are type $\mathrm{I}$ factors,
and hence%
\[
\widetilde{\mathcal{O}}_{d}^{\prime\prime}\vee\mathcal{O}_{d}^{\prime\prime
}\simeq\widetilde{\mathcal{O}}_{d}^{\prime\prime}\otimes\mathcal{O}%
_{d}^{\prime\prime}.
\]
Recall that
$\mathcal{O}_{d}^{\prime\prime\,\tau_{H}}=\operatorname*{UHF}
\nolimits_{d}^{\prime\prime}$, 
and that $\tau_{H}$ is inner. It then follows that
$\mathcal{O}_{d}^{\prime\prime}
\cap\operatorname*{UHF}\nolimits_{d}^{\prime}$ is abelian,
and we 
conclude 
that
\[
\widetilde{\operatorname*{UHF}}
\vphantom{\operatorname*{UHF}}
_{d}^{\prime\prime}\vee\operatorname*{UHF}%
\nolimits_{d}^{\prime\prime}\simeq\widetilde{\operatorname*{UHF}}
\vphantom{\operatorname*{UHF}}
_{d}%
^{\prime\prime}\otimes\operatorname*{UHF}\nolimits_{d}^{\prime\prime}%
\]
and 
further that 
the commutant of this algebra is 
again 
abelian. If (\ref{Thm7.1(1)}) holds, it
follows from Lemma \ref{Lem7.11} and Proposition \ref{Pro7.3} that
$\operatorname*{UHF}\nolimits_{d}^{\prime\prime}$ and 
$\displaystyle \widetilde
{\operatorname*{UHF}}
\vphantom{\operatorname*{UHF}}
_{d}^{\prime\prime}$ are factors, and hence%
\[
\widetilde{\operatorname*{UHF}}
\vphantom{\operatorname*{UHF}}
_{d}^{\prime\prime}\vee\operatorname*{UHF}%
\nolimits_{d}^{\prime\prime}=\mathcal{B}\left(  \mathcal{H}\right)  .
\]
This implies that the subgroup $H$ of $\mathbb{T}$ must be trivial, and hence
\[
\operatorname*{PSp}\left(  \bs|_{\mathcal{M}}\right)  \cap\mathbb{T}=\left\{
1\right\}
\]
by Lemma \ref{Lem7.10}, so (\ref{Thm7.1(2)}) holds.
\medskip

Finally, we argue that
\medskip

(\ref{Thm7.1(3)}) $\Leftrightarrow$ (\ref{Thm7.1(4)}): If $\mathcal{N}%
=\pi_{\psi}\left(  \mathcal{O}_{d}\right)  ^{\prime\prime}\cong\mathcal{O}%
_{d}^{\prime\prime}$, then $\mathcal{N}^{\tau|_{H}}=\operatorname*{UHF}%
\nolimits_{d}^{\prime\prime}$ by Lemma \ref{Lem6.10}. But, by Corollary 8.10.5
in \cite{Ped79}, it follows that $\mathcal{N}^{\tau|_{H}}=\operatorname*{UHF}%
\nolimits_{d}^{\prime\prime}$ is a factor if and only if the Connes spectrum
of the extension of $\tau|_{H}$ to $\mathcal{O}_{d}^{\prime\prime}$ is equal
to $\hat{H}$. The rest of the proof of (\ref{Thm7.1(3)}) $\Leftrightarrow$
(\ref{Thm7.1(4)}) is as above.
\end{proof}

\begin{remark}
\label{Rem7.15}If $\mathcal{N}$ is a type $\mathrm{I}$ factor in the last
paragraph of the preceding proof, then $H=\left\{  1\right\}  $ since $H$ is
$\mathbb{T}$ or a
cyclic
group%
, and we get the trivial peripheral spectrum. If $\mathcal{N}$ is a type
$\mathrm{III}$ factor, the subgroup $H$ could in principle be nontrivial,
with some nontrivial $\tau_{t}$ inner, i.e., 
$\operatorname*{PSp}\left(  \bs|_{\mathcal{M}}\right)  \neq \left\{
1\right\} $. This is because 
$\tau_{t}$ could be implemented by a unitary in
$\mathcal{N}$ which would then not be fixed by $H$, and, in this case, the
Connes spectrum still could be $\hat{H}$.
See, e.g., \cite{Spi90} for examples of that.
\end{remark}

\chapter{\label{Hilbert}Remarks on duality}

In this \chaptername{} we will establish several more properties of the duality
theory of the objects $\left(  \mathcal{M},\varphi,V_{1},\dots,V_{d}%
,\bs\right)  $ introduced between Lemmas \ref{Lem7.4} and \ref{Lem7.5},
and thereby also extend \cite[Lemma 6.3]{BrJo97a}.
So
we assume that
$\left(  \mathcal{M},\varphi,V_{1},\dots,V_{d},\bs\right)  $ satisfies the
general hypotheses in Theorem \ref{Thm7.1}, and let the dual object
\[
\vphantom{\widetilde{\mathcal{M}},
\tilde{\varphi},\tilde{V}_{1},\dots,\tilde{V}_{d},\tilde{\bs}}\left(
\vphantom{\mathcal{M}V_{1}V_{d}}\smash{\widetilde{\mathcal{M}},
\tilde{\varphi},\tilde{V}_{1},\dots,\tilde{V}_{d},\tilde{\bs}}\right)  
\]
be
defined as after Lemma \ref{Lem7.4}, i.e.,%
\begin{align*}
\widetilde{\mathcal{M}}  &  =\mathcal{M}^{\prime},\\
\tilde{V}_{j}^{{}}  &  =\left(  J\sigma_{\frac{i}{2}}\left(  V_{j}^{\ast
}\right)  J\right)  \overline{\quad\mathstrut},\\
\tilde{\varphi}\left(  X\right)   &  =\ip{\Phi}{X\Phi},
\end{align*}
where $\Phi$ is the separating and cyclic vector for $\mathcal{M}$ defining
the state
$\varphi$.

Recall the notation%
\begin{align*}
\mathcal{M}^{\bs}  &  =\left\{  X\in\mathcal{M}\mid\bs\left(  X\right)
=X\right\}  ,\\
\mathcal{M}_{\ast}^{\bs}  &  =\left\{  \eta\in\mathcal{M}_{\ast}\mid\eta
\circ\bs=\eta\right\}  ,
\end{align*}
where $\mathcal{M}_{\ast}$ is the predual of $\mathcal{M}$.

\begin{lemma}
\label{Lem8.1}If $\mathcal{M}^{\bs}=\mathbb{C}%
\mkern2mu%
\openone$, then $\mathcal{M}%
_{\ast}^{\bs}=\mathbb{C}%
\mkern2mu%
\varphi$.
\end{lemma}

\begin{proof}
Since $\varphi\circ\bs=\varphi$, and $\bs$ is ergodic, this is established as
in Lemma \ref{Lem6.8} and its proof.
\end{proof}

We use Lemma \ref{Lem8.1} to establish:

\begin{proposition}
\label{Pro8.2}The map $\bs$ is ergodic if and only if the dual map $\tilde
{\bs}$ is ergodic.
\end{proposition}

\begin{proof}
Assume that $X^{\prime}\in\mathcal{M}^{\prime}$ and that%
\[
\tilde{\bs}\left(  X^{\prime}\right)  =X^{\prime}.
\]
This means that%
\[
\sum_{j}V_{j}^{\ast}\Delta^{\frac{1}{2}}JX^{\prime}J\Delta^{\frac{1}{2}}%
V_{j}^{{}}=\Delta^{\frac{1}{2}}JX^{\prime}J\Delta^{\frac{1}{2}},
\]
for example as sesquilinear forms on $\mathcal{M}\Phi$. Define%
\[
X=JX^{\prime}J\in\mathcal{M}.
\]
If $A\in\mathcal{M}$ is an entire element for the modular group $\sigma_{t}$,
we have%
\begin{align*}
\ip{\Delta^{\frac{1}{2}} X ^{\ast}\Phi}{\bs\left( A\right) \Phi} &
=\vphantom{\sum_{j}}\ip{\Delta^{\frac{1}{2}} X ^{\ast}\Phi}{\vphantom{\sum
}\smash{\sum_{j}}V_{j}^{{}}AV_{j}^{\ast}\Phi}\\
&  =\vphantom{\sum_{j}}\ip{X ^{\ast}\Phi}{\vphantom{\sum}\smash{\sum_{j}%
}\sigma_{-\frac{i}{2}}\left( V_{j}^{{}}\right) \sigma_{-\frac{i}{2}}%
\left( A\right) \Delta^{\frac{1}{2}}V_{j}^{\ast}\Phi}\\
&  =\sum_{j}\ip{\Phi}{X\sigma_{-\frac{i}{2}}\left( V_{j}^{{}}\right
) \sigma_{-\frac{i}{2}}\left( A\right) \Delta^{\frac{1}{2}}V_{j}^{\ast}\Phi}\\
&  =\sum_{j}\ip{\Phi}{X\sigma_{-\frac{i}{2}}\left( V_{j}^{{}}\right
) \sigma_{-\frac{i}{2}}\left( A\right) J^{2}_{{}}\Delta^{\frac{1}{2}}%
V_{j}^{\ast}\Phi}\\
&  =\sum_{j}\ip{\Phi}{X\sigma_{-\frac{i}{2}}\left( V_{j}\right
) \sigma_{-\frac{i}{2}}\left( A\right) JV_{j}J\Phi}\\
&  =\sum_{j}\ip{JV_{j}^{\ast}\Phi}{X\sigma_{-\frac{i}{2}}\left( V_{j}^{{}}%
\right) \sigma_{-\frac{i}{2}}\left( A\right) \Phi}\\
&  =\sum_{j}\ip{J
\Delta^{\frac{1}{2}}
\sigma_{\frac{i}{2}}\left( V_{j}^{\ast}\right) 
\Phi}{X\sigma_{-\frac{i}{2}}\left( V_{j}^{{}}%
\right) \sigma_{-\frac{i}{2}}\left( A\right) \Phi}\\
&  =\sum_{j}\ip{\Phi}{\sigma_{\frac{i}{2}}\left( V_{j}^{\ast}\right
) X\sigma_{-\frac{i}{2}}\left( V_{j}^{{}}\right) \sigma_{-\frac{i}{2}}%
\left( A\right) \Phi}\\
&  =\sum_{j}\ip{\Phi}{\Delta^{-\frac{1}{2}}V_{j}^{\ast}\Delta^{\frac{1}{2}%
}X\Delta^{\frac{1}{2}}V_{j}^{{}}A\Phi}.  
\end{align*}
But 
$\sum_{j}V_{j}^{\ast}\Delta^{\frac{1}{2}}X\Delta^{\frac{1}{2}}V_{j}^{{}}
=\Delta^{\frac{1}{2}}X\Delta^{\frac{1}{2}}$
by the identity above, so furthermore:
\begin{align*}
\ip{\Delta^{\frac{1}{2}} X ^{\ast}\Phi}{\bs\left( A\right) \Phi} &  =\ip{\Phi
}{X\Delta^{\frac{1}{2}}A\Phi}\\
&  =\ip{\Delta^{\frac{1}{2}} X ^{\ast}\Phi}{A\Phi}.
\end{align*}
Thus the functional $\eta\in\mathcal{M}_{\ast}$ defined by%
\[
\eta\left(  A\right)  =\ip{\Delta^{\frac{1}{2}} X ^{\ast}\Phi}{A\Phi}%
\]
is in $\mathcal{M}_{\ast}^{\bs}$. Since $\bs$ is ergodic, it follows from
Lemma \ref{Lem8.1} that $\eta\in\mathbb{C}%
\mkern2mu%
\varphi$, and so $\Delta^{\frac
{1}{2}}X^{\ast}\Phi\in\mathbb{C}%
\mkern2mu%
\Phi$, which implies $X\Phi\in\mathbb{C}%
\mkern2mu%
\Phi$,
or $X\in\mathbb{C}%
\mkern2mu%
\openone$. Thus $X^{\prime}\in\mathbb{C}%
\mkern2mu%
\openone$, which
shows that $\tilde{\bs}$ is ergodic. Since $\Tilde{\Tilde{\bs}}=\bs$, the
other implication follows.
\end{proof}

Our next aim is to show that the two sets $\operatorname*{PSp}\left(
\bs\right)  \cap\mathbb{T}$ and 
$\vphantom{\tilde{\bs}}
\operatorname*{PSp}\left(  
\vphantom{\bs}
\smash{\tilde{\bs}}\right)  \cap\mathbb{T}$ 
are equal. First we need a lemma.

\begin{lemma}
\label{Lem8.3}If $\eta\in\mathcal{M}_{\ast}\diagdown\left\{  0\right\}  $ and
\[
\eta\circ\bs=t\eta
\]
for
some
$t\in\mathbb{T}$, then $t$ is an eigenvalue of $\bs$.
\end{lemma}

\begin{proof}
Again we invoke the technique in the proof of Lemma \ref{Lem6.8}: Pick an
$X\in\mathcal{M}$ such that $\eta\left(  X\right)  \neq0$, and take a limit
point $L$ of the sequence%
\[
X_{n}=\frac{1}{n}\sum_{k=0}^{n-1}t^{-k}\bs^{k}\left(  X\right)
\]
in the weak operator topology. Then it follows that $\eta\left(  L\right)
=\eta\left(  X\right)  \neq0$, and $\bs\left(  L\right)  =tL$. Thus $t$ is an
eigenvalue for $\bs$.
\end{proof}

\begin{proposition}
\label{Pro8.4}$\operatorname*{PSp}\left(  \bs\right)  \cap
\mathbb{T}=
\vphantom{\tilde{\bs}}
\operatorname*{PSp}\left(  
\vphantom{\bs}
\smash{\tilde{\bs}}\right)  \cap\mathbb{T}$.
\end{proposition}

\begin{proof}
Let $t\in
\vphantom{\tilde{\bs}}
\operatorname*{PSp}\left(  
\vphantom{\bs}
\smash{\tilde{\bs}}\right)  \cap\mathbb{T}$, and
let $X^{\prime}\in\mathcal{M}^{\prime}$ be a corresponding eigenvector:%
\[
\tilde{\bs}\left(  X^{\prime}\right)  =tX^{\prime}.
\]
This means that, with $X=JX^{\prime}J\in\mathcal{M}$,
\[
\sum_{j}V_{j}^{\ast}\Delta^{\frac{1}{2}}X\Delta^{\frac{1}{2}}V_{j}^{{}}%
=\bar{t}\Delta^{\frac{1}{2}}X\Delta^{\frac{1}{2}}.
\]
In the same way as in the proof of Proposition \ref{Pro8.2} one uses this
identity to establish%
\[
\ip{\Delta^{\frac{1}{2}} X ^{\ast}\Phi}{\bs\left( A\right) \Phi}=\bar{t}%
\ip{\Delta^{\frac{1}{2}} X ^{\ast}\Phi}{A\Phi}%
\]
for all $A\in\mathcal{M}$. Thus the linear functional $\eta\left(  A\right)
=\ip{\Delta^{\frac{1}{2}} X ^{\ast}\Phi}{A\Phi}$ satisfies
\[
\eta\circ\bs
=\bar{t}\eta,
\]
and therefore
\[
\bar{t}\in\operatorname*{PSp}\left(  \bs\right)
\cap\mathbb{T}
\]
by Lemma \ref{Lem8.3}. Since $\operatorname*{PSp}\left(
\bs\right)  $ is invariant under complex conjugation, we obtain that
$\vphantom{\tilde{\bs}}
\operatorname*{PSp}\left(  
\vphantom{\bs}
\smash{\tilde{\bs}}\right)  \cap\mathbb{T}\subset
\operatorname*{PSp}\left(  \bs\right)  \cap\mathbb{T}$. As 
$\Tilde{\Tilde{\bs}}=\bs$, the other implication follows.
\end{proof}

\begin{acknowledgements}
This work was done
over several time periods
when the last three
named authors visited the University
of Oslo. 
The generous hospitality from the UO Mathematics Institute is
gratefully acknowledged.
We thank
Professor Mark Fannes for helpful conversations.
Expert typesetting and
author coordination by Brian Treadway in Iowa is
greatly appreciated.
The research 
was
supported by the Norwegian Research Council, the
University of Oslo,
and the
Scandinavia--Japan Sasakawa Foundation.
\end{acknowledgements}

\ifx\undefined\bysame
\newcommand{\bysame}{\leavevmode\hbox to3em{\hrulefill}\,}
\fi

\end{document}